\documentclass[11pt]{amsart}

\usepackage{amssymb}
\usepackage{crossreftools}
\usepackage{enumerate, xspace}
\usepackage[normalem]{ulem}
\usepackage{bbm}       %allows \mathds{E 1}, without sans is less heavy
\usepackage[backgroundcolor=blue!20!white, linecolor=blue!20!white,
textsize=footnotesize]{todonotes}
\usepackage{comment}
\usepackage{verbatim}
\numberwithin{equation}{section}
\usepackage{stmaryrd}
\usepackage{enumitem}
\setlist[description]{style=nextline, font=$\displaystyle$, labelsep=0pt}

\usepackage[colorlinks=true,linkcolor=blue,citecolor=magenta]{hyperref}

\DeclareMathOperator*{\bigtimes}{\vartimes}

\makeindex
\newtheorem{theorem}{Theorem}[section]
\newtheorem{lemma}[theorem]{Lemma}
\newtheorem{proposition}[theorem]{Proposition}
\newtheorem{corollary}[theorem]{Corollary}

\newtheorem{definition}[theorem]{Definition}
\newtheorem{remark}[theorem]{Remark}

\newtheorem{notation}{Notation}[section]

\newcommand{\cal}{\mathcal}

\newcommand{\bs}[1]{{\boldsymbol #1}}

\newcommand\R{{\mathbb R}}

\newcommand\Z{{\mathbb Z}}

\newcommand\ve{\varepsilon}

\newcommand\vf{\varphi}

\newcommand \ga{\gamma}
\newcommand \om{\omega}

\newcommand\fR{\mathfrak R}

\newcommand\fe{\mathfrak e}
\newcommand\fE{{\mathfrak E}}
\newcommand\ff{\mathfrak f}

\renewcommand{\ge}{\geqslant}
\renewcommand{\le}{\leqslant}

\renewcommand{\hat}{\widehat}
\renewcommand{\tilde}{\widetilde}
\renewcommand{\bar}{\overline}
\newcommand{\bol}[1]{\boldsymbol{#1}}

\newcommand{\lf}[1]{\lfloor #1 \rfloor} % Floor function

\newcommand\cA{{\mathcal A}}

\newcommand\cC{{\mathcal C}}
\newcommand\cD{{\mathcal D}}
\newcommand\cE{{\mathcal E}}

\newcommand\cG{{\mathcal G}}

\newcommand\cL{{\mathcal L}}
\newcommand\cM{{\mathcal M}}

\newcommand\cO{{\mathcal O}}

\newcommand\cR{{\mathcal R}}

\newcommand\bA{{\mathbb A}}
\newcommand\bB{{\mathbb B}}
\newcommand\bC{{\mathbb C}}

\newcommand\bE{{\mathbb E}}

\newcommand\bI{{\mathbb I}}
\newcommand\bF{{\mathbb F}}
\newcommand\bG{{\mathbb G}}
\newcommand\bN{{\mathbb N}}

\newcommand\bR{{\mathbb R}}
\newcommand\bT{{\mathbb T}}
\newcommand\bV{{\mathbb V}}

\newcommand\bZ{{\mathbb Z}}

\newcommand\fF{{\mathfrak F}}
\newcommand\fQ{{\mathfrak Q}}
\newcommand\fS{{\mathfrak S}}

\newcommand\Id{{\mathbbm{1}}}
\newcommand\bmp{{\mathbbm{p}}}

\DeclareMathAlphabet{\mymathbb}{U}{BOONDOX-ds}{m}{n}
\newcommand{\zerob}{{\mymathbb{0}}}
\newcommand{\err}{\fR}

\newcommand\Const{C_\#}

\newcommand{\Ene}{\mathcal E}

\newcommand{\mynegspace}{\hspace{-0.09 em}}
\newcommand{\lvvvert}{\left\rvert\mynegspace\left\rvert\mynegspace\left\rvert}
\newcommand{\rvvvert}{\right\rvert\mynegspace\right\rvert\mynegspace\right\rvert}
\newcommand{\trn}[1]{{\lvvvert #1 \rvvvert}}
\begin{document}
\title{Heat equation from a deterministic dynamics}

\author{Giovanni Canestrari}
\address{\hskip-12pt Department of Mathematics, University of Toronto, 
40 St.~George Street, Toronto, ON M5S 2E4, Canada.}
  \email{giovanni.canestrari@utoronto.ca}

\author{Carlangelo Liverani}
\address{\hskip-12pt Università di Roma Tor Vergata\\
  Via della Ricerca Scientifica 1, Roma\\
 \emph{and} University of Maryland\\
 4176 Campus Drive,
College Park, MD 20742-4015 }
    \email{liverani@mat.uniroma2.it}
  
\author{Stefano Olla}
\address{\hskip-12pt Universit\'e Paris-Dauphine, PSL Research University\\
  CNRS, CEREMADE\\ 
  75016 Paris, France \\ \emph{and}
 Institut Universitaire de France\\ \emph{and} GSSI, L'Aquila}
  \email{olla@ceremade.dauphine.fr}

\date{\today}

\begin{abstract}
We derive the heat equation for the thermal energy under diffusive space-time scaling for a purely deterministic microscopic dynamics satisfying Newton's equations perturbed by an external deterministic chaotic force acting like a magnetic field.
\end{abstract}


\maketitle

\vspace*{-1cm}

\tableofcontents

\section{Introduction}
\renewcommand{\thefootnote}{\fnsymbol{footnote}}
\footnotetext{\hskip-12pt This work was supported by the PRIN Grants ``Regular and stochastic behavior in dynamical systems" (PRIN 2017S35EHN), ``Stochastic properties of dynamical systems" (PRIN 2022NTKXCX); by the MIUR Excellence Department Projects MatMod@TOV awarded to the Department of Mathematics, University of Rome Tor Vergata; and by the Institut Universitaire de France. GC and CL acknowledge membership to the GNFM/INDAM. GC and CL acknowledge the hospitality of the Université Paris-Dauphine, where part of this work was done. We thank Cedric Bernardin for key discussions at the very beginning of this work.  Finally, we thank the referees for a very careful reading of the paper and for many valuable suggestions, and Samarpan Das for pointing out that we missed a factor $2$ in the diffusion coefficient.} 
\renewcommand{\thefootnote}{\arabic{footnote}}
An outstanding conceptual issue, {  going back, at least, to } Zeno and Democritus 2400 years ago,
stems from the fact that the world around us behaves like a continuum
(we describe it using partial differential equations), yet it consists of atoms, hence it is discrete in nature
(as conclusively proven by Einstein \cite{Ei05} who showed how Brownian motion may emerge from the microscopic dynamics).

The derivation of macroscopic PDE evolution laws from microscopic ones has finally become
a well-posed problem with the foundation of statistical mechanics by
Boltzmann, Gibbs, and Maxwell.
Thanks to the axiomatization of probability by Kolmogorov,
this problem can now be seen as a purely mathematical one.
The relevance of this line of research in mathematics has been acknowledged by
Hilbert in his famous sixth problem \cite{Hi902}. 

The mathematical problem is outlined as follows:
for dynamics with a large number of degrees of freedom, some conserved quantities (e.g. energy, momentum, density of particles) evolve slowly and, properly rescaling space and time,
the evolution of their distributions is expected to converge to the solution of 
macroscopic hydrodynamic equations (such as Euler or heat equations). 
In contrast, all other observables are expected to exhibit fast fluctuations for generic systems. 
The observed macroscopic equation depends on the space-time scaling chosen.
In dynamics where momentum is conserved, the first natural space-time scale
is the hyperbolic one, and
in this limit arises the Euler system of equations for a compressible gas.
The first \emph{formal} derivation of the Euler equation from the microscopic
Hamiltonian dynamics of interacting particles was done by Morrey \cite{Morrey55}.

Currently, there is no rigorous derivation of the Euler equation
from a completely deterministic Hamiltonian dynamics.
The main problem is the lack of good ergodic properties
of the large microscopic dynamics that guarantee that
energy, momentum, and density of particles are the only conserved quantities that evolve on the macroscopic scale.
Progress has been made only by adding to the dynamics conservative noise terms that destroy all other integrals
of motion, like local random exchanges
of velocities among particles \cite{OVY,FFL94,LO96,FLO97}.

The mathematical problem is even more difficult in diffusive scaling.
This arises in dynamics that do not conserve momentum
or have initial conditions with null momentum.
In this scaling, if {density is fixed and} energy is the only conserved quantity,
the expected macroscopic equation is the heat equation
\begin{equation}
  \label{eq:13}
  \partial_t e = \nabla D(e) \nabla e,
\end{equation}
where $D(e)$ is the thermal diffusivity given by the Green-Kubo formula associated with the space-time correlations of the microscopic energy currents.
For completely deterministic models, even the proof of the existence of the thermal diffusivity $D$ is missing.
In fact, \eqref{eq:13} is not always valid, and in one dimensional systems with
momentum conservation it is expected a superdiffusion of the energy driven 
by a fractional Laplacian, as confirmed by numerical simulation and similar superdiffusion proven with conservative noise
(cf. \cite{Lepri03, Basile16} and references within).

Important progress has been achieved in the so-called kinetic theory of gases through rarefaction density limits (the so-called Grad limit).  
We can mention the derivation of the Boltzmann equation for a gas in the low-density regime
\cite{La75}, the derivation of the hydrodynamic limit from the
Boltzmann equation \cite{Demasi89,Bardos93,Laure09}
and important recent progress concerning the direct derivation of the heat equation and
fluctuations in the kinetic limit
(cf. \cite{BGS-R, Bodineau19,Bodineau20, DHX, DHX2, TGSS} and references within).
However, these results are relevant only for extremely rarefied gases.

For systems at fixed particle densities, some important results have been obtained in diffusive limits assuming the microscopic dynamics to be stochastic.
This has led to a well-developed theory of the hydrodynamic limit for stochastic microscopic dynamics,
starting with the work of Fritz \cite{Fritz88}, Presutti and collaborators \cite{Demasi88} and
Varadhan et al. \cite{GPV,Va}.
See \cite{Sp} for a more complete description of the large-scale limits and the main results and
\cite{KL99} for an exposition of some of the main relevant mathematical techniques.
Still, when interactions are non-linear, the non-equilibrium diffusive behavior of energy
remains an open problem even for stochastic dynamics.
One of the main difficulties is that known methods often rely on entropy
estimates that fail to control quadratic non-linearities (in particular energy).
When stochastic perturbations are present, some results on the existence of thermal diffusivity  $D(e)$ exist \cite{Bernardin11}, and  equilibrium fluctuations of the energy
following linearized heat equations \cite{sasada13} are established.

Diffusive scaling limits have also been established for the densities of {\em independent} particles with very chaotic dynamics, starting with the seminal work of Bunimovich and Sinai on the evolution of the density of a periodic Lorentz gas \cite{BS80}. See also \cite{Bardos97, CD09}.

Even more special are (deterministic)
\emph{extended completely integrable systems} \cite{Spohn23}.
In particular, for a one dimensional system of hard rods, hydrodynamic
limits have been proved in hyperbolic
and diffusive timescales (\cite{Dobru83,Dobru90,Bodri97,Pablo23}).
Since in these integrable extended systems there are infinitely many conserved quantities,
in the diffusing scaling, after recentering on the hyperbolic evolution, a diffusive system of infinitely many equations
arises for the evolution of the energies associated with each conserved quantity.

Another example of an extended, completely integrable system is a chain of harmonic oscillators. When this dynamics is perturbed by an energy conserving noise
(like the random flip of the sign of the velocities), completely integrability is destroyed, and it is possible to prove the diffusive behavior of the energy \eqref{eq:13} with a thermal diffusivity $D$
independent of $e$ and explicitly computable \cite{Bernardin05,bernardin07,klo22}.
Such noise is multiplicative, but the linearity of the dynamics allows to perform the proof by computing the covariances.
Unfortunately, this yields a result only on the asymptotic behavior of
the average density of the energy, and not as a law of large numbers, as no known technique allows the control of higher moments.

To date, no direct hydrodynamic result for the heat equation (energy conservation)
at high-density
is available for purely deterministic (Hamiltonian) interacting systems.
In this paper, we introduce a mechanical model for which
we can perform a diffusive scaling
and show that, on average, the energy density evolves according to the heat Equation  \eqref{eq:13}.

To our knowledge, this is the first proof that a mechanical model yields
the heat equation in some appropriate scaling limit (without rarefaction), thereby showing that the heat equation can indeed be a consequence of deterministic microscopic classical mechanics.

The model is a harmonic chain in which each particle is subjected to a fast evolving external
force acting like a magnetic field. Such a field has a deterministic but chaotic evolution.
In particular, each particle has two degrees of freedom, and
the external force acts by rotating the velocity vector.
In this way, it does not change the energy of the particles.
The basic idea is to consider a system with three time scales:
the macroscopic scale in which we establish the heat equation,
the microscopic one in which the evolution of the positions
and velocities of the particles takes place,
and a faster timescale describing the evolution of the external field (the {\em environment}). 
Hence, the microscopic dynamics can be seen as a fast-slow system. 
This particular mechanical model was suggested to us by Cedric Bernardin.\\

Our key idea is that the fast chaotic external forcing acts on the harmonic chain
similarly to a stochastic noise that randomly rotates
the velocities (similar to the noise in \cite{Bernardin05,bernardin07,klo22}).
This phenomenon is intensively studied in fast-slow dynamical systems
{(under the name of averaging or homogenization for fast-slow systems)},
and some powerful techniques have been developed \cite{Do05, DSL15, KM17, CFKMZ}. 
In fact, we prove that, under proper scaling on the fast external field, in the diffusive limit the
covariance matrix converges to the one of the corresponding stochastic dynamics.
A main obstacle is the need to control, uniformly in the size of the system, the difference between the two time evolutions for a time much longer than that achieved in any existing result.

To obtain our result, we combine \emph{standard pairs}, a technique inspired by the work of
Dolgopyat \cite{Do05} (see \cite{DSL15}
for a simple presentation) and techniques coming from the theory of hydrodynamic limits of stochastic dynamics \cite{klo22}.
Standard pairs are much more flexible than Markov partitions and are well adapted to combining dynamics with probabilistic ideas, see the introduction of Section 4.3 for more comments on this.
In particular, we succeeded in controlling the speed at which the deterministic dynamics converges to the stochastic one
in a precise quantitative manner and for times that increase with the size of the system
(since we take a diffusive scaling). In a different context,
a related point of view has been explored in \cite{DSL16, DSL18, DLPV, Li18, CL22}.
However, in all previous work, the fast-slow systems were always systems
with few degrees of freedom, with the notable exception of \cite{DoLi11}, for which, contrary to the present work,
no estimate on the speed of convergence is available. Here, we succeed in controlling the dependence of the speed of convergence on the system size. This is necessary in order to perform a hydrodynamic limit.

Although we could model the exterior force evolution by any sufficiently chaotic system, to simplify the presentation, we chose the simplest possible model
(a smooth expanding map of the circle).
This has the disadvantage that the dynamics is not reversible,
but it has the big advantage of considerably reducing the technicalities,
hence making the argument accessible also to people with no prior knowledge of dynamical systems. 
Indeed, this allowed us to make the paper essentially self-contained by adding only a couple of appendices; hence catering to a much larger audience.
The alternative to having the evolution of the magnetic-like field modeled by a reversible Anosov map
or an Anosov flow would have required the introduction of many other technical ideas, which would have clouded the main argument without changing the conceptual relevance of the result.

For the same reasons, we use periodic boundary conditions, the simplest possibility. However, our results can be extended to different boundary conditions, such as open systems in contact with a heat bath or periodic forcing, where proofs already exist for the corresponding stochastic dynamics \cite{klo22}. 
The result can also be extended to the unpinned dynamics (that is $\omega_0 = 0$ in \eqref{eq:det}), where more conserved quantities will appear, cf. \cite{kos1} for the corresponding stochastic dynamics.

The effect of stationary magnetic fields acting on harmonic chains
has also been studied in \cite{saito18,cane21,bhat22}. Our situation is very different, as the forcing is time-dependent
but does not exchange energy with the system. Hence, contrary to the constant case, we do not consider the energy of this exterior field as part of the system. Finally, in \cite{giardina05} is proposed and numerically studied a discrete time dynamics where a time-dependent external field exchanges chaotically
the momentum between nearest-neighbor particles conserving the energy, but without
any other interaction in the dynamics. 

%%%%%%%%%%%%%%%%%%%%%%%%%%%
\section{The deterministic model and the result}\label{sec:determinsitic}
For each $N\in\bN$, \(x \in \bZ_N\),\footnote{We assume \textit{periodic boundary conditions}, in the sense that the site \(x=N+1\) coincides with the site \(x=1\).}\; \(\theta = (\theta_x) \in \bT^N\), $\ve\in (0,1)$, we denote by \(q_x^{N,\ve}(t)\), \(p_x^{N,\ve}(t) \in \bR^2\) the solution of the following ODE on $\bR^{2N}\times\bR^{2N}$,
\begin{equation}
    \label{eq:det}
   \begin{split}
    &\dot{q}^{N,\ve}_x = p_x^{N,\ve},\\
    &\dot{ p}_x^{N,\ve} = (\Delta_{ \bZ_N} - \omega_0^2) q_x^{N,\ve} + \ve^{-\frac 12}b \bigl( \theta_x^{\ve}(t)\bigr )Jp_x^{N,\ve},
   \end{split} 
\end{equation}
where\footnote{ To simplify notations, we set $\bT=\bT^1$.} 
\[
\begin{split}
&b \in \cC^{\infty}(\bT,\bR),\quad
     J = \begin{pmatrix}
        0 & 1\\
        -1&0
     \end{pmatrix}, \\
     &\Delta_{\bZ_N} q_x^{N,\ve} = q^{N,\ve}_{x+1}+q^{N,\ve}_{x-1}-2q^{N, \ve}_x  \;;\quad  \omega_0 >0,
\end{split}   
\]
and
\begin{equation}\label{eq:thetadyn}
    \theta_x^{\ve}(t) = f^{\lf{t\ve^{-1}}}(\theta_x),
\end{equation}
where \(f \in \cC^{\infty}(\bT,\bT)\), \(f' \geq \lambda> 1\).\\
Equation \eqref{eq:det} corresponds to the classical evolution of a periodic harmonic chain in which each particle is subject to a fast evolving external magnetic-like field.  In particular, each particle has two degrees of freedom, and its motion is restricted to a plane, while the external ``magnetic" field is perpendicular to the planes of the particles.\\

The choice of an expanding map is motivated by its strong statistical properties. In particular, by the fast decay of correlations.  Moreover, the discrete-time dynamical system \( (f, \bT)\) has a unique absolutely continuous invariant probability measure (called SRB measure); see Appendix \ref{sec:TO}, e.g. Lemma \ref{lem:dec_cor} and Remark \ref{rem:uniqueabs}, for details. That is, calling \(\rho_{\star}\) the invariant density, for each $\vf\in\cC^0(\bT,\bR)$,
\begin{equation}\label{eq:inve_den}
\int_{\bT}\vf( f(\theta))\rho_\star(\theta) d\theta= \int_{\bT} \vf(\theta)\rho_\star(\theta) d\theta.
\end{equation}
For simplicity, we assume that
\begin{equation}
    \label{eq:zero-average}
\int_{\mathbb{T}}b(\theta)\rho_{\star}(\theta)d\theta = 0,    
\end{equation}
i.e, the fluctuations of the external field are zero-average.
\begin{remark}
If $b$ were not zero average, this would be equivalent to the presence of an additional very strong static magnetic field. Hence, the particles would move in very tight spirals. This would not change the essence of the result, as it would not contribute to the energy transport, but would substantially complicate the exposition.
\end{remark}
The following sum of correlations will play an essential role
\begin{equation}
    \label{eq:green-Kubo}
    \gamma = \frac 12\int_{\mathbb{T}}b(\theta)^2\rho_{\star}(\theta)d\theta + \sum_{k=1}^{\infty}\int_{\mathbb{T}}b(\theta)b(f^k(\theta))\rho_{\star}(\theta)d\theta.
\end{equation}
Since \(b \in \cC^\infty(\bT, \bR)\), Lemma \ref{lem:dec_cor} (with $\rho=b\rho_\star$ and $g=b$), and Equation \eqref{eq:zero-average} imply that \(\gamma < \infty\). Notice that $2\gamma$ is the asymptotic variance of the Birkhoff sums of $b$
(see Lemma \ref{lem:L2menogamma}). In particular, $\ve^{-\frac 12}\int_0^t b(\theta_x^\ve(s))ds$ has asymptotic variance $2\gamma t$, consistently with the noise amplitude $\sqrt{2\gamma}$ in \eqref{eq:2st} and the generator in \eqref{eq:5}, see Lemma \ref{lem:covariance}.\\
In addition, if we assume that there exists a periodic orbit $\{x_n\}_{n=0}^p$, \(p \in \bN\), $x_{n+1}=f(x_n)$, $x_p=x_0$, such that
\begin{equation}\label{eq:no_co}
\sum_{n=1}^p b(x_n)\neq 0,
\end{equation}
then Lemma \ref{lem:coboundaries} (see also the subsequent discussion in Appendix \ref{subsec:coboundaries}) implies that \(\gamma >0\). 

\begin{notation}\label{notation:z}In the following we will drop $N,\ve$ from the notation, and we denote $\boldsymbol{p} = (p_x)_{x\in\mathbb{Z}_N}$,
    $\boldsymbol{q} = (q_x)_{x\in\mathbb{Z}_N}$ and
    $z = (\boldsymbol{q}, \boldsymbol{p})$, unless it creates ambiguities.
\end{notation}

We define the single particle energy and  the total energy,
\begin{equation}\label{eq:energy}
    \fe_x(z) = \frac{p_x^2}{2} + \frac{1}{2}(q_x - q_{x-1})^2 + \frac{\omega^2_0 q_x^2}{2}\,; \quad \cE_N(z)= \sum_{x \in \bZ_N}\fe_x(z).
\end{equation}
Here and below, given any vectors $v,w\in\bR^2$, we often
write $v^2$ for $\|v\|^2$ and $vu$ for $\langle v,u\rangle$.\label{page1}\\
Since Equations \eqref{eq:det} are linear, we can rescale $q,p$ so that $ \cE_N=N$. Consequently, given that the total energy is conserved (see Lemma \ref{lem:energy}), the evolution takes place on the energy shell
\begin{equation}\label{eq:shell}
    \Sigma_N = \{ z \in \bR^{4N}: \cE_N(z) = N\}.
\end{equation}
\begin{notation}\label{not:2}
Let $ (z(t,\bar z,\bar\theta),\theta(t,\bar \theta))=\phi^t(\bar z,\bar \theta)$ be defined by \eqref{eq:det} and \eqref{eq:thetadyn} with initial conditions $(z(0,\bar z,\bar\theta),\theta(0,\bar \theta))=(\bar z,\bar \theta)$. If the initial conditions are clear from the context, we may write $z(t),\theta(t)$ to simplify the notation. Also, given $A\in\cC^0(\Sigma_N\times\bT^N,\bR)$ we will write 
\[
A(t)=A(z(t),\theta(t))=A(z(t,\bar z,\bar \theta), \theta(t,\bar \theta))=A\circ \phi^t(\bar z,\bar \theta),
\]
and, given a measure $\nu$ on $\Sigma_N\times\bT^N$, we use the notions
\[
\bE_\nu(A(t))=\nu(A(t))=\int_{\Sigma_N\times\bT^N} A(z(t,\bar z,\bar \theta), \theta(t,\bar \theta))\nu(d\bar z, d\bar\theta).
\]
Finally, $\phi_{\star}^{t}\nu$ is the usual push forward of the measure $\nu$ by the map $\phi^{t}$.
\end{notation}
The last crucial ingredient to specify is the class of initial conditions. As we are looking for statistical properties, we have to start with an initial condition described by a probability measure. This is necessary since the dynamics is deterministic, and if one considers a deterministic initial condition, there would be no stochastic behavior whatsoever. Moreover, we want measures that behave reasonably well in the thermodynamic limit. 
\begin{definition}\label{def:initial0}
  Let $\widetilde \cM_{\textrm{init}} =  \widetilde \cM_{\textrm{init}}(N,C_0)$, $ C_0>0$, $N\in\bN$, be the weak closure of the probability measures $\mu_N$ on $\Sigma_N\times\bT^N$, \(N\in\bN\), such that
\begin{equation}\label{eq:random}
d\mu_N(\bol q, \bol p,\theta)= \bar\mu_{N}(d\bol q,d\bol p) \rho_N(\theta) d\theta,
\end{equation}
where $\rho_N(\theta)=\prod_{x\in\bZ_N}\rho_{x,N}(\theta_x)$, with $\rho_{x,N}\in \mathcal C^1(\bT)$ such that
\begin{equation}\label{eq:Czero}
  \sup_{x\in\bZ_N}\sup_{\theta_x\in\bT}
  \frac{|\partial_{\theta_x} \rho_{x,N}(\theta_x)|}{\rho_{x,N}(\theta_x)}\leq C_0,
\end{equation}
and \(\bar \mu_{N}\) is any probability measure on \(\Sigma_N\). Furthermore, we denote by $\bE_{\mu_N}$ the expectation with respect to $\mu_N$.
\end{definition}

\begin{remark}\label{nota:2}
In Definition \ref{def:initial0} we only demand some regularity in the initial distribution of $\theta$, while we can choose $\bar\mu_{N}(dq,dp)$ as singular as desired, even deterministic. The product structure of the measure is used only to simplify the presentation. In fact, since product measures are not invariant by the dynamics, we will be forced to introduce a much larger class of initial conditions. Unfortunately, its introduction requires several technical definitions. Hence, here we prefer to present a simple, but already quite general, class of measures and postpone to Definition \ref{def:initial} the precise definition of the set $\cM_{\textrm{init}}\supset \widetilde \cM_{\textrm{init}}$ of measures for which our results are valid (the inclusion is proven in Lemma \ref{lem:density}). Notably, we will have \(\phi_{\star}^{t}\cM_{\textrm{init}} \subset \cM_{\textrm{init}}\), for all $t\geq 0$, (Theorem \ref{thm:standard-pair-invariance}).
\end{remark}
Let \(\mathcal P(\bT)\) be the set of probability measures on \(\bT\), endowed with the weak topology.
We are interested in sequences of initial distributions $\mu_N \in \cM_{\textrm{init}}(N)$
corresponding to a  macroscopic energy distribution
\(\Ene_0(du) \in \mathcal P(\bT)\). See Definition \ref{def:initial} for the precise definition of \(\cM_{\textrm{init}}\supset \widetilde \cM_{\textrm{init}}\), and  Remark \ref{rem:C0_choice} for our use of the notation.
\begin{definition}\label{def:initialstar}
The class of initial conditions \(  \cM^{\star}_{\textrm{init}}\) consists of the sequences $(\mu_N)_{N\in\bN}$, $\mu_N\in \cM_{\textrm{init}}(N)$, such that
\begin{equation}
  \label{eq:initial-conditions}
  \begin{split}
   & \lim_{N \rightarrow \infty}\bE_{\mu_N}\left[\frac{1}{N}\sum_{x \in \bZ_N} \varphi\left(\frac{x}{N}\right)\fe_x(z(0))\right] = \int_{\bT}\varphi(u)\Ene_0(du),
\end{split}
\end{equation}
for some $\Ene_0\in\mathcal{P}(\bT)$
and any \(\varphi \in \cC^0(\bT,\bR)\).
\end{definition}
We need to introduce a quantitative version of \eqref{eq:initial-conditions} to state our result. Given \({ (\mu_N)_{N\in\bN}} \in \cM^{\star}_{init}\), for \(k\in \bZ\), \(N \in \bN\), we set
\begin{equation}
\label{eq:varpi}
\begin{split}
      &\varpi \left(k,N,\mu_N\right) = \left|\bE_{\mu_N}\left[\frac{1}{N}\sum_{x\in \bZ_N} e^{2\pi i k \frac{x}{N}}\fe_x(z(0))\right] - \int_{\bT}e^{2\pi i k u}\Ene_0(du)\right|,\\
&    \Psi \left(N; \mu_N\right) = \sum_{k\in \bZ\setminus\{0\}}{k^{-7}}\varpi(N,k,\mu_N).
\end{split}
\end{equation}
We call \(\Psi(N; \mu_N)\) the {\em convergence rate to the energy profile associated with the sequence of measures \((\mu_N)_{N\in\bN} \in \cM^{\star}_{init}\)}. 
\begin{remark}
    If the initial profile of energies has a smooth density, i.e. 
    $\Ene_0(du) = e_0(u) du$, and the initial distribution satisfies
    $\bE_{\mu_N} \left(\fe_x(0) \right) = e_0\left(\frac xN\right)$, then 
    $\Psi(N; \mu_N) \le \frac CN$, with a constant $C$ depending on the regularity of $e_0$.
\end{remark}
More generally, Equation \eqref{eq:initial-conditions}, with the choice \(\vf(x) = e^{2\pi i kx}\), implies that \(\lim_{N \to \infty} \varpi \left(k,N,\mu_N\right) = 0\) and,  consequently,
\begin{equation}\label{eq:PsiN}
\lim_{N \rightarrow \infty}\Psi(N; \mu_N) =  0.
\end{equation}
However, the convergence in \eqref{eq:initial-conditions} could be arbitrarily slow. Since we aim for a quantitative result, we need to quantify the speed of convergence, and hence the utility of \(\Psi\).\\
In order to formulate our main result, we introduce the measure-valued heat equation with initial condition $\Ene_0\in \mathcal P(\bT)$,
\begin{equation}
  \label{eq:PDE}
  \begin{split}
    &\partial_t \Ene = \frac{D}{2\gamma}\partial_u^2 \Ene\\
    &\Ene(0,du) = \Ene_0(du),
  \end{split}
\end{equation}
where \(\gamma\) is given by \eqref{eq:green-Kubo} and \(D = \frac{2}{2+\omega_0^2 +\omega_0\sqrt{\omega_0^2 +4}}\), see Appendix \ref{sec:the-diff-expr}.
\begin{definition}\label{def:init}
We say that a function \(\Ene:[0,+\infty) \rightarrow \mathcal P (\bT)\) is a (weak, measure-valued) solution of \eqref{eq:PDE} if: it belongs to \(\cC^0([0,+\infty), \mathcal P(\bT))\) and for any 
\(\varphi \in \cC^2(\bT,\bR)\), we have
\begin{equation}\label{eq:wheat}
   \int_\bT \varphi(u) \Ene(t,du) - \int_\bT \varphi(u)\Ene_0(du) = \frac{D}{2\gamma}\int_0^t ds \int_\bT \varphi''(u)\Ene(s,du).  
 \end{equation}
\end{definition}
It is well known that for any $t>0$ the measure $\Ene(t,du)$ has a smooth density,  and that the solution of \eqref{eq:wheat} is unique.
Hence, for $t>0$, we have $\Ene(t,du) = \Ene(t,u) du$, where  $\Ene(t,u)$ is a strong solution of \eqref{eq:PDE}  such that the weak limit of $\Ene(t,u) du$ for $t\to 0$ is  $\Ene_0(du)$. Nevertheless, it is convenient to use the formulation \eqref{eq:wheat} for all times since it is the sense in which the convergence in our main Theorem takes place.\\
We can finally state our main result.
\begin{theorem}[Main Theorem]
  \label{thm:heat-equation}
For each sequence $(\ve_N)_{N\in\bN}$, such that $\ve_N\leq N^{-\alpha}$, $\alpha>6$, and for any sequence of initial probability measures  $(\mu_N)_{N\in\bN}\in \cM^{\star}_{init}$ there exist \(C_{\star} , \tilde D\in \bR_{+}\) such that, for each $N\in\bN$ and $\vf\in\cC^7(\bT,\bR)$, we have
\[
\begin{split}
&\sup_{t\ge 0} \left|\bE_{\mu_N}\left[\frac{1}{N}\sum_{x \in \bZ_N} \varphi\left(\frac{x}{N}\right)\fe_x(z(N^2 t))\right] - \int_{\bT}\varphi(u)\Ene(t,du)\right|
 \\
&\phantom{\hskip3.8cm}
\leq C_{\star} \|\vf\|_{\cC^7} \left\{e^{-\tilde D t}\Psi(N; \mu_N)+ 
 N^{-\beta_{\star}} \left(\ln N\right)^2\right\},
\end{split}
\]
where \(\Ene(t,du)\) is the unique solution of \eqref{eq:PDE}, with initial data 
\(\Ene_0\in \mathcal P(\bT)\);
$z(t)$ is the solution of \eqref{eq:det} with initial measures $(\mu_N)_{N\in\bN}$ that satisfy \eqref{eq:initial-conditions} with the same \(\Ene_0\); $\Psi(N; \mu_N)$ is defined in \eqref{eq:varpi} and $\beta_{\star} = \min\left\{1,\frac{\alpha}6-1\right\}$.
\end{theorem}

The rest of the paper is devoted to the proof of Theorem \ref{thm:heat-equation}. Such a proof, after several preliminaries, is concluded in section \ref{sec:atlast}. \\
Before starting, let us state an important corollary of Theorem \ref{thm:heat-equation}.\\
The corollary states that on time scales shorter than $N^{2}$ the energy does not evolve; on the time scale $N^2$ the energy evolves according to the heat equation, and on longer time scales the energy is uniformly distributed.
\begin{corollary}\label{cor:long_time}
For any $ (\mu_N)_{\bN}\in\cM^{\star}_{\textrm{init}}$, $\ve_N\leq N^{-\alpha}$, $\alpha>6$, $t\in\bR_+$, and $\beta>0$ we have
\[
\begin{split}
& \lim_{N\to\infty}\frac{1}{N}\sum_{x \in \bZ_N}\bE_{\mu_N}( \fe_x(z(N^{2-\beta}  t)))\delta_{\frac{x}{N}}
= \Ene_0\\
&\lim_{N\to\infty}\frac{1}{N}\sum_{x \in \bZ_N}\bE_{\mu_N}( \fe_x(z(N^{2} t)))\delta_{\frac{x}{N}}= \Ene(t)\\
&\lim_{N\to\infty}\frac{1}{N}\sum_{x \in \bZ_N}\bE_{\mu_N}( \fe_x(z(N^{2+\beta}  t)))\delta_{\frac{x}{N}}
= \operatorname{Leb},
\end{split}
\]
where the limits are meant in the weak topology.
\end{corollary}
 \begin{proof}
For each $\vf\in\cC^0 (\bT, \bR)$ and $\delta>0$ let 
$\vf_\delta\in\cC^7 (\bT, \bR)$ be such that $\|\vf-\vf_\delta\|_\infty\leq \delta$. By Theorem \ref{thm:heat-equation}, for each $\beta\in\bR$,
 \[
 \begin{split}
& \left|\bE_{\mu_N}\left[\frac{1}{N}\sum_{x \in \bZ_N} \varphi\left(\frac{x}{N}\right)\fe_x(z(N^{2+\beta} t))\right] - 
\int_{\bT}\varphi(u)\Ene(N^{\beta}t,du)\right|\leq 
2\delta \\
&\phantom{|\bE_{\mu_N}\frac{1}{N}\sum_{x \in \bZ_N} \varphi}
 +C_{\star} \left\{\Psi(N; \mu_N)+ N^{-\beta_{\star}} \left(\ln N\right)^2\right\}\|\vf_\delta\|_{\cC^7}.
\end{split}
 \]
Therefore, recalling \eqref{eq:PsiN} and the arbitrariness of $\delta$, we have
\[
\lim_{N\to\infty} \left|\bE_{\mu_N}\left[\frac{1}{N}\sum_{x \in \bZ_N} \varphi\left(\frac{x}{N}\right)\fe_x(z(N^{2+\beta} t))\right] - 
\int_{\bT}\varphi (u)\Ene(N^{\beta}t,du)\right|=0.
\]
If $\beta=0$, the second statement of the Lemma follows. 
If $\beta<0$, then
\begin{equation*}
  \lim_{N\to\infty} \int_{\bT}\varphi(u)\Ene(N^{\beta} t,du)= \int_{\bT}\varphi(u) \Ene_0(du)
\end{equation*}
yielding the first statement.
If $\beta>0$, then the last statement follows, since\footnote{ One can easily check it by studying the behavior of Equation \eqref{eq:PDE} when applied to the Fourier modes $e^{2\pi ikx}$.}
 \[
 \lim_{N\to\infty} \int_{\bT}\varphi(u)\Ene(N^{\beta}t,du)= \int_{\bT}\varphi(u) du.
\]
\vskip-24pt \end{proof}

\bigskip

%%%%%%%%
\subsection{Comments on the results}\label{subsec:comments-results}\ \\
We conclude this section with some comments on the ideas behind the result and its possible generalizations.\\
The microscopic Equation \eqref{eq:det} depends on two parameters $\ve, N$. The thermodynamics limit corresponds to the limit $N\to\infty$, therefore the final objective would be to study the evolution of the energy at time $N^2 t$ in the limit $N\to\infty$ for $\ve$ small, but fixed. This is an extremely hard problem that we do not know how to tackle. On the other extreme, there is the possibility of first taking the limit $\ve\to 0$ at $N$ fixed and then the limit $N\to\infty$. The first limit should yield the stochastic differential equations 
\begin{equation}
  \label{eq:2}
d\tilde q_x = \tilde p_x dt, \qquad  d \tilde p_x =  \left(\Delta - \om_0^2\right) \tilde q_x dt - \gamma \tilde p_x dt
  + \sqrt {2\gamma} J\tilde p_x dw_x(t) ,
\end{equation}
where $\Delta \tilde q_x = \tilde q_{x+1}+\tilde q_{x-1} - 2\tilde q_x$ with periodic boundary, and $\{w_x(t)\}_{x=1}^ N$ are independent
one-dimensional standard Wiener processes.
The SDEs \eqref{eq:2} are written in Ito's formulation, which is the formulation we will exclusively work with in the following. However, it is interesting to note that in the Stratonovitch formulation, equations \eqref{eq:2} read
\begin{equation}
  \label{eq:2st}
 d\tilde q_x = \tilde p_x dt, \qquad  d \tilde p_x =  \left(\Delta - \om_0^2\right) \tilde  q_x dt
  + \sqrt{2\gamma} J \tilde p_x \circ dw_x(t).
\end{equation}
which is more suggestive of the relation with \eqref{eq:det} and may be more familiar to some readers.\\
The study of the hydrodynamic limit for \eqref{eq:2} is similar to existing works in the hydrodynamic limit
(for example, in \cite{klo22} the hydrodynamic limit is proven for a related stochastic dynamics with a different noise term). Such a strategy has been pursued in more general models (e.g. \cite{DoLi11, LO12}). Unfortunately, for such models, the limit $N\to\infty$ still presents obstacles, even the Green-Kubo relation can be obtained only formally \cite{BHLLO}.\\
However, the latter limit is rather non-physical. 
One can interpolate between the above two cases by choosing $\ve=N^{-\alpha}$, this is what we do. Clearly, the smaller $\alpha$ the closer to the ideal goal.\\
The condition $\alpha>6$ in Theorem \ref{thm:heat-equation} is not optimal. Including more terms in the expansion stated in Lemma \ref{lem:Taylor-0} one could probably establish the same result for a smaller $\alpha$, at the cost of considerably more work. However, to have a result in which $\ve$ does not vanish when $N\to \infty$ seems to require substantially new ideas. So, at this point, we did not push for optimality.\\
Our results would also apply to a model where $f,b$ depend on $x$, as long as they all satisfy some obvious uniformity conditions (e.g. \eqref{eq:zero-average} and \eqref{eq:green-Kubo} with \(\gamma_x = \gamma\)). We refrain from such a generalization to simplify the exposition.  In contrast, choosing a specific $f$, e.g., a doubling map, would not simplify the paper in any relevant way.\\
Another valid critique is that the results contained in Theorem \ref{thm:heat-equation} and its corollary concern only the expected value of the energy density field. In order to obtain a law of large numbers (i.e. 
a convergence in probability), we would need some control of higher moments, while we only have bounds on the second moments. This problem affects all existing results of diffusive hydrodynamic limits with energy conservation, even for stochastic dynamics (e.g. \cite{klo22}).

%%%%%%%%%%
\section{Dynamics: basic facts}
Recalling \(z \!= \!(\bol q, \bol p)\), \(\bol q, \bol p \!\in\! \bR^{2N}\),  \(\!\theta \!\in\! \bT^{N}\!\!\), we can write Equation \eqref{eq:det} as
\begin{equation}
    \label{eq:det-matrix}
    \dot{z}(t) = \left (-\mathbb{A}+ \ve^{-\frac 12}\mathbb{B}_{\ve}(t, \theta) \right ) z(t),
\end{equation}
where
\begin{equation}
    \label{eq:unperturbed-matrices}
\begin{split}
    &\mathbb{A} = \begin{pmatrix}
        \zerob & -\mathbbm{1}\\
        -\Delta + \omega^2_{0}\mathbbm{1}& \zerob
    \end{pmatrix}_{4N \times 4N}, \quad \bB_{\ve}(t,\theta) = \begin{pmatrix}
    \zerob & \zerob\\
    \zerob & \bol{B}_{\ve}(t, \theta)
\end{pmatrix}_{4N \times 4N}.
\end{split}
\end{equation}
Here, \(\zerob\) is the null \(2N \times 2N\) matrix, \(\Delta\) and $\bol B_{\ve}(t, \theta)$ are \(2N \times 2N\) dimensional, and \(\Id\) is the identity matrix of the appropriate dimension. More precisely,
\[
\Delta = \begin{pmatrix}
      -2\mathbbm{1} & \mathbbm{1}& 0&0&\cdots & 0 & \mathbbm{1}\\
      \mathbbm{1}&  -2\mathbbm{1} &  \mathbbm{1} &0&\cdots &0&0\\
      0& \mathbbm{1}& -2\mathbbm{1}& \mathbbm{1}&\cdots&0&0\\
      \vdots & \cdots &\cdots & \cdots &\cdots &\cdots&\vdots \\
      \mathbbm{1}&0&0&0&\cdots& \mathbbm{1}& -2\mathbbm{1}
   \end{pmatrix},
\]
\[
\bol B_{\ve}(t, \theta) = \begin{pmatrix}
      b(f^{\lf{t\ve^{-1}}}\theta_1)J & 0 & \cdots &0\\
      0 &  b(f^{\lf{t\ve^{-1}}}\theta_2)J & \cdots&0\\
      \vdots&\vdots&\vdots&\vdots\\
      0&0&\cdots&b(f^{\lf{t\ve^{-1}}}\theta_N)J
   \end{pmatrix}. 
\]
Recall the configuration set $\Sigma_N$ defined by \eqref{eq:shell}, and denote by 
\[
z(t, \bar z, \bar\theta )= 
(\bol q(t,\bar z,\bar\theta), \bol p(t,\bar z,\bar\theta))
\] 
the solution of \eqref{eq:det-matrix} with initial conditions \(\bar z \in \Sigma_N\) and \(\bar\theta \in \bT^N\). 

\begin{lemma}
    \label{lem:energy}
For all \(\bar z \in {\Sigma_N}\), \(\bar\theta\in \bT^N\) 
and \(t \in \bR^{+}\), we have $z(t,\bar z,\bar\theta)\in \Sigma_N$.
\end{lemma}
\begin{proof}
We define the current as
\begin{equation}\label{eq:current}
     j_{x,x+1} =-p_x(q_{x+1}-q_{x}).  
\end{equation}
Recalling the definition of the energy \eqref{eq:energy} and the Equations \eqref{eq:det}, we have
\begin{equation}
  \label{eq:3}
  \frac{d}{dt} \fe_x = j_{x-1,x} - j_{x,x+1}. 
\end{equation}
Note that the above local energy balance does not depend explicitly on the magnetic field. The above implies $ \frac{d}{dt} \sum_{x\in\bZ_N} \fe_x = 0$ and 
consequently \(\mathcal{E}_N(z(t,\bar z,\bar\theta)) = 
\mathcal{E}_N(\bar z) = N\).
\end{proof}
%%%%%%%%%%%%%%%%%%%%%%%%%%%%%
\subsection{Functions controlled by the energy}\ \\
Given the usual norm in $\bR^{4N}$
\[
\|z\| = \sqrt{\langle \bol q, \bol q\rangle + \langle \bol p, \bol p\rangle},
\]
where \(\langle\cdot ,\cdot\rangle\) is the standard scalar product, we have, for $z\in\Sigma_N$,
\begin{equation}\label{eq:enenorm}
    \Const^{-1}\sqrt{N} \le \|z\| \le \Const\sqrt{N}.
\end{equation}
Let \(A \in \cC^{\infty}(\Sigma_N, \bR)\) and define the following seminorms,\footnote{ \(DA=(\partial_{z_i}A)\) and \(D^2 A=(\partial_{z_i}\partial_{z_j}A)\).}
\begin{equation}
\label{eq:functions-controlled-by-energy-norms}
\begin{split}
&\|A\|_{0} = \sup_{z \in \Sigma_{N}}|A(z)|,\\
&\|A\|_{1} =\sup_{z \in \Sigma_{N}}\|D A(z)\|, \\
&\|A\|_{2} =\sup_{z \in \Sigma_N}
\sup_{\bol v \in \bR^{4N}}
\frac{\|D^2 A(z)\bol v\|}{\|\bol v\|}.
\end{split}
\end{equation}

\begin{definition}\label{def:contro-energy}
Let \( \mathfrak{E}_N\) 
be the set of degree-two polynomials on { \(\Sigma_N\)}. For each \(A \in \mathfrak{E}_N\) we define
\begin{equation}
\label{eq:norm-a}
\lvvvert A\rvvvert = \max\left\{\frac{\|A\|_{0}}{N}, 
    \frac{\|A\|_{1}}{\sqrt N},
\| A\|_{2}\right\}.
\end{equation}
\end{definition}
Denote by $S_N = S_N(z)$ the covariance matrix
\begin{equation}
\label{S1ts}
S_N
=\left[
  \begin{array}{cc}
    {S^{(q)}_N}&S^{(q,p)}_N\\
   S^{(p,q)}_N& S^{(p)}_N
  \end{array}
\right]_{2N \times 2N},
\end{equation}
where
\begin{align}
\label{S1ts1}
  &\left[S^{(q)}_N\right]_{x,y}=\langle q_x, q_y\rangle,
    \quad \left[S^{(q,p)}_N\right]_{x,y}=\langle q_x, p_y\rangle,\notag\\
&\\
&
      \left[S^{(p)}_N\right]_{x,y}=\langle p_x, p_y\rangle\quad \mbox{and}
      \quad S^{(p,q)}_N=\Big[S^{(q,p)}_N\Big]^T. \notag
\end{align}
We will need to consider only the following subset of elements of $\fE_N$.\\
Given any \(2N \times 2N\) constant matrix $\bol{\beta}_N$ let
  \begin{equation}
    \label{eq:45}
    A_N(z) = \text{Trace}(\bol{\beta}_N S_N(z)),
  \end{equation}
clearly $A_N\in\fE_N$. In addition,
  $$\|A_N\|_{0} \le \Const \|\bol{\beta}_N\| N,\quad
  \|A_N\|_{1} \le  \Const\|\bol{\beta}_N\| \sqrt N, \quad
  \|A_N\|_{2} \leq  \Const\|\bol{\beta}_N\|.
  $$
  Thus,
  \begin{equation}\label{eq:norm}
        \lvvvert A_N \rvvvert \le \Const \|\bol{\beta}_N\|.
  \end{equation}
 \begin{remark}\label{rem:constant}
By the linearity of \eqref{eq:45}, the bound 
\eqref{eq:norm} also holds for $\bar A_N = \int A_N(z) d \mu(z)$, for any probability measure $\mu$ on $\Sigma_N$.
 \end{remark}
We conclude the subsection with a simple fact.
\begin{lemma}
  \label{lem:matrix-norms}
 For all \(N \in \bN\), \(t \in \bR_+\), \(\theta \in \bT^N\), \( \max\bigl\{\left\|\bA\right\|,\left\|\bB_{\ve}(t, \theta)\right\|\bigr\} \le \Const\).
\end{lemma}
\begin{proof}
Since \((q_{x+1}+q_{x-1}-(\omega_0^2 +2)q_x)^2 \le 3\left(q^2_{x+1} + q^2_{x-1} + (\omega^2_0 +2)^2q^2_{x}\right)\),
  \[
      \begin{split}
      &\|\bA z\|^2 = \sum_{x\in \bZ_N}p^2_{x} + \left(q_{x+1}+q_{x-1}-\left(\omega_0^2 +2\right)q_x\right)^2\le 9(\omega^2_0 +2)^2\|z\|^2.
      \end{split}
  \]
Also, since \(\left( b(f^{t\ve^{-1}}\theta_x)Jp_x\right)^2 \le \|b\|^2_{\infty}p_{x}^2\),
\[
\|\bB_{\ve}(t, \theta)z\|^2 =  \sum_{x\in \bZ_N}(b(f^{\lf{t\ve^{-1}}}\theta_x)Jp_x)^2 \le \|b\|_{\infty}^2\|z\|^2.
\]
\end{proof}
%%%%%%%%%%%%%%%%%%%%
\subsection{Short time dynamics}\ \\
The next Lemma is a boring, but useful, expansion.
\begin{lemma}
    \label{lem:Taylor-0}
  For any \(A \in \mathfrak{E}_N\),  \(\bar z \in \Sigma_N \), \(\bar\theta \in \bT^N\), $h\in  (\ve ,\ve^{\frac 12})$, 
    \begin{equation}\label{eq:duham}
         \begin{split}
                A\left(z(h, \bar z, \bar\theta)\right) &= A(\bar z) - \langle DA(\bar z),  \bA \bar z\rangle h \\
                &+\ve^{-\frac 12}\int_0^h  \left\langle DA (\bar z),  \bB_{\ve}(s, \bar\theta) \bar z\right\rangle ds  \\
                &+ \ve^{-1} \int_{0}^h \int_0^s\langle DA (\bar z), \bB_{\ve}(s,  \bar\theta)\bB_{\ve}(u, \bar\theta) \bar z\rangle dsdu \\
                &+ \ve^{-1}\int_{0}^h\int_0^s \langle D^2 A(\bar z) \hspace{0,05cm}\bB_{\ve}(u, \bar\theta) \bar z,\bB_{\ve}(s, \bar\theta) \bar z\rangle dsdu\\
                &+  \cO\left(\lvvvert A\rvvvert N \left( \ve^{- \frac 32}h^3\right)\right).
            \end{split}
    \end{equation}
\end{lemma}
\begin{proof}
Since the initial condition \( \bar\theta\) will not play any role in this proof, we will omit it.  
By the fundamental theorem of calculus and the chain rule,
    \begin{equation}
        \label{eq:generator}
 \begin{split}
        A(z(h,\bar z)) &= A(\bar z) 
        + \int_{0}^h \!\!\left \langle  D A(z(s,\bar z)), 
        \left(-\bA + \ve^{-\frac 12}\bB_{\ve}(s)\right)z(s,\bar z)\right\rangle ds.
\end{split}
\end{equation}
Using the previous equation for $A_i(z)= z_i$, \(i \in \{1,...,4N\}\), by Lemma \ref{lem:matrix-norms}, we have
\begin{equation}\label{eq:rough_est}
\left\|z(h,\bar z) -\bar z\right\|\leq
\Const h\ve^{-\frac 12}\sqrt N.
\end{equation}
Next, for $A\in \fE_N$, we apply again the chain rule to \( DA(z(s,\bar z))\) and use that \(D^2A\) is constant, obtaining
    \begin{equation}
        \label{eq:loooong}
    \begin{split}
        &A(z(h,\bar z)) = A(\bar z) - \int_{0}^{h}\langle DA (\bar z), \bA z(s,\bar z)\rangle ds \\
        &+ \ve^{-\frac 12} \int_{0}^h\langle DA (\bar z), \bB_{\ve}(s)z(s,\bar z)\rangle ds \\
        &+ \ve^{-1}\int_{0}^h\int_0^s \langle  D^2 A \hspace{0,05cm}\bB_{\ve}(u)z(u,\bar z),\bB_{\ve}(s) z(s,\bar z)\rangle dsdu\\
 &+ \int_{0}^h\int_0^s \langle D^2 A \hspace{0,05cm}\bA z(u,\bar z),\bA z(s,\bar z)\rangle dsdu \\
        &- \ve^{- \frac 12}\int_{0}^h\int_0^s  \langle  D^2 A \hspace{0,05cm} \bA z(u,\bar z),\bB_{\ve}(s) z(s,\bar z) \rangle dsdu \\
        &- \ve^{- \frac 12}\int_{0}^h\int_0^s \langle  D^2 A \hspace{0,05cm}\bB_{\ve}(u) z(u,\bar z), \bA z(s,\bar z) \rangle dsdu.
    \end{split}    
\end{equation}
Let us study each term separately. First, note that
\[
    \begin{split}
 \int_{0}^{h}\langle DA (\bar z),\bA z(s,\bar z)\rangle ds  
=\langle DA (\bar z), \bA \bar z\rangle h + 
\cO(h^2 \ve^{-\frac 12}N \trn{ A})
    \end{split}
\]
where we have used \eqref{eq:rough_est}
and Lemma \ref{lem:matrix-norms}. The second line of \eqref{eq:loooong} is equal to
\[
\begin{split}
    &\ve^{-\frac 12}\int_{0}^h\langle DA (\bar z),  \bB_{\ve}(s)z(s,\bar z)\rangle ds = \ve^{- \frac 12}\int_{0}^h\langle DA (\bar z),  \bB_{\ve}(s)\bar z\rangle ds\\
    &\phantom{=} +\ve^{-\frac 12}\int_{0}^h \int_0^s \left\langle  DA (\bar z), \bB_{\ve}(s)\left(-\bA + \ve^{-\frac 12}\bB_{\ve}(u)\right)z(u,\bar z)\right\rangle dsdu\\
   &= -\ve^{-\frac 12}\int_{0}^h\langle DA (\bar z),  \bB_{\ve}(s)\bA \bar z)\rangle ds +\ve^{-1}\int_{0}^h \int_0^s \langle  DA (\bar z), \bB_{\ve}(s)\bB_{\ve}(u) \bar z\rangle dsdu\\
   &\phantom{=} +\cO\left(\lvvvert A\rvvvert N  h^3\ve^{-\frac 32}
   \right).
  \end{split}  
\]
The third line of \eqref{eq:loooong} is equal to, by \eqref{eq:rough_est},
\[
\begin{split}
     & \ve^{-1}\int_{0}^h\int_0^s \langle D^2 A \hspace{0,05cm}\bB_{\ve}(u)z(u,\bar z),\bB_{\ve}(s) z(s,\bar z)\rangle dsdu\\
& = \ve^{-1}\int_{0}^h\int_0^s \langle D^2 A \hspace{0,05cm}\bB_{\ve}(u)\bar z,\bB_{\ve}(s) \bar z\rangle dsdu+\cO(\|A\|_{2} N h^4 \ve^{-2}),\\
& = \ve^{-1}\int_{0}^h\int_0^s \langle D^2 A \hspace{0,05cm}\bB_{\ve}(u)\bar z,\bB_{\ve}(s) \bar z\rangle dsdu+\cO(\lvvvert A \rvvvert N h^3 \ve^{-\frac{3}{2}}), 
\end{split}
\]
where, in the last line, we have used \(h \in (\ve, \ve^{\frac{1}{2}})\). The remaining terms of \eqref{eq:loooong} can be estimated similarly and yield smaller errors.
\end{proof}
 We conclude the section with a simple consequence.
 \begin{corollary}
    \label{lem:stochastic-Taylor}
 For any \(A \in \mathfrak{E}_N\),  \(\bar z \in \Sigma_N\), $\bar\theta\in\bT^N$, $t\in\bR_+$, $h\in  (\ve,\ve^{\frac 12})$,
     \[
    \int_{t}^{t+h} A(z(u,\bar z,\bar \theta)) du = A(z(t,\bar z, \bar \theta))h + \cO\left(\lvvvert A\rvvvert N h^2 \ve^{-\frac 12}\right)
    \]
\end{corollary}

%%%%%%%%%%%%%%%
\section{Time evolution of averages}\label{sec:averages}
Our goal is to study the evolution of energy averages. To this end, it is useful to have a class of measures that behaves well under the dynamics and allows us to do something equivalent to conditioning for stochastic processes. 

This class of measures has a long history and multiple incarnations, starting at least with the work of Pesin and Sinai \cite{PS83}. However, it has reached maturity in the formulation of Dolgopyat \cite{Do05}, under the name of {\em standard pairs}. In Section \ref{subsec:standardpairs} we introduce the idea of standard pairs and show how to adapt it to the present context.

We will show that the evolution of the average with respect to a standard pair is similar
to a stochastic evolution.
Hence, to state the result, we need first to describe the generator of such a stochastic process.

\subsection{The stochastic generator and the main estimate}\ \\
The stochastic generator associated with the process defined by \eqref{eq:2} is 
\begin{equation}
  \label{eq:5}
  \begin{split}
    \cG &= \mathcal A + \gamma \mathcal S\\
    \mathcal A &=
     \sum_{x\in \bZ_N} \left(
      p_x \partial_{q_x} + \left(\Delta q_x - \omega_0^2 q_x\right)\partial_{p_x}\right)
   \\
    \mathcal S &=  \sum_{x\in \bZ_N}
    \left( p_{x,2}\partial_{p_{x,1}} - p_{x,1}\partial_{ p_{x,2}}\right)^2
    =  \sum_{x\in \bZ_N} \left( J p_x  \partial_{ p_x}\right)^2,
\end{split}
\end{equation}
where $\gamma > 0$ is given by \eqref{eq:green-Kubo}.
We will not directly use the SDEs \eqref{eq:2} or \eqref{eq:2st}, but only the generator \(\cG\).
A useful first estimate on the generator is provided by the following Lemma.

\begin{lemma}
    \label{lem:stochastic-econtrolled}
Let \(A \in \mathfrak{E}_N\), then $\cG A\in\mathfrak{E}_N$ and \(\lvvvert\cG A\rvvvert \leq \Const \lvvvert A\rvvvert \).
 \end{lemma}
\begin{proof}
By Equations \eqref{eq:5} we can write
\begin{equation}\label{eq:gen_shape}
\begin{split}
 \cG A(z) =&  \langle  DA (z), \bA z\rangle  - \gamma\! \sum_{x \in \bZ_N} \left\langle  \frac{\partial A}{\partial p_x},  p_x\right\rangle 
 + \gamma\! \sum_{x\in \bZ_N} \left\langle \frac{\partial^2 A}{\partial p_x ^2 } Jp_{x}, J p_{x}\right\rangle. 
\end{split}
\end{equation}
and the Lemma follows noticing that $\cG A$ is still a degree two polynomial and recalling Definition \ref{def:contro-energy} and Lemma \ref{lem:matrix-norms}.
\end{proof}
The main result of this section is the following theorem, whose proof is postponed to Section \ref{sec:th_main_proof}. Note that the result is stated for the set of initial conditions $\cM_{\textrm{init}}\supset \widetilde \cM_{\textrm{init}}$ that will be defined in Definition \ref{def:initial}.\\
\begin{theorem}\label{thm:basic_fact}
For all \(N \in \bN\), \(\ve \in  (0, N^{-6})\), \(h \in (\ve^{\frac{5}{6}}, 2\, \ve^{\frac 56}) \cap\ve\bN\), $\mu_N\in \cM_{\textrm{init}}$, $t\geq h$, $t\in\ve\bN$, and \(A \in \mathfrak{E}_N\),\footnote{Recall notation \ref{not:2}.}
    \[
   \left|   \mu_{N}\left(A(t+h)\right) -   \mu_{N}\left(A(t)\right)-\int_{t}^{t+h}\mu_{N}\left(\cG A (s)\right)ds\right|
\leq \Const  \err_A(h,\ve, N) h,
\]
 where
\begin{equation}\label{eq:RA}
\err_A(h,\ve, N)= \lvvvert A\rvvvert N \Big(  h^2\ve^{-\frac 32} + \ve   (\ln\ve^{-1} )^2 h^{-1} \Big).
\end{equation}
\end{theorem}
The above theorem quantifies the sense in which the deterministic process is close to the stochastic process \eqref{eq:2} determined by the generator $\cG$. Indeed, if $\err_A$ were zero for all $A$ in the domain of $\cG$, then Dynkin's formula would imply that the process is exactly \eqref{eq:2}.

Before proving Theorem \ref{thm:basic_fact}, we first explain how it will be used and then introduce the basic tool needed in the proof.

\begin{corollary}
\label{cor:accellerated_generator}
Let \(\mu_N\), \(A\) and \(\ve\) be as in Theorem \ref{thm:basic_fact}. For any \(t \in \bR_+\),
\[
\begin{split}
\int_0^t\mu_N(\cG A(N^2s)) ds =&\frac{\mu_N(A(N^2 t))-\mu_N(A(0))}{N^2}\\
&+\cO\left(\lvvvert A\rvvvert \left(t \ve^{\frac 16}(\ln\ve^{-1})^2N
  + \ve^{\frac 13}N^{-1}\right)\right).
\end{split}
\]
\end{corollary}
\begin{proof}
For each \(h \in (\ve^{\frac{5}{6}}, 2\, \ve^{\frac 56}) \cap\ve\bN\), we can write
\begin{equation}
\label{eq:subdivision}
\begin{split}
\int_0^t\mu_N(\cG A(N^2s)) ds=&
\frac{1}{N^{2}}\left(\int_0^{\min\{h,tN^2\}}  \mu_N(\cG A(s)) ds\right.\\
  &\phantom{\frac{1}{N^{2}}(}
  \left.+ \int_{\min\{h,tN^2\}}^{tN^2}\mu_N(\cG A(s))ds\right).
\end{split}
\end{equation}
Using Lemma \ref{lem:stochastic-econtrolled}, we have
\begin{equation}
\label{eq:rough}
\begin{split}
\left|\int_0^{\min\{h,tN^2\}}\frac{\mu_N(\cG A(s)) }{N^2}ds \right| 
&\leq N^{-2}\|\cG A\|_\infty \min \{h,tN^2\}\\
&\le N^{-1}\lvvvert A\rvvvert \min \{h,tN^2\} \le \frac{2\ve^{ \frac 56}}{N} \lvvvert A\rvvvert.
\end{split}
\end{equation}
By \eqref{eq:enenorm} we have $\|z(s,\bar z)\|\leq \Const \sqrt N$. By Equations \eqref{eq:generator} and \eqref{eq:norm-a}, one gets that 
\begin{equation}
\label{eq:very-rough}
\begin{split}
  \left|\frac{\mu_N(A(\min\{h,tN^2\}))-\mu_N(A(0))}{N^2}\right|& \le \Const
  h \ve^{-\frac 12}\frac{\lvvvert A\rvvvert}{N} \le \Const \ve^{\frac 13}\frac{\lvvvert A\rvvvert}{N}.
\end{split}
\end{equation}
If \(t\le \ve^{\frac 56}N^{-2} < h N^{-2}\), Equations \eqref{eq:rough} and \eqref{eq:very-rough} prove the statement; the same argument applies whenever $N^2t<h$.\\
 If \(t > \ve^{\frac 56}N^{-2}\), then we can choose 
 \(h \in (\ve^{\frac 56}, 2\ve^{\frac 56})\cap\ve\bN\) and set $K=\lfloor N^2 t h^{-1}\rfloor$. The integral of $N^{-2}\mu_N(\cG A(s))$ over $[Kh,N^2t]$, an interval of length smaller than $h$, and the difference $N^{-2}\left[\mu_N(A(N^2t))-\mu_N(A(Kh))\right]$ are estimated exactly as in \eqref{eq:rough} and \eqref{eq:very-rough}, hence they are of order $\ve^{\frac 13}N^{-1}\trn{A}$. Accordingly, we can assume without loss of generality that $N^2t=Kh$.
 Then, by Theorem \ref{thm:basic_fact}, the second term of \eqref{eq:subdivision} is equal to
\begin{equation}
\label{eq:iteration}
\begin{split}
&\frac{1}{N^2}\sum_{k=1}^{K-1}\int_{kh}^{(k+1)h}\mu_N(\cG A(s)) ds =\\
&=\sum_{k=1}^{K-1}\frac{\left[\mu_N(A((k+1)h))-\mu_N(A(kh))
  +h \cO(\err_A(h,\ve, N))\right]}{N^2}\\
&=\frac{\mu_N(A(N^2 t))-\mu_N(A(h))}{N^2}
+\cO\left(\frac{h (K-1)}{N^2}\err_A(h,\ve, N)\right).
\end{split}
\end{equation}
Equation \eqref{eq:RA} implies
\begin{equation}\label{eq:errortot}
\begin{split}
  \frac{h (K-1)}{N^2}\err_A(h,\ve,N)&\leq
  t \lvvvert A\rvvvert N\left(h^2 \ve^{-3/2} + \ve (\ln\ve^{-1})^2 h^{-1}\right) \\
 & \leq t   \lvvvert A\rvvvert N \ve^{\frac 16}(\ln\ve^{-1})^2 .
\end{split}
\end{equation}

Therefore, by \eqref{eq:iteration} and \eqref{eq:rough},
\begin{equation}
\label{eq:last-to-conclusion}
\begin{split}
&\int_0^t\mu_N(\cG A(N^2s)) ds =N^{-2}\int_0^{N^2t}\mu_N(\cG A(s)) ds\\
&=\frac{\mu_N(A(N^2 t))-\mu_N(A(h))}{N^2}\\
&\phantom{==}
+ \cO\left(\lvvvert A\rvvvert \left(t\ve^{\frac 16}(\ln \ve^{-1})^2 N
  + N^{-1}\ve^{\frac 5 6}\right)\right),
\end{split}
\end{equation}
concluding the proof.
\end{proof}

\subsection{On the strategy for the proof of Theorem \ref{thm:basic_fact} }\ \\
\label{sec:the-strategy-prove}
Suppose that we start, at $t=0$, with an initial distribution $\mu_N\in \widetilde \cM_{\textrm{init}}$ of the type
\eqref{eq:random} but with a deterministic configuration
$\bar z = (\bar {\bf q}, \bar{\bf p}) \in \Sigma_N$,
i.e. $\mu_{N}(d {\bf q}, d{\bf p}) = \delta_{\bar z}$.
Thanks to the regularity of $\rho_N(\theta)$, we can prove Theorem \ref{thm:basic_fact} even if the support of $\rho_N(\theta)$ is very small, let's say
of the order of $\ve^{16}$. This happens because the time interval $(0,h)$
is much longer than $\ve$ (since we choose $h\sim \ve^{5/6}$),
so that there is plenty of time for the fast dynamics to converge to its invariant measure and produce the correlations appearing in
\eqref{eq:duham} (this will happen in a time of order $\ve\ln\ve^{-1}$, since $f$ is a smooth expanding map, see Lemma \ref{lem:dec_cor}). On the other hand $h<\!\!<\ve^{1/2}$, so that the configuration
of $({\bf q}, {\bf p})$ at time $h$ changes very little from the initial configuration.

To extend the argument at times of order $N^2$, we have to iterate the
estimate for $N^2 h^{-1}$ time intervals of length $h$, as we have done in the proof of
Corollary \ref{cor:accellerated_generator}, see \eqref{eq:iteration}.
The difficulty is that at the time $h$ the configuration of the positions and velocities is a function of the initial $\theta$:
\begin{equation}\label{eq:st_primitive}
  z(h) = G(\bar z,\theta) \in \Sigma_N.
\end{equation}
Although $G(\bar z,\theta)$ is very close to $\bar z$, to iterate, we would need the distribution of $\theta$ conditioned to a particular configuration in $z(h)$ to have the density with regularity as the original $\rho_N$. This is false, in general. So, to iterate, we cannot simply condition on the $({\bf q}, {\bf p})$ configuration.

Accordingly, we need to find a set of probability measures $\mu$ on
$\Sigma_N\times \bT^N$ that have some regularity properties (so we can prove Theorem \ref{thm:basic_fact}),
and are invariant for the dynamics.
In addition, it seems natural that such a set of probability measures should contain measures obtained by iterating a measure in which the $({\bf q}, {\bf p})$ are fixed, as above. That is, measures supported on graphs of the form $G(\bar z,\theta)$, as in \eqref{eq:st_primitive}.
In the next Section, we construct exactly such a class of measures.

%%%%%%%%%%%%%
\subsection{Standard pairs and their evolution}\label{subsec:standardpairs}\ \\
The main problem in dealing with deterministic systems stems from the fact that the only randomness comes from the initial conditions. As already explained, conditioning on the $({\bf q}, {\bf p})$ configuration may completely kill the randomness of the initial conditions.
\\
However, conditioning is a basic tool for dealing with stochastic processes. Thus, it is necessary to introduce some form of conditioning that does not kill the randomness of the initial conditions.\\
This can be done in various ways; one possibility is to introduce Markov partitions. This has the disadvantage of obliterating the differential structure of the dynamics since the conjugation with the associated symbolic dynamics is, typically, only Hölder. In addition, Markov partitions are a very rigid structure that is often not simple to use or generalize.\\
The latter problems can be avoided by employing {\em standard pairs}. The supports of standard pairs are, roughly speaking,
 the finest class of events on which it is possible to condition the evolution without killing the randomness present in the initial conditions. 
They do not form a $\sigma$-algebra, but can be combined to create {\em standard families}. To each standard family is associated a measure which (in some loose sense) takes the place of conditional expectation. In particular, standard families retain a Markov-like structure, thanks to their invariant properties with respect to the dynamics, that makes it possible to implement iterative arguments. Finally, they can be used in very general systems for which Markov partitions are not available.
%%%%%%%%%%%%%%%%%%%%%%%
\subsubsection{Standard pairs: definitions}\ \\ \label{eqref:sp_defin}
Our first goal is to define measures as close as possible to a $\delta$ function for which the deterministic dynamics still exhibit good statistical properties.\\
Choose $\delta_0\in\bR_+$ so that the function $f$, see \eqref{eq:thetadyn}, is invertible on each interval of size $\delta_0$ and set\footnote{ The choice of $16$ is arbitrary, any power large enough would do.} 
\begin{equation}\label{eq:delta_choice}
\delta_{\ve} =\delta_0 \ve^{16}.
\end{equation}
For each number \(a_x, b_x \in \bT\),  \(b_x - a_x \in [\delta_{\ve}/2, \delta_{\ve}]\), \(x\in\bZ_N\), 
we define { the domain} $\bI(a,b)=\bigtimes_{x \in \bZ_N}[a_x,b_x]$, where \(\bigtimes\) stands for the Cartesian product and \(a=(a_x)_{x\in \bZ_N}\) and \(b = (b_x)_{x\in \bZ_N}\). We consider the surfaces defined by the graph of functions \(G \in \cC^2 \left(\bI(a,b), \bR^{4N}\right)\) with \(\text{Im}(G) \subseteq \Sigma_N\), whose components are denoted by
\begin{equation*}
  G(\theta) = \left \{ G_x(\theta) =
    \bigl(G_{q^1_x}(\theta),G_{q^2_x}(\theta),G_{p^1_x}(\theta),G_{p^2_x}(\theta)\bigr) \right\}_{x \in \bZ_N},
\end{equation*}
and set,
\[
    \begin{split}
    &\left\|\partial_\theta G\right\|_{\infty} = \sup_{x \in \bZ_N}\sup_{\theta_x \in [a_x,b_x]}\left(\sum_{y=1}^{N}\left(\frac{\partial G_y(\theta)}{\partial \theta_x}\right)^2\right)^{\frac 12}, \\
    & \left(\frac{\partial G_y(\theta)}{\partial \theta_x}\right)^2 = \sum_{j=1,2}\left(\frac{\partial G_{q_y^j}(\theta)}{\partial \theta_x}\right)^2 + \left(\frac{\partial G_{p_y^j}(\theta)}{\partial \theta_x}\right)^2.
    \end{split}
\]
We can now define the set of graphs that we will use:
for \(C_0 \in \bR_+\), \(\ve\in (0,1)\), and \(N \in \bN\), consider the set of functions 
\begin{equation}\label{eq:Graphdef}
\begin{split}
\mathfrak S(N, \ve, C_0) = \Bigg\{G \in \cC^2\left(\bI(a,b), \Sigma_N\right)\;\bigg | \;& b_x - a_x \in [\delta_{\ve}/2, \delta_{\ve}];\; \\
&\|\partial_\theta G\|_{\infty} \le C_0\sqrt{ N\ve} \Bigg\}.
    \end{split}
\end{equation}
To construct a class of measures supported on such surfaces, we need to define a set of densities: we define the standard probability densities as
\begin{equation}\label{eq:density_def}
    \begin{split}
    \cD_{C_0}(a,b) = \biggl\{\rho \in \cC^1&\left( \bI(a,b), \bR_+\right)\;\bigg |\; \rho(\theta) = \prod_{x \in \bZ_N} \rho_x(\theta_x), \\
& \int_{a_x}^{b_x} \!\!\!\rho_x(\theta_x)d\theta_x = 1,\;\biggl\|\frac{\rho'_x}{\rho_x}\biggr\|_{\cC^0}\! \le C_0\biggr\}.
    \end{split}
\end{equation}
\begin{remark} The choice to restrict $\rho$ to be a product is not necessary, we make it only to simplify some of the following arguments.
\end{remark}
We will call the graphs in \eqref{eq:Graphdef} the {\em support of the standard pair}. Since the precise values of \(a\) and \(b\) do not play any role in the estimates, we will write $\cD_{C_0}$ instead of $\cD_{C_0}(a,b)$, if there is no ambiguity. Similarly, we may write $\mathfrak S$ instead of $\mathfrak S(N,\ve, C_0)$.
\begin{definition}\label{defin:standardp}
A standard pair \(\ell\) is given by \(\ell = (\bI_\ell,G_{\ell}, \rho_{\ell})\), where \(G_{\ell} \in \mathfrak S(N,\ve, C_0)\) and \(\rho_{\ell} \in \cD_{C_0}\) are defined on the same domain $\bI_\ell=\bI(a_\ell,b_\ell)$.\\
A standard pair \(\ell\) induces the probability measure \(\mu_{\ell}\) on \(\Sigma_N\times\bT^N\) by
\[
    \mu_{\ell}(g) = \int_{\bI_\ell}g(G_{\ell}(\theta), \theta)\rho_{\ell}(\theta)d\theta,
\]
for any continuous function \(g\in\cC^0(\Sigma_N\times \bT^N,\bR)\). Here, \(d\theta\) is the Lebesgue measure on \(\bT^N\).
\end{definition}
The measure $\mu_\ell$ is supported on the graph of $G_\ell$, see Figure 1.\\
Since the diameter of $\bI_\ell$ is bounded by $\delta_0N^{\frac 12}\ve^{16}$, the condition \(\|\partial_\theta G\|_{\infty} \le C_0\sqrt{N \ve}\) implies that if \(N \ve^{16} \ll 1\), then the surfaces are almost flat. In other words, the measures $\mu_\ell$ almost fix the positions of the particles but do not precisely fix the state of the external field.\footnote{ If $\partial_\theta G=0$, then the positions and velocities of the particles are completely fixed.}

Let \(\mathfrak S^*(N, \ve, C_0)\) be the set of standard pairs (\(\mathfrak S^* \) if no confusion can arise). \\

\vskip 3pt
\hskip1cm
\begin{tikzpicture}
\draw (-4,0)--(4,0);
\draw (-3.5,-0.5)--(-3.5,4);
\draw[ultra thick] (-1,0)--(1,0);
\draw [ cyan, very thick] plot [smooth] coordinates {(-1,2)  (-.5,2.05) (.5, 1.8) (1, 1.9)};
\draw[->]  (2,3)--(0,2.1);
%%%leggenda
\node at (-4, 3) {$\bR^{4N}$};
\node at (3.5, -.5) {$\bT^N$};
\node at (0, -.5) {$\bI(a,b)=\bI_\ell$};
\node at (3.5,3) {support of $\mu_\ell$};
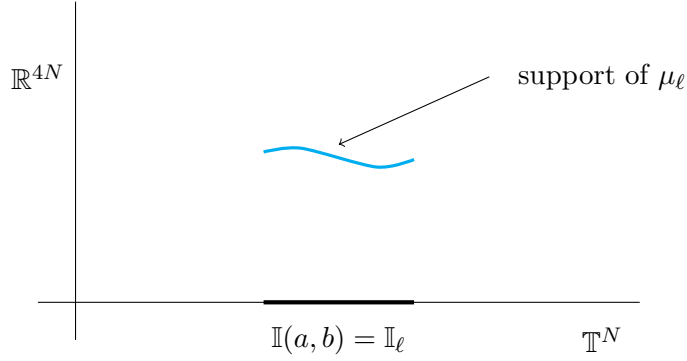
\node at (0.5,-1.5) [align = left] {Figure 1: The graph of $G_\ell$ (in cyan) supporting the measure $\mu_\ell$.}; 
\end{tikzpicture}
\vskip12pt
Next, we define ``convex combination" of standard pairs.
\begin{definition}\label{defin:standardf}
A {\em standard family} is a finite collection of standard pairs \(\{ \ell_j\}_{j=1}^K\subset \mathfrak S^*(N,\ve, C_0) \), together with a probability measure \(\bmp = (p_1, ..., p_K)\), on $\{1,\dots,K\}$. Any standard family \(\ff=\{(\ell_j, p_j)\}_{j=1}^K\) naturally induces a probability measure \(\mu\) on \({\Sigma_N}\times \bT^N\) through the convex combination
\begin{equation}\label{eq:sf_rep}
        \mu_\ff(g) = \sum_{j=1}^Kp_j\mu_{\ell_j}(g), \quad g \in \cC^0({\Sigma_N}\times \bT^N, \bR).
\end{equation}
\end{definition}
We call $\bF_{N,\ve, C_0}$ the set of standard families and $\fF_{N,\ve, C_0}$ the set of probability measures determined by a standard family using the formula \eqref{eq:sf_rep}. Note that different standard families may define the same measure.

%%%%%%%%%%%%%%%%%%%%%%%%%
\subsubsection{Standard pairs: properties}\ \\
The first relevant fact concerning standard families is that they are general enough to describe a wide range of measures. 
\begin{definition} \label{def:initial}
Let $\cM_{\textrm{init}}(N,\ve,C_0)$ be the weak closure of the set \(\mathfrak F_{N,\ve, C_0}\).
\end{definition}
Recall that $\cM^{\star}_{\textrm{init}} (C_0)$ consists of the sequences $(\mu_N)_{N\in\bN}$, $\mu_N\in \cM_{\textrm{init}}(N,C_0)$, satisfying Equation \eqref{eq:initial-conditions} in Definition \ref{def:initialstar}.
\begin{lemma}\label{lem:density}
For any \(\ve\ge 0\) and $N\in\bN$ the set $\widetilde \cM_{\textrm{init}}(N,C_0)$ (see \ref{def:initial0} for the definition) is contained in $\cM_{\textrm{init}}(N,\ve,C_0)$.
\end{lemma}
The proof can be found in Appendix \ref{sec:Standard-Pairs}.\\
Recall that \(\phi^t\) is the flow defined by  \eqref{eq:det}, \eqref{eq:thetadyn}  and the push-forward of measures is \(\phi^t_{\star}\).
In Appendix \ref{sec:Standard-Pairs} we prove the following fundamental fact that justifies our choice of $\cM_{\textrm{init}}(N,\ve,C_0)$.
\begin{theorem}
    \label{thm:standard-pair-invariance}
There exist \(C_*>0\) large enough and \(\ve_0>0\) small enough, depending only on $\omega_0, b$, and $f$, such that, for all $C_0\geq C_*$, \(\ve \in (0, \ve_0)\), $t\in\ve\bN$ and \(N\in\bN\) we have $\phi^t_\star\cM_{\textrm{init}}(N,\ve,C_0)\subset \cM_{\textrm{init}}(N,\ve,C_0)$.
\end{theorem}
\begin{remark}\label{rem:C0_choice}
Let $\bar C_*$ be the inf of $C_*$ for which Theorem \ref{thm:standard-pair-invariance} applies. From now on, we fix $C_0=2\bar C_*$; hence, we will remove it from the notation. Since $\ve$ and $N$ are often clear from the context, we may simply write $\cM_{\textrm{init}}$ or $\cM_{\textrm{init}}(N)$ for $\cM_{\textrm{init}}(N,\ve,C_0)$.
\end{remark}
\begin{remark}\label{rem:grid}
Since $\theta$ jumps only at times belonging to $\ve\bN$, see \eqref{eq:thetadyn}, the family $\{\phi^t\}_{t\ge 0}$ is not a flow: $\phi^{t+s}=\phi^s\circ\phi^t$ holds if $t\in\ve\bN$, but not in general. In order to recover the flow property, one can add to the phase space a clock $\tau\in[0,\ve)$, evolving according to $\dot\tau=1$ with the reset $(z,\theta,\ve)\sim(z,f(\theta),0)$: the resulting dynamics on $\Sigma_N\times\bT^N\times[0,\ve)$ is a flow, and $\phi^t$ is its restriction to the fibre $\{\tau=0\}$. However, to keep the notation lighter, we decided to state some results only for $t\in\ve\bN$ (as for example Theorem \ref{thm:standard-pair-invariance}). This entails no loss of generality for the quantities studied in Sections \ref{sec:hydrodynamic} and \ref{sec:final}. Indeed, by \eqref{eq:3}, \eqref{eq:current} and Lemma \ref{lem:energy},
\[
\sum_{x\in\bZ_N}\left|\frac{d}{dt}\fe_x\right|\leq 2\sum_{x\in\bZ_N}|p_x||q_{x+1}-q_x|\leq 2\cE_N=2N.
\]
Given $s\in\bR_+$, let $s_\ve=N^{-2}\ve\lceil N^2s\ve^{-1}\rceil$, so that $N^2s_\ve\in\ve\bN$ and $0\le N^2 s_\ve- N^2 s\le \ve$. Hence, anticipating the definition \eqref{eq:6} of \(\xi_N\),
\[
\begin{split}
|\xi_N(s_\ve,\vf)-\xi_N(s,\vf)|
&=\Big|\frac1N\sum_{x\in\bZ_N}\vf\big(\tfrac xN\big)\,\bE_{\mu_N}\Big[\int_{N^2s}^{N^2s_\ve}\frac{d}{du}\fe_x(u)\,du\Big]\Big|\\
&\le \frac{\|\vf\|_{\cC^0}}{N}\int_{N^2s}^{N^2s_\ve}\bE_{\mu_N}\Big[\sum_{x\in\bZ_N}\Big|\frac{d}{du}\fe_x(u)\Big|\Big]du
\le2\ve\|\vf\|_{\cC^0}.
\end{split}
\]
Moreover, the time integrals appearing in Section \ref{sec:hydrodynamic} change by the boundary term $|s-s_{\ve}|\|\xi_N(\cdot, \vf'')\|_{\infty} = \cO(\ve N^{-2}\|\vf\|_{\cC^2})$ when $s$ is replaced by $s_\ve$. These errors are both negligible with respect to $R(N,\ve,\cdot)$, see \eqref{eq:err}.
\end{remark}

%%%%%%%
\subsection{Proof of Theorem \ref{thm:basic_fact}}\ \\
\label{sec:th_main_proof}
It suffices to prove Theorem \ref{thm:basic_fact} for an initial condition consisting of the measure associated to a generic standard pair $\ell\in \mathfrak S^*$.  Indeed, by Definition \ref{def:initial} of \(\cM_{\textrm{init}}\), for each $N\in\bN$, $\mu_N\in \cM_{\textrm{init}}$, and any \(\ve > 0\), there exists a sequence $(\ff_n)_{n\in\bN}$, with elements in $\mathfrak F_{N,\ve}$, of standard families such that the associated measures $\mu_{\ff_n}$ converge weakly to $\mu_N$ as \(n \to \infty\). Accordingly, assuming that the theorem is true when $\mu_N=\mu_\ell$, for an arbitrary standard pair $\ell$, we have\footnote{ See Notation \ref{not:2} for the meaning of $A(t)$ and $z(t)$.}
\[
\begin{split}
&\mu_N(A(t+h))=\lim_{n\to\infty} \sum_{(\ell, p_\ell)\in\ff_n} p_\ell \mu_\ell(A(t + h))\\
& =\lim_{n\to\infty} \sum_{(\ell, p_\ell)\in\ff_n} p_\ell  \left\{\mu_{\ell}\left(A(t)\right)+\int_{t}^{t+h}\mu_{\ell}\left(\cG A (s)\right)ds
      +\cO(\err_A(h,\ve, N) )h\right\}\\
&= \mu_{N}\left(A(t)\right)+\int_{t}^{t+h}\mu_{N}\left(\cG A (s)\right)ds
      +\cO(\err_A(h,\ve, N)) h,  
\end{split}
\]
where $\err_A(h,\ve, N)$ is defined by \eqref{eq:RA}. Since \(\cG A(s)\) is bounded, we can exchange the limit and the integral on the second line of the above equation.\\
It remains to prove the statement for the case in which $\mu_N=\mu_\ell$ for a standard pair $\ell\in \fS^*$. By Theorem \ref{thm:standard-pair-invariance}, for all $s\in\ve\bN$, the push-forward \(\phi_{\star}^s \mu_{\ell}\) is equal to \(\mu_{\ff_s}\) for some standard family \(\ff_s \in \bF_{N,\ve}\). We can then use $\mu_{\ff_{t-h}}$ as an initial condition and write, for any \(t \ge h\),
\[
\begin{split}
     &\mu_{\ell}\left(A(z(t+h))\right) =\mu_{\ff_{t-h}}\left(A(z(2h))\right)= \sum_{(\ell', p_{\ell'}) \in \ff_{t-h}}p_{\ell'}\mu_{\ell'}\left(A(z(2h))\right).
 \end{split}
 \]
If the statement is true for \(t =h\), then the equation above is equal to
 \[    
\begin{split}
     &\sum_{(\ell', p_{\ell'}) \in \ff_{t-h}}p_{\ell'}\left( \mu_{\ell'}(A(z(h)))+\int_{h}^{2h}\mu_{\ell'}\left(\cG A(z(s))\right)ds + \cO(\err_A(h,\ve, N ))h \right)\\
&= \mu_{\ell}(A(z(t)))+\int_{h}^{2h}\mu_{\ell}\left(\cG A(z(t-h+s))\right)ds + \cO(\err_A(h,\ve, N ))h  \\
&=\mu_{\ell}(A(z(t)))+\int_{t}^{t+h}\mu_{\ell}\left(\cG A(z(s))\right)ds + \cO(\err_A(h,\ve, N ))h,
\end{split}
\]
which yields the statement for all $t\geq h$. We are then left with the study of \(\mu_{\ell}\left(A(z(2h))\right)\) which, by Theorem \ref{thm:standard-pair-invariance} again, can be written as
   \[
    \begin{split}
    \mu_{\ell}\left(A(z(2h))\right) &=\sum_{(\ell', p_{\ell'}) \in \mathfrak{f}_{h}}p_{\ell '}\mu_{\ell'}\left(A(z(h)\right)\\
    &= \sum_{(\ell', p_{\ell'}) \in \mathfrak{f}_{h}}p_{\ell '} \int_{\bI_{\ell '}}A\left(z\left(h, G_{\ell '}(\theta),\theta\right)\right)\rho_{\ell '}(\theta)d\theta.
    \end{split}
   \]

\begin{remark} The integrals appearing in the above equation are nothing else than the deterministic evolution starting with initial conditions distributed according to the standard pair \(\ell'\). The problem thus reduces to comparing the deterministic evolution with the stochastic one over a time interval $h$, starting from the measure associated with a standard pair.
\end{remark}

By Lemma \ref{lem:Taylor-0}, with \(z(0) =G_{\ell '}(\theta)\) , we can write the above as
   \begin{equation}
    \label{eq:support}
    \begin{split}
        \sum_{\ell' \in \mathfrak{f}_{h}}&p_{\ell '} \int_{\bI_{\ell '}}\rho_{\ell'}(\theta)\Biggl(A(G_{\ell'}(\theta)) - \langle DA(G_{\ell'}(\theta)),\bA G_{\ell'}(\theta)\rangle h \\
        &+ \ve^{-\frac 12}\int_0^h\langle DA (G_{\ell'}(\theta)),  \bB_{\ve}(s, \theta) G_{\ell'}(\theta)\rangle ds  \\
        &+ \ve^{-1} \int_{0}^h \int_0^s  \langle DA (G_{\ell'}(\theta)),\bB_{\ve}(s, \theta)\bB_{\ve}(u, \theta)G_{\ell'}(\theta) \rangle dsdu \\
        &+  \ve^{-1}\int_{0}^h\int_0^s \langle D^2\!A\phantom{!}\bB_{\ve}(u, \theta)G_{\ell'}(\theta),\bB_{\ve}(s, \theta) G_{l'}(\theta)\rangle dsdu\Biggr)d\theta\\
      &+\cO\left( \ve^{- \frac 32}h^3\trn{A} N\right).
    \end{split}
\end{equation}
We want to recognize in \eqref{eq:support} the action of the generator; see \eqref{eq:gen_shape}. To this end, we analyze the terms one by one. The second term in the first line is part of the generator.  The second line of \eqref{eq:support} turns out to be negligible. Indeed, since \(  \sum_{\ell' \in \mathfrak f_h}p_{\ell'}\mu_{\ell'} = \phi^h_{\star}\mu_{\ell}\), Lemma \ref{lem:annoying-term} in Section \ref{sec:aveft}, imply
\begin{equation}
\label{eq:secon_support}
\begin{split}
    &\ve^{-\frac 12}\sum_{\ell' \in \mathfrak f_h}p_{\ell'}\int_{\bI_{\ell'}}\int_0^{h}\left\langle DA (G_{\ell'}(\theta)),  \bB_{\ve}(s, \theta) G_{\ell'}(\theta)\right\rangle\rho_{\ell'}(\theta)ds d\theta \\
&=\ve^{-\frac12}\mu_{\ell}\left(\int_h^{2h}\left \langle DA\left(z(h)\right),  \bB_{\ve}(s) z(h)\right\rangle ds\right)\\
&= \cO( \trn{A} N (\ve (\ln \ve^{-1})^2+ h^{3}\ve^{-\frac 32}).
\end{split}
\end{equation}
\subsubsection{Ito drift}\ \\
The third line of \eqref{eq:support} yields the second term of \eqref{eq:gen_shape} plus an error. To see this, it is convenient to set, for each \(x,y\in \bZ_N\), \( \mathfrak m_x  (\bol q, \bol p)= \left\langle  \frac{\partial A}{\partial p_{x}}, p_{x}\right\rangle\) and \(\mathfrak l_{x,y} (\bol q,\bol p)  = \left\langle \frac{\partial^2 A}{\partial p_{x}\partial p_{y}}Jp_{y}, J p_{x}\right\rangle\).
Observe that \(\mathfrak m_x, \mathfrak l_{x,y} \in \fE_N\), and\footnote{
Setting \(\beta_{x,y} =  \text{sgn}(\mathfrak l_{x,y}) (Jp_y) \otimes (Jp_x)\) and \(\alpha_{x,y} = \frac{\partial^2 A}{\partial p_x \partial p_y} \), we have 
\[
\begin{split}
\sum_{x\in \bZ_N} \sum_{y\in \bZ_N}|\mathfrak l_{x,y}| &=  \text{Tr} ( \alpha \beta)\leq \sqrt{\text{Tr} (\beta \beta^T)\text{Tr} (\alpha \alpha^T)}\leq \|p\|^2 \sqrt{\sum_{x,y} \left(\frac{\partial^2 A}{\partial p_x \partial p_y}\right)^2}\\
&\leq 2N^{\frac 32}\|A\|_2,
\end{split}
\]
since \(\|p\|^2 \le 2 N\) and the Frobenius norm is bounded by $\sqrt N$ the operator norm.
}
\begin{equation}
\label{eq:dominated-energy}
\begin{split}
&\sum_{x,y \in \bZ_N} |\mathfrak l_{x,y}| \le \Const\trn{A}N^{\frac{3}{2}}, \phantom{==} \sum_{x\in \bZ_N} |\mathfrak l_{x,x}| + |\mathfrak m_x| \le \Const\trn{A}  N,\\
& \trn{\mathfrak m_x}+ \trn{\mathfrak l_{x,y}} \le \Const \trn{A}.
\end{split}
\end{equation}
In addition, by \eqref{eq:5} and \eqref{eq:gen_shape},
  \begin{equation}
    \label{eq:46}
\gamma \cal{S} A = \gamma \sum_{x \in \bZ_N} \left(\mathfrak l_{x,x} - \mathfrak m_x  \right).
  \end{equation}
Using \eqref{eq:unperturbed-matrices}, one has
\[
    \begin{split}
    &\langle  DA (G_{\ell'}(\theta)), \bB_{\ve}(s,\theta)\bB_{\ve}(u, \theta)G_{\ell'}(\theta) \rangle \\
    &= - \sum_{x\in \bZ_N} b(f^{\lf{s\ve^{-1}}}\theta_x)b(f^{\lf{u\ve^{-1}}}\theta_x)\mathfrak m_x(G_{\ell'}(\theta)),
    \end{split} 
\]
so that the expression in the third line of the Equation \eqref{eq:support} is equal to 
\begin{equation}
\label{eq:one}
    \begin{split}
    -\ve^{-1}\sum_{\ell' \in \mathfrak{f}_{h}}p_{\ell'} \sum_{x\in \bZ_N}&\int_{\bI_{\ell'}}\int_{0}^h \int_0^s\mathfrak m_x(G_{\ell'}(\theta))\\
    &\times b(f^{\lf{s\ve^{-1}}}\theta_{x})b(f^{\lf{u\ve^{-1}}}\theta_x) \rho_{\ell'}(\theta)\; d \theta\; ds\; du.
    \end{split}
\end{equation}
Since \(|I^x_{\ell'}|\leq \delta_\ve\), see \eqref{eq:delta_choice}, the intermediate value theorem yields 
\[
    \begin{split}
    &\int_{\bI_{\ell'}}\int_{0}^h \int_0^s \mathfrak m_x(G_{\ell'}(\theta)) b(f^{\lf{s\ve^{-1}}}\theta_{x})b(f^{\lf{u\ve^{-1}}}\theta_x)
    \rho_{\ell'}(\theta)d \theta dsdu\\
    &=  \int_{\bI_{\ell'}}\int_{0}^h \int_0^s b(f^{\lf{s\ve^{-1}}}\theta_{x})b(f^{\lf{u\ve^{-1}}}\theta_x)
    \rho_{\ell'}(\theta)d \theta dsdu\\
    &\phantom{=+}\times\int_{\bI_{\ell'}}\!\!\mathfrak m_x(G_{\ell'}(\theta')) \rho_{\ell'}(\theta')d\theta'
    + \cO\left(\ve^{\frac 12}h^2\delta_{\ve}N\trn{A}\right)\\
    &=  \int_{I^x_{\ell'}}\int_{0}^h \int_0^s b(f^{\lf{s\ve^{-1}}}\theta_{x})b(f^{\lf{u\ve^{-1}}}\theta_x)\rho_{\ell',x}(\theta_x)d \theta_x dsdu\\
    &\phantom{=+}\times\int_{\bI_{\ell'}}\!\!\mathfrak m_x(G_{\ell'}(\theta')) \rho_{\ell'}(\theta')d\theta'
    + \cO\left(\ve^{\frac 12}h^2\delta_{\ve}N\trn{A}\right),
\end{split}
\]
where we first used Lemma \ref{lem:BV-estimate}, and then the fact that $\rho_{\ell'}$ is a product of probability densities.\\
Accordingly, \eqref{eq:one} equals
\[
\begin{split}
 &-\ve^{-1}\sum_{\ell' \in \mathfrak{f}_{h}}p_{\ell'} \sum_{x\in \bZ_N}\int_{I^{x}_{\ell'}}\int_{0}^h \int_0^s b(f^{\lf{s\ve^{-1}}}\theta_{x})b(f^{\lf{u\ve^{-1}}}\theta_x)
 \rho_{\ell',x}(\theta)d \theta_x dsdu\\
&\times \int_{\bI_{\ell'}}\!\!\mathfrak m_x(G_{\ell'}(\theta')) \rho_{\ell'}(\theta')d\theta' + \cO\left(N^2\ve^{15}h^2\trn{A}\right).
\end{split}
\]
By Lemma \ref{lem:covariance}, applied, for each $x\in\bZ_N$, with $I=I^x_{\ell'}$ and $\rho=\rho_{\ell',x}$ to the first factor of the corresponding term in the sum above, and the first line of \eqref{eq:dominated-energy}, for each $\ell'$ we have
\begin{equation}
\label{eq:drift-strato}
\begin{split}
&\ve^{-1} \int_{\bI_{\ell'}}\int_{0}^h \int_0^s\left\langle DA (G_{\ell'}(\theta)), \bB_{\ve}(s, \theta)\bB_{\ve}(u, \theta)G_{\ell'}(\theta) \right\rangle \rho_{\ell'}(\theta)d\theta\, dsdu\\
&=-h\gamma \int_{\bI_{\ell'}} \sum_{x\in \bZ_N}\mathfrak m_x(G_{\ell'}(\theta)) \rho_{\ell'}(\theta)d\theta \\
&\phantom{=}+ \cO\left(\ve\trn{A} N (\ln \ve^{-1})^2\right)+\cO\left(N^2\ve^{15}h^2\trn{A}\right).
    \end{split}
\end{equation}
%%%%%%%%%
\subsubsection{Ito-like diffusion} \ \\
Finally, the fourth line of \eqref{eq:support} yields the last term of \eqref{eq:gen_shape}. We can write it as
\begin{equation}
\label{eq:two}
   \begin{split}
      &\sum_{\ell' \in \mathfrak f_h}p_{\ell'}\sum_{x\in \bZ_N} \int_0^h \int_0^s\int_{\bI_{\ell'}}\frac{b(f^{\lf{u\ve^{-1}}}\theta_{x})b(f^{\lf{s\ve^{-1}}}\theta_{x})}\ve \\ &\phantom{\sum_{\ell' \in \mathfrak f_h}p_{\ell'}\sum_{x\in \bZ_N} \int_0^h \int_0^s\int_{\bI_{\ell'}}}
      \times\mathfrak l_{x,x}(G_{\ell'}(\theta))\rho_{\ell'}(\theta)d\theta dsdu\\
     &+ \sum_{\ell' \in \mathfrak f_h}p_{\ell'}\sum_{\substack{x,y \in \bZ_N \\ x \neq y}}\!\int_0^h \!\!\int_0^s\!\!\int_{\bI_{\ell'}}\frac{b(f^{\lf{u\ve^{-1}}}\theta_{x})b(f^{\lf{s\ve^{-1}}}\theta_{y})}\ve \\
&\phantom{\sum_{\ell' \in \mathfrak f_h}p_{\ell'}\sum_{x\in \bZ_N} \int_0^h \int_0^s\int_{\bI_{\ell'}}}
\times \mathfrak l_{x,y}(G_{\ell'}(\theta))\rho_{\ell'}(\theta)d\theta dsdu.
   \end{split} 
\end{equation}
We start with the second term of \eqref{eq:two}. We step back in time and
use again \(  \sum_{\ell' \in \mathfrak f_h}p_{\ell'}\mu_{\ell'} = \phi^h_{\star}\mu_{\ell}\), obtaining
\[
\begin{split}
    &\sum_{\ell' \in \mathfrak f_h}p_{\ell'}\sum_{\substack{x,y \in \bZ_N \\ x \neq y}}\!\int_0^h \!\!ds \int_0^s\!\!du\int_{\bI_{\ell'}}\frac{b(f^{\lf{u\ve^{-1}}}\theta_{x})b(f^{\lf{s\ve^{-1}}}\theta_{y})}\ve\mathfrak l_{x,y}(G_{\ell'}(\theta))\rho_{\ell'}(\theta)d\theta  =\\
    &=\!\!\!\sum_{\substack{x,y \in \bZ_N \\ x \neq y}}\!\int_0^h \!\!\int_0^s\!\!\int_{\bI_{\ell}}\frac{b(f^{\lf{(u+h)\ve^{-1}}}\theta_{x})b(f^{\lf{(s+h)\ve^{-1}}}\theta_{y})}\ve\mathfrak l_{x,y}(z(h,G_{\ell}(\theta)))\rho_{\ell}(\theta)d\theta dsdu.
\end{split}
\]
Using Lemma \ref{lem:Taylor-0} for any \(x,y\), and \eqref{eq:dominated-energy} we have
\[
\begin{split}
\Big\|\mathfrak l_{x,y}(z(h,G_{\ell}(\theta))) - &\mathfrak l_{x,y}(G_{\ell}(\theta))
-\int_0^h\langle \nabla l_{x,y}(G_{\ell}(\theta)), (-\bA+\ve^{-\frac{1}{2}}\bB_\ve(s)) G_{\ell}(\theta)\rangle ds\Big\| \\
&\le \Const \trn{A} N ( h^2\ve^{-1} + h^3 \ve^{-\frac{3}{2}})\le \Const \trn{A} N h^2\ve^{-1},
\end{split}
\]
where we have used that $h\in (\ve^{\frac 56}, 2\ve^{\frac 56})$.
Treating the term containing \(\bA\) as an error, the second term of \eqref{eq:two} then becomes
\begin{equation}\label{eq:ii-420}
\begin{split}
&\sum_{\substack{x,y \in \bZ_N \\ x \neq y}}\!\!\int_0^h \!\!\!ds\!\int_0^s\!\!\! du\!\int_{\bI_{\ell}}\!\tfrac{b(f^{\lf{(u+h)\ve^{-1}}}\theta_{x})b(f^{\lf{(s+h)\ve^{-1}}}\theta_{y})}\ve\mathfrak l_{x,y}(G_{\ell}(\theta))\rho_{\ell}(\theta)d\theta\\
&+\ve^{-\frac 32}\sum_{\substack{x,y,z \in \bZ_N \\ x \neq y}}\!\int_0^h ds \!\!\int_0^s\!\!du\int_0^h\!\!d\tau\int_{\bI_{\ell}}b(f^{\lf{(u+h)\ve^{-1}}}\theta_{x})b(f^{\lf{(s+h)\ve^{-1}}}\theta_{y})\\
&\times\langle \partial_{p_z} l_{x,y}(G_{\ell}(\theta)),JG_{\ell, p_z}(\theta)\rangle b(f^{\lf{\tau\ve^{-1}}}\theta_{z}) \rho_{\ell}(\theta)d\theta \\
&+\trn{A}N^3\cO( h^4 \ve^{-2} +h^3\ve^{-1}).
\end{split}
\end{equation}
By Lemma \ref{lem:BV-estimate}, the first term in \eqref{eq:ii-420} can be written as
\begin{equation}\label{eq:frist_421}
\begin{split}
&\ve^{-1}\sum_{\substack{x,y \in \bZ_N \\ x \neq y}}\!\int_0^h \!\!ds\int_0^s \!\!du 
\Bigg\{ \int_{I_{\ell}^x}\!\!\!\! b(f^{\lf{(u+h)\ve^{-1}}}\theta_{x})\rho_{\ell,x}(\theta_x)d\theta_x \\ &\phantom{=====}\times\int_{I_{\ell}^y}\!\!\!\!b(f^{\lf{(s+h)\ve^{-1}}}\theta_{y}) \rho_{\ell,y}(\theta_y) d\theta_y\Bigg\}
\int_{\bI_{\ell}}\mathfrak l_{x,y}(G_{\ell}(\theta))\rho_{\ell}(\theta)\\
&\phantom{=====}+ \cO(\trn{A}h^2N^3\ve^{-\frac{1}{2}}\delta_\ve)
\end{split}
\end{equation}
Moreover, by Lemma \ref{lem:dec_cor} and \ref{lem:BV-standard-pair}, for any \(u \ge 0\),
\begin{equation}\label{eq:exp-small-avg}
\begin{split}
\int_{I^x_{\ell}}b(f^{\lf{(u+h)\ve^{-1}}}\theta_x)\rho_{\ell, x}(\theta_x)d\theta_x  &\le \Const \left\|\rho_{\ell, x} \Id_{I^x_{\ell}}\right\|_{BV}\nu^{\lf{(h+u)\ve^{-1}}} \\
&\le \Const  \delta_\ve^{-1}\nu^{\ve^{-\frac 16}} .
\end{split}
\end{equation}
Hence, using the second of \eqref{eq:dominated-energy}, we estimate the first term of \eqref{eq:two} as
\begin{equation}\label{eq:second_line} 
\begin{split}
\sum_{\ell' \in \mathfrak f_h}p_{\ell'}\sum_{\substack{x,y \in \bZ_N \\ x \neq y}}&\int_0^h \int_0^s\int_{\bI_{\ell'}}\frac{b(f^{\lf{(u+h)\ve^{-1}}}\theta_{x})b(f^{\lf{(s+h)\ve^{-1}}}\theta_{y})}\ve \\
&\times\mathfrak l_{x,y}(G_{\ell'}(\theta))\rho_{\ell'}(\theta)d\theta  dsdu=  \cO(\trn{A}h^2N^3\ve^{-\frac{1}{2}}\delta_\ve).
\end{split}
\end{equation}
Next, we can treat the second term of \eqref{eq:ii-420} similarly and write it as
\begin{equation}\label{eq:ii-41}
\begin{split}
\ve^{-\frac 32}\sum_{\substack{x,y \in \bZ_N \\ x \neq y}}\sum_{z\in \bZ_N}\!\int_{[0,h]^2}\hskip-6pt &ds\, d\tau \int_0^s\!\!du
\Bigg\{ \int_{\bI_{\ell}}b(f^{\lf{(u+h)\ve^{-1}}}\theta_{x})\\
&\times b(f^{\lf{(s+h)\ve^{-1}}}\theta_{y}) b(f^{\lf{\tau\ve^{-1}}}\theta_{z}) \rho_{\ell}(\theta)d\theta\Bigg\}\\
&\times\int_{\bI_{\ell}}\langle \partial_{p_z} l_{x,y}(G_{\ell}(\theta)),JG_{\ell, p_z}(\theta)\rangle  \rho_{\ell}(\theta)d\theta\\
&\hskip-3cm+\cO(\ve^{- 1}h^3N^3\delta_\ve)=\cO(\ve^{ -1}h^3N^3\delta_\ve)
\end{split}
\end{equation}
where, in the last equality, we have used \eqref{eq:exp-small-avg} again since either $x$ or $y$ must be different from $z$.
Collecting the estimates in \eqref{eq:ii-420}, \eqref{eq:second_line} and \eqref{eq:ii-41} we have
\begin{equation}\label{eq:second_term_final}
\begin{split}
&\sum_{\ell' \in \mathfrak f_h}p_{\ell'}\!\!\!\sum_{\substack{x,y \in \bZ_N \\ x \neq y}}\!\int_0^h \!\!ds\int_0^s\!\! du\int_{\bI_{\ell'}}\!\!\tfrac{b(f^{\lf{u\ve^{-1}}}\theta_{x})b(f^{\lf{s\ve^{-1}}}\theta_{y})}\ve \mathfrak l_{x,y}(G_{\ell'}(\theta))\rho_{\ell'}(\theta)d\theta\\
&=\cO( \trn{A}N^3h^4 \ve^{-2} )\leq \cO( \trn{A}N h^3 \ve^{-\frac 32} ),
\end{split}
\end{equation}
where, in the last inequality, we have used $h\leq 2\ve^{\frac 56}$ and $\ve\leq N^{-6}$.\\
It remains to consider the first term in \eqref{eq:two}. Using again Lemma \ref{lem:BV-estimate},
\begin{equation}
\label{eq:three}
    \begin{split}
    &\sum_{x\in \bZ_N}\int_0^h \int_0^s\int_{\bI_{\ell'}}\frac{b(f^{\lf{u\ve^{-1}}}\theta_{x})b(f^{\lf{s\ve^{-1}}}\theta_{x})}\ve \mathfrak l_{x,x}(G_{\ell'}(\theta))\rho_{\ell'}(\theta)d\theta dsdu \\
&=\sum_{x\in \bZ_N}\int_0^h \int_0^s\int_{I^x_{\ell'}}\frac{b(f^{\lf{u\ve^{-1}}}\theta_{x})b(f^{\lf{s\ve^{-1}}}\theta_{x})}\ve \rho_{\ell', x}(\theta_x)d\theta_x  dsdu \\
&\phantom{=}
\times\int_{\bI_{\ell'}}\mathfrak l_{x,x}(G_{\ell'}(\theta))\rho_{\ell'}(\theta)d\theta + \cO(N^3\ve^{15} h^2\trn{A}).
\end{split}
\end{equation}
Finally, by Lemma \ref{lem:covariance}, applied as in \eqref{eq:drift-strato}, and the first of \eqref{eq:dominated-energy}, Equation \eqref{eq:three} is equal to
\begin{equation}
\label{eq:Ito-term}
\begin{split}
     h\gamma\int_{\bI_{\ell'}} \sum_{x\in \bZ_N}\mathfrak l_{x,x}&(G_{\ell'}(\theta))\rho_{\ell'}(\theta)d\theta+ \cO(N\ve (\ln \ve^{-1})^2) \trn{A}) \\
&+ \cO(N^3\ve^{15} h^2\trn{A}).
    \end{split}
\end{equation}
We now collect all the estimates \eqref{eq:secon_support}, \eqref{eq:drift-strato}, \eqref{eq:second_term_final} and \eqref{eq:Ito-term} together. Recalling \eqref{eq:46}, yields
\[
      \mu_{\ell}(A\left(z(2h))\right) = \mu_{\ell}\left( A\left(z(h)\right)\right)+ \mu_{\ell}\left(\cG A\left(z(h)\right)\right)h + \err_A(h,\ve, N)h ,
\]
where we have used \(h \in (\ve^{\frac{5}{6}}, 2\ve^{\frac{5}{6}})\) and \(\ve \in (0, N^{-6})\).
Hence, by Lemma \ref{lem:stochastic-econtrolled} and Corollary \ref{lem:stochastic-Taylor},
\begin{equation}
\label{eq:2h}
\begin{split}
     \mu_{\ell}(A\left(z(2h))\right) =& \mu_{\ell}\left( A\left(z(h)\right)\right)+ \int_{h}^{2h}\mu_{\ell}\left(\cG A\left(z(s)\right)\right)ds\\
     & +\err_A(h,\ve, N) h,
\end{split}
\end{equation}
concluding the proof.
%%%%%%%%%
\subsection{Average of a first order term}\label{sec:aveft}\ \\
Here, we show that the contribution of the second line of \eqref{eq:support}  is negligible.
\begin{lemma}
    \label{lem:annoying-term}
For \(\ve \in (0, N^{-6})\), \(h \in (\ve^{\frac{5}{6}}, 2\ve^{\frac 56})\cap\ve\bN\), \(\ell \in \mathfrak S^*(N,\ve)\), and  \(A \in \mathfrak{E}_N\),\\
we have
\[
\begin{split}
&\ve^{-\frac12}\mu_{\ell}\left(\int_h^{2h} \!\!\!\!\left\langle  DA(z(h)),  \bB_{\ve}(s) z(h)\right\rangle ds\right) = \cO\left( \trn{A} N \left[\ve (\ln \ve^{-1})^2  + h^{3}\ve^{-\frac 32} \right]\right).
\end{split}
\]
\end{lemma}
\begin{proof}
Expanding \( DA (z(h))\) and \(z(h)\) as in \eqref{eq:generator},
for \(s \in [h,2h]\),
\[
\begin{split}
\left\langle  DA(z(h)),  \bB_{\ve}(s) z(h)\right\rangle \! &= \! \biggl\langle DA+ \int_{0}^h \!\!\! D^2 A \left(-\bA + \frac{\bB_{\ve}(u)}{\sqrt{\ve}}\right)z(u) du,\\
&\qquad\qquad\bB_{\ve}(s) \biggl(z + \int_0^h \left(-\bA + \frac{\bB_{\ve}(u)}{\sqrt{\ve}}\right)z(u)du\biggr)\biggr \rangle.
\end{split}
\]
By \eqref{eq:rough_est} and Lemma \ref{lem:matrix-norms} we can estimate the difference between \(z(u)\) and the initial condition \(z(0)=z\), and write
\begin{equation}\label{eq:annoying_term_expansion}
\begin{split}
 &\frac{\mu_{\ell}\left(\int_h^{2h}\left\langle  DA(z(h)),  \bB_{\ve}(s) z(h)\right\rangle ds\right)}{\sqrt{\ve}}=\frac{\mu_{\ell} \left( \int_h^{2h}\left\langle  DA(z) ,  \bB_{\ve}(s) z\right\rangle ds \right)}{\sqrt{\ve}} \\
&+\ve^{-\frac 12}\mu_{\ell} \left(\int_0^h   \!\!\int_{0}^h  \!\! \left\langle  DA(z), \bB_{\ve}(h+s)  \left(-\bA \!+ \!\ve^{-\frac 12}\bB_{\ve}(u)\right)z\right\rangle dsdu\right)\\
&+\ve^{-\frac 12}\mu_{\ell}\left(\int_{0}^{h} \! \!\int_0^h \! \!\left\langle  D^2A\left(-\bA +\ve^{-\frac 12}\bB_{\ve}(u)\right)z,  \bB_{\ve}(h+s) z\right\rangle dsdu\right)\\ 
& + \cO\left(\trn{A}Nh^3\ve^{-\frac 32}\right).
\end{split}
\end{equation}
We analyze the above terms one by one. Let us set \(\mathfrak g_x(\bol q, \bol p) = \langle\partial_{p_x} A(\bol q, \bol p), Jp_x\rangle\), \(x \in \bZ_N\). Note that \(\mathfrak g_x \in \fE_N\), moreover
\begin{equation}
\label{eq:controlled-energy-2}
\sum_{x \in \bZ_N} |\mathfrak g_x| \le \trn{A}  N \quad \text{and} \quad \trn{\mathfrak g_x} \le \trn{A}.
\end{equation}
Furthermore,
\[
        \begin{split}
       &\langle  DA(G_{\ell}(\theta)), \bB_{\ve}(h+s)G_{\ell}(\theta) \rangle =  \sum_{x\in \bZ_N}b(f^{\lf{(h+s)\ve^{-1}}}\theta_x) \mathfrak{g}_x(G_{\ell}(\theta)).
        \end{split}
\]
Using the second of \eqref{eq:controlled-energy-2}, \eqref{eq:delta_choice} and Lemma \ref{lem:BV-estimate}, the first term of \eqref{eq:annoying_term_expansion} equals
    \begin{equation}
\label{eq:lem-BB-1}
        \begin{split}
& \ve^{-\frac 12}\sum_{x\in \bZ_N} \int_{0}^h \int_{\bT^{N}}b (f^{\lf{(h+s)\ve^{-1}}}\theta_x)  \mathfrak{g}_x(G_\ell(\theta)) \Id_{\bI_{\ell}}(\theta)\rho_{\ell}(\theta)d\theta=\\
&= \ve^{- \frac 12}\sum_{x\in \bZ_N} \int_{0}^h \int_{\bT}b(f^{\lf{(h+s)\ve^{-1}}}\theta_x) \Id_{I^x_{\ell}}(\theta)\rho_{\ell, x}(\theta)d\theta_x\\
&\phantom{=\sum_{x\in \bZ_N}}
\times \int_{\bI_{\ell}}\mathfrak g_x (G_{\ell}(\theta))\rho_{\ell}(\theta)d\theta
+ \cO(h \delta_\ve N^{2}\trn{A}),
        \end{split}
   \end{equation}
and by \eqref{eq:exp-small-avg} the sum above integrals is exponentially small in \(\ve\).
Therefore, also using the first of \eqref{eq:controlled-energy-2}, it follows that
\begin{equation}\label{eq:primo_termine}
\ve^{-\frac 12}\left|\mu_{\ell}\left( \int_h^{2h}\left\langle  DA(z), \left( \bB_{\ve}(s) \right)z\right\rangle ds \right)\right| 
\leq \Const  \delta_\ve hN^2 \trn{A}.
\end{equation}
Using \eqref{eq:rough_est}, the second term of \eqref{eq:annoying_term_expansion} reads
\[
\begin{split}
&\ve^{-\frac 12}\mu_\ell \left(\int_0^h \int_0^h \langle  DA(z),  \bB_{\ve}(h+s)  \left(-\bA + \ve^{-\frac 12}\bB_{\ve}(u)\right)z \rangle dsdu \right)\\
&= \ve^{-1}\mu_\ell \left( \int_0^h \int_0^h\left\langle  DA(z),  \bB_{\ve}(h+s) \bB_{\ve}(u)z \right\rangle dsdu\right)+\cO\left(\ve^{-\frac 32}h^3 N \trn{A}\right),
\end{split}
\]
where, in the last equality, we have used $h>\ve$.
Let $\mathfrak n_x(\bol q , \bol p)=-\langle p_x, \partial_{p_x}A\rangle$, \(x \in \bZ_N\). Note that \(\mathfrak n_x \in \fE_N\), and, 
\begin{equation}
\label{eq:controlled-energy-3}
\begin{split}
\sum_{x \in \bZ_N}&|\mathfrak n_x|\le \trn{A}N \quad \text{and}\quad \trn{\mathfrak n_x} \le \trn{A}.
\end{split}
\end{equation}
We have
\[
\begin{split}
&\ve^{-1}\int_{0}^h\int_0^h\mu_{\ell} \left(\langle  DA(z),  \bB_{\ve}(h+s)   \bB_{\ve}(u)z\rangle\right) dsdu \\
&=\ve^{-1}\sum_{x \in \bZ_N}\int_{0}^{h} \int_{0}^h 
\int_{\bI_{\ell}}b(f^{\lf{(h + s)\ve^{-1}}}\theta_x)b(f^{\lf{u\ve^{-1}}}\theta_x)\mathfrak n_x(G_{\ell}(\theta))\rho_{\ell}(\theta)d\theta dsdu.
\end{split}
\]
Using Lemma \ref{lem:BV-estimate} and the second of \eqref{eq:controlled-energy-3}, this is equal to
\begin{equation}
\label{eq:gamma-appears}
\begin{split}
&\ve^{-1}\sum_{x \in \bZ_N}\int_{0}^{h} \int_{0}^h  \int_{I^x_{\ell}}b(f^{\lf{(h + s)\ve^{-1}}}\theta_x)b(f^{\lf{u\ve^{-1}}}\theta_x)\rho_{\ell, x}(\theta_x)d\theta_xdsdu\\
&\times \int_{\bI_{\ell}}\mathfrak n_x(G_{\ell}(\theta))\rho_{\ell}(\theta)d\theta + \cO(\delta_\ve\ve^{-1}h^2 N^{2} \trn{A}).
\end{split}
\end{equation}
By Lemma \ref{lem:covariance2} and the first of \eqref{eq:controlled-energy-3}, the first term of the above is bounded by  
\(\Const \trn{A}N\ve(\ln \ve^{-1})^2\). Hence,
\begin{equation}\label{eq:secondo_termine}
\begin{split}
&\ve^{-\frac 12}\left|\mu_{\ell} \left(\int_0^h\int_0^h\left\langle  DA(z), \bB_{\ve}(h+s)  \left(\bA + \ve^{-\frac 12}\bB_{\ve}(u)\right)z(u)\right\rangle dsdu\right)\right|\\
&\phantom{\frac{1}{\sqrt\ve}\mu_{\ell'}}
\leq \Const N \trn{A}\left(\ve(\ln \ve^{-1})^2 + \ve^{-\frac 32}h^3\right).
\end{split}
\end{equation}
Next, we analyze the third term in \eqref{eq:annoying_term_expansion}. By \eqref{eq:rough_est}, we can write it as
\begin{equation}
\label{eq:BB}
\begin{split}
\ve^{-1}\mu_{\ell}&\left(\int_{0}^{h}\int_{0}^{h}\left\langle D^2 A \bB_{\ve}(s) z, \bB_{\ve}(h+u) z \right\rangle dsdu\right)\\
&+\cO(\trn{A}\ve^{-\frac 32}h^3 N ).
\end{split}
\end{equation}
For \(x,y \in \bZ_N\), recall that \(\mathfrak l_{x,y} = \left\langle \frac{\partial^2 A}{\partial p_{x}\partial p_{y}}Jp_{y}, J p_{x}\right\rangle\) and the estimates in \eqref{eq:dominated-energy}.
With the above notation we can write the first line of \eqref{eq:BB} as
\[
\begin{split}
&\ve^{-1}\sum_{x \in \bZ_N}\sum_{y \in \bZ_N}\int_0^h \int_0^h  \int_{\bI_{\ell}}  \mathfrak l_{x,y}(G_{\ell}(\theta)) b(f^{\lf{s\ve^{-1}}}\theta_y)b(f^{\lf{(u+h)\ve^{-1}}}\theta_x)\rho_{\ell}(\theta)d\theta dsdu.
\end{split}
\]
By Lemma \ref{lem:BV-estimate} and the second of \eqref{eq:dominated-energy}, this is equal to
\begin{equation}
\label{eq:BB2}
\begin{split}
&\ve^{-1}\sum_{x \in \bZ_N}\sum_{y \in \bZ_N}\int_0^h \int_0^h  \int_{\bI_{\ell}} b(f^{\lf{s\ve^{-1}}}\theta_y)b(f^{\lf{(u+h)\ve^{-1}}}\theta_x)\rho_{\ell}(\theta)d\theta dsdu\\
&\times\int_{\bI_{\ell}} \mathfrak l_{x,y}(G_{\ell}(\theta))  \rho_{\ell}(\theta)d\theta +  \cO\left(N^{4} \delta_\ve\ve^{-1}\trn{A} h^2\right).
\end{split}
\end{equation}
For \(x \neq y\), since \(\rho_{\ell}\) is a product measure,
\[
\begin{split}
\int_{\bI_{\ell}} b(f^{\lf{s\ve^{-1}}}\theta_y)&b(f^{\lf{(u+h)\ve^{-1}}}\theta_x)\rho_{\ell}(\theta)d\theta \\
&= \int_{I^y_{\ell}}b(f^{\lf{s\ve^{-1}}}\theta_y)\rho_{\ell, x}(\theta_y)d\theta_y\int_{I^x_{\ell}}b(f^{\lf{(u+h)\ve^{-1}}}\theta_x)\rho_{\ell, x}(\theta_x)d\theta_x.
\end{split}
\]
However, by Lemma \ref{lem:dec_cor} and \ref{lem:BV-standard-pair}, for any \(u \ge 0\),
\[
\begin{split}
\int_{I^x_{\ell}}b(f^{\lf{(u+h)\ve^{-1}}}\theta_x)&\rho_{\ell, x}(\theta_x)d\theta_x = \int_{\bT}b(\theta_x) \cL^{\lf{(h+u)\ve^{-1}}}\left(\rho_{\ell, x}(\theta_x) \Id_{I^x_{\ell}}(\theta_x)\right)d\theta_x \\
& \le \Const \left\|\rho_{\ell, x} \Id_{I^x_{\ell}}\right\|_{BV}\nu^{\lf{(h+u)\ve^{-1}}} \le \Const  \delta_\ve^{-1}\nu^{\ve^{-\frac 16}} .
\end{split}
\]
For \(x=y\),
\[
\begin{split}
\ve^{-1}\int_0^h \int_0^h \int_{\bI_{\ell}}&b(f^{\lf{s\ve^{-1}}}\theta_y)b(f^{\lf{(u+h)\ve^{-1}}}\theta_x)\rho_{\ell}(\theta)d\theta dsdu=\\
&=\ve^{-1}\int_0^h \int_0^h \int_{I^x_{\ell}} b(f^{\lf{s\ve^{-1}}}\theta_x)b(f^{\lf{(u+h)\ve^{-1}}}\theta_x)\rho_{\ell, x}(\theta_x)d\theta_x dsdu,
\end{split}
\]
and by Lemma \ref{lem:covariance2} this is of order \(\cO(\ve (\ln \ve^{-1})^2)\).\\
Therefore, recalling the definition \eqref{eq:delta_choice} of $\delta_\ve$, Equations \eqref{eq:BB}, \eqref{eq:BB2} and using the first of \eqref{eq:dominated-energy}, we have
\begin{equation}\label{eq:terzo_termine}
\begin{split}
\ve^{-\frac 12}\mu_{\ell}&\left(\int_{0}^{h}\int_0^h\left\langle  D^2A\left(-\bA +\ve^{-\frac 12}\bB_{\ve}(s)\right)z(s),  \bB_{\ve}(h+u) z\right\rangle dsdu\right)\\
&\leq  \Const \trn{A}N\left(  \ve (\ln \ve^{-1})^2+\ve^{-\frac 32}h^3 \right).
\end{split}
\end{equation}
Collecting the estimates \eqref{eq:primo_termine}, \eqref{eq:secondo_termine} e \eqref{eq:terzo_termine} the lemma follows.
\end{proof}
%%%%%%%%%%%%%%

\section{Hydrodynamic limit}
\label{sec:hydrodynamic}
For any test function $\varphi\in\cC^0(\bT,\bR)$, we define
\begin{equation}
\label{eq:6}
\xi_N(t,\varphi) = \frac 1N \sum_{x\in\bZ_N} \varphi\left(\frac xN\right) \bE_{\mu_N}(\fe_x(N^2 t)).
\end{equation} 
Alternatively, we can see $\xi_N(t)$ as a probability measure on $\bT$.\\
In order to shorten notations, we denote the error terms as
\begin{equation}
  \label{eq:err}
  R(N,\ve,t) := t \left[N^{-1}+N \ve^{\frac 16}(\ln\ve^{-1})^2 \right]+ N^{-1} .
\end{equation}

\begin{theorem}
\label{prop:fundamental}
For each  $\ve \in(0,N^{-6})$, $N\in\bN$, and initial probability measure 
$\mu_N\in\cM_{\textrm{init}}$ we have, for any $s,t\in\R_+$, \(s < t\), and $\vf\in\cC^5(\bT,\bR)$,
\[
\begin{split}
  \xi_N(t,\varphi) -  \xi_N(s,\varphi) =\frac{D}{2\gamma} \int_s^t  \xi_N(u,\varphi'') du
  + \|\vf\|_{\cC^5}\cO\left(R(N,\ve,t-s) \right).
  \end{split}
\]
\end{theorem}

\begin{proof}
Equations \eqref{eq:current}, \eqref{eq:3} imply, 
\[
\begin{split}
  {\xi_N (t,\varphi)-\xi_N(s,\varphi)}&=
  N\int_s^t \!\!\sum_{x\in\bZ_N} \!\!\varphi\left(\frac xN\right) \bE_{\mu_N}\Big(j_{x-1,x}(N^2 u)-j_{x,x+1}(N^2 u)\Big) du\\
  &=N\int_s^t \sum_{x\in\bZ_N} \nabla_{\bZ_N}\varphi\left(\frac xN\right) \bE_{\mu_N}(j_{x,x+1}(N^2 u)) du
\end{split}
\]
where, for $g:\bZ_N\to\bR$, 
\[
\nabla_{\bZ_N} g(x)=g(x+1)-g(x).
\]
In Proposition \ref{prop:second_der}, Section \ref{sec:fluct-diss-relat}, we use a fluctuation-dissipation relation to obtain that
\begin{equation}\label{eq:intermediate-heat-step}
\begin{split}
 \xi_N(t,\varphi) -  \xi_N(s,\varphi) =&\frac 1{2\gamma} \frac 1N \sum_{x\in\bZ_N} \varphi''\left(\frac xN\right) \int_s^t \bE_{\mu_N}(\mathfrak u_x (N^2 u)) du\\
 &+\cO(\|\vf\|_{\cC^3}R(N,\ve,t-s)),
\end{split}
\end{equation}
where
\[
\mathfrak u_x =  p_x^2 - \omega_0^2 q_x^2 + (q_{x+1} - q_x) (q_x - q_{x-1}).
\]
Next, we want to establish a local equilibrium results. 
This means that \emph{locally} the distribution of the $(q_x,p_x)$ is \emph{close} to 
the equilibrium measure of the infinite dynamics on $\Z$ at the temperature $T$ corresponding to 
the local average of the kinetic energy $p_x^2$. This equilibrium measure $\nu_T$
is given by the Gibbs--Gaussian measure on $\bZ$ with covariance, for $x,y\in \bZ$ 
and $i,j=1,2$,
\begin{equation}\label{eq:gaussian-equilibrium-chain}
\begin{split}
    &\mathbb E_{\nu_T}\left(q_{x}^{i} q_{y}^{j}\right) = T \delta_{i,j}\bG_{\omega_0}(x-y) ,\\
   & \mathbb E_{\nu_T}\left(p_{x}^{i} p_{y}^{j}\right) = T \delta_{i,j}\delta_{x,y},\\
   & \mathbb E_{\nu_T}\left(q_{x}^{i} p_{y}^{j}\right) = 0,
\end{split}
\end{equation}
 where $\bG_{\omega_0}(x)$ is the Green function of the operator $\left(\omega_0^2 - \Delta\right)^{-1}$, that is
\begin{equation}\label{eq:green_def}
\left[\left(\omega_0^2 - \Delta\right)^{-1}\vf\right](x)=\sum_{y\in\bZ}\bG_{\omega_0}(x-y)\vf(y),
\end{equation}
see Appendix \ref{sec:the-diff-expr} for more details. The expectation of  $\mathfrak u_x$ with respect to $\nu_T$
is then given by
\begin{equation}
  \label{eq:D1}
\begin{split}
      &T\left[ 1 - \omega_0^2 \bG_{\omega_0}(0) + 2\bG_{\omega_0}(1) - \bG_{\omega_0}(0) - \bG_{\omega_0}(2)\right]=
  \\
  &\phantom{=================}=  T\left[1 - \omega_0^2\left(\bG_{\omega_0}(0) + \bG_{\omega_0}(1)\right)\right],
\end{split}
\end{equation}
where we have used Equation \ref{eq:28}.
This identifies the diffusion coefficient 
\begin{equation}\label{eq:D2}
D= \left[1 - \omega_0^2\left(\bG_{\omega_0}(0) + \bG_{\om_0}(1)\right)\right].
\end{equation}
Proposition \ref{prop:loc-equi} in Section \ref{sec:local-equilibrium}, together with \eqref{eq:intermediate-heat-step}, \eqref{eq:D2} and \eqref{eq:D1}, implies
\[
\begin{split}
 &\xi_N(t,\varphi) -  \xi_N(s,\varphi) = \frac 1{2\gamma N} \sum_{x \in \bZ_N} \varphi''\left(\frac{x}{N}\right)\\
 &\times
    \int_s^t \bE_{\mu_N}\left(\left[1-\omega_0^2\bG_{\omega_0}(0)+2\bG_{\omega_0}(1)-\bG_{\omega_0}(0)-\bG_{\omega_0}(2)\right]p_x^2 (N^2 u) \right) du \\
 &+\cO(\|\vf\|_{\cC^5}R(N,\ve,t-s))\\
 =& \frac D{2\gamma N} 
    \int_s^t \sum_{x \in \bZ_N}\varphi''\left(\frac{x}{N}\right)\bE_{\mu_N}\left(p_x^2 (N^2 u) \right) du +\cO(\|\vf\|_{\cC^5}R(N,\ve,t-s)).
\end{split}
\] 
It remains to show that the expectation of $p_x^2$ equals the expectation of $\fe_x$. This is a form of equipartition, and it is proven in Proposition \ref{cor012912-21}, Section \ref{sec:equipartition}, yielding 
\[
\begin{split}
\xi_N(t,\varphi) -  \xi_N(s,\varphi) =&\frac D{2\gamma N} 
    \int_s^t \sum_{x \in \bZ_N} \varphi''\left(\frac{x}{N}\right)\bE_{\mu_N}\left(\fe_x (N^2 u) \right) du \\
    &+\cO(\|\vf\|_{\cC^5}R(N,\ve,t-s))\\
    &= \frac{D}{2\gamma} 
    \int_s^t \xi_N(N^2 u, \vf'') du +\cO(\|\vf\|_{\cC^5}R(N,\ve,t-s))
\end{split}
\]
which concludes the proof.
\end{proof}
The rest of the section is devoted to proving the building blocks that we have used in the proof of Theorem \ref{prop:fundamental}.

%%%%%%%%%
\subsection{The fluctuation-dissipation relation.}
\label{sec:fluct-diss-relat}\ \\
 Our first result is the following \emph{fluctuation-dissipation} relation.
\begin{lemma}\label{lem:fluct_diss}
For each $x\in\bZ_N$, we have\footnote{ Recall that \((\nabla_{\bZ_N}g)_x = g_{x+1} -g_x\) and by \((\Delta_{\bZ_N}g)_x = g_{x+1} + g_{x-1} - 2g_x\).}
\begin{equation}
  \label{eq:FD}
  \begin{split}
    &j_{x,x+1} = - \frac 1{2\gamma} \nabla_{\bZ_N} \mathfrak u_x  + \mathcal G F_x,\\
    &\mathfrak u_x =  p_x^2 - \omega_0^2 q_x^2 + (q_{x+1} - q_x)  (q_x - q_{x-1})\\
     & F_x = \frac 1{2\gamma} \left(q_{x+1} - q_x\right)  \left(p_{x+1} + p_x\right) + \frac 14 \left(q_{x+1} - q_x\right)^2.
  \end{split}
\end{equation}
\end{lemma}
\begin{proof}
A direct computation yields
\begin{equation*}
  \gamma \mathcal S F_x = -  \frac 1{2} \left(q_{x+1} - q_x\right)  \left(p_{x+1} + p_x\right) =  + j_{x,x+1} -
  \frac 1{2} \left(q_{x+1} - q_x\right)  \left(p_{x+1} - p_x\right),
\end{equation*}
and
\begin{equation*}
  \begin{split}
    \mathcal A F_x &= \frac 1{2\gamma} \left(p_{x+1}^2 - p_x^2\right) +  \frac 1{2} \left(q_{x+1} - q_x\right)  \left(p_{x+1} - p_x\right)\\
    &\qquad + \frac 1{2\gamma} \left(q_{x+1} - q_x\right)  \left[ \left(\Delta_{\bZ_N} - \omega_0^2\right) \left(q_{x+1} + q_x\right) \right]\\
    &= \frac 1{2\gamma} \nabla_{\bZ_N}\left( p_x^2 - \omega_0^2 q_x^2 + (q_{x+1} - q_x)  (q_x - q_{x-1}) \right)\\
    &\qquad +  \frac 1{2} \left(q_{x+1} - q_x\right)  \left(p_{x+1} - p_x\right)
\end{split}
\end{equation*}
from which the statement follows.
\end{proof}
Due to the above decomposition, we obtain the following result, in which the second derivative of the test function appears.
\begin{proposition}\label{prop:second_der}
For each $\vf\in\cC^3(\bT,\bR)$, $s,t\in\bR_+$, \(s<t\), $N\in\bN$, $\ve\in (0,N^{-6})$, and initial probability measure $\mu_N\in\cM_{\textrm{init}}$, we have
\[
\begin{split}
&N\int_s^t \sum_{x\in\bZ_N} \nabla_{\bZ_N}\varphi\left(\frac xN\right) \bE_{\mu_N}(j_{x,x+1}(N^2 s)) ds=\\
&=\frac 1{2\gamma N} \sum_{x\in\bZ_N} \varphi''\left(\frac xN\right) \int_s^t \bE_{\mu_N}(\mathfrak u_x (N^2 s)) ds+\cO(\|\vf\|_{\cC^3}R(N,\ve,t-s)).
\end{split}
\]
\end{proposition}
\begin{proof}
By Lemma \ref{lem:fluct_diss} we have
\begin{equation}\label{eq:evol1}
\begin{split}
N&\int_s^t \sum_{x\in\bZ_N} \nabla_{\bZ_N}\varphi\left(\frac xN\right) \bE_{\mu_N}(j_{x,x+1}(N^2 s)) ds =\\
 =&\frac 1{2\gamma N} \sum_{x\in\bZ_N} N^2\Delta_{\bZ_N}\varphi\left(\frac xN\right) \int_s^t \bE_{\mu_N}(\mathfrak u_x (N^2 s)) ds\\
 &+ \sum_{x\in\bZ_N} N \nabla_{\bZ_N}\varphi\left(\frac xN\right) \int_s^t \bE_{\mu_N}(\mathcal G F_x (N^2 s)) ds.
\end{split}
\end{equation}
Defining $A_N= \sum_x F_x$, we have $\lvvvert A_N\rvvvert \le \Const$ for all $N$. Hence, by Theorem \ref{thm:standard-pair-invariance}, we can apply Corollary \ref{cor:accellerated_generator} to the last term of \eqref{eq:evol1} obtaining
\begin{equation}
  \label{eq:last}
  \begin{split}
    \|\varphi'\|_{\infty}  \left|\int_s^t \bE_{\mu_N}( \mathcal G A_N (N^2 s)) ds\right| \le \Const R(N,\ve,t-s).
  \end{split}
\end{equation}
That is, the last term of \eqref{eq:evol1} is negligible. The Proposition follows since $\|N^2\Delta_{\bZ_N}\varphi -\varphi''\|_\infty\leq \|\vf\|_{\cC^3}N^{-1}$ and
$\sum_x |\mathfrak u_x| \le \Const N$. 
\end{proof}

\subsection{Local equilibrium}
\label{sec:local-equilibrium}\ \\
The goal of this section is to prove that the space correlation of the positions is determined, modulo an error, by the Gaussian equilibrium distributions \eqref{eq:gaussian-equilibrium-chain} in which the Green function \eqref{eq:green_def} appears.
\begin{proposition}
  \label{prop:loc-equi}
For each $\varphi\in \cC^5(\bT,\bR)$, $t> s\in\bR_+$, $N\in\bN$, $\ve\in (0,N^{-6})$, $w,w'\in\{0,1\}$, and initial probability measure $\mu_N\in\cM_{\textrm{init}}$ we have
\[
\begin{split}  \frac 1N\! \!\sum_{x \in \bZ_N}\! \!\varphi''\left(\frac{x}{N}\right)\!&
    \int_s^t\! \bE_{\mu_N}\left(q_{x+w} (N^2 u) q_{x-w'}(N^2 u) - \bG_{\omega_0}(w+w') p_x^2 (N^2 u)\right) du=
  \\ &=  \cO(\|\varphi\|_{\mathcal C^5} R(N,\ve, (t-s))).
      \end{split}
\]
\end{proposition}
\begin{proof}
By Remark \ref{rem:grid} we can assume $N^2s\in\ve\bN$. Then, by Theorem \ref{thm:standard-pair-invariance}, $\phi^{N^2s}_\star\mu_N\in\cM_{\textrm{init}}$, so it is enough to prove the statement for $s=0$.
The space correlations between $q$ at time $s$ are elements of the covariance matrices $S^{(q)}_N(s) := S^{(q)}_N(z(s))$ defined in Equations \eqref{S1ts1}.
We must study the time integral of $S^{(q)}_N(s)$. To this end, it is necessary to consider the full covariance matrix:
\begin{equation}
  \label{eq:14}
\begin{split}
 & {\bs S}^{\alpha} (t)= \int_0^t  \bE_{\mu_N}\left(S^{\alpha}_N(N^2 s)\right) ds,\qquad (\alpha) \in \{ (q), (p),(q,p), (p,q)\}\\
 &{\bs S}(t)=\begin{pmatrix} {\bs S}^{(q)}(t)&{\bs S}^{(q,p)}(t)\\
 {\bs S}^{(p,q)}(t)&{\bs S}^{(p)}(t)
 \end{pmatrix}
\end{split}
\end{equation}
In this proof the time $t$ is fixed, hence we will suppress it if no confusion arises, e.g. we will write ${\bs S}$ for ${\bs S}(t)$.\\
It turns out to be convenient to compute using the Fourier transform.
We use the orthonormal base of eigenvectors of $\Delta_{\bZ_N}$: for all $j\in\bZ_N$,
\begin{equation}
  \label{eq:29}
  \begin{split}
&  \psi_j(x) = \frac{e^{i2\pi  j x/N}}{\sqrt N}, \; \psi^*_j(x) = \frac{e^{-i2\pi  j x/N}}{\sqrt N},\qquad
v_j = \omega_0^2 + 4 \sin^2\left(\frac{\pi j}{N}\right),\\
  &(\omega_0^2 - \Delta_{\bZ_N}) \psi_j= v_j \psi_j.
\end{split}
\end{equation}
In other words, the Laplacian is diagonalized by the discrete Fourier transform  $N\times N$ unitary matrix $U$ defined, for each $\phi:\bZ_N\to\bC$, by
\begin{equation}\label{eq:fourier}
\begin{split}
&\hat\phi_j:=(U \phi)_j=\sum_{x\in\bZ_N}\phi_x\psi^*_j(x)\\
&\phi_x=(U^{-1} \hat\phi)_x=\sum_{j\in\bZ_N}\hat\phi_j\psi_j(x).
\end{split}
\end{equation}
For any $N\times N$ matrix \(\cA: \bZ_N  \times \bZ_N \rightarrow \bC\),~let
\begin{equation}
\label{eq:fourier-matrix}
\begin{split}
&\hat \cA_{j,j'} := (U^{-1} \cA U)_{j,j'} = \sum_{x,x'\in\bZ_N}\cA_{x,x'}\psi_j(x)\psi^{*}_{j'}(x')\\
&\cA_{x,x'}= (U\hat{\cA}U^{-1})_{x,x'} = \sum_{j,j'\in\bZ_N}\hat{\cA}_{j,j'}\psi^{*}_j(x)\psi_{j'}(x').
\end{split}
\end{equation}
Note that if $\phi\in\cC^k(\bT,\bC)$, then, for each $N\in\bN$, it induces naturally functions from $\bZ_N\to\bC$, which we call with the same name, i.e. $\phi_x = \phi\left(\frac xN\right)$. Remark that, in the limit $N\to\infty$, we obtain   
\[
\lim_{N\to\infty}\frac 1{\sqrt N} \hat\phi_j= \int_\bT  \phi(u) e^{-i2\pi j u} du,
\]
the usual Fourier transform of $\phi$. Note also that, for any \(q \in \bZ_N\),
\begin{equation}
\label{eq:krone}
\sum_{x \in \bZ_N} \psi_{q}(x) = N^{\frac 12}\delta_{q,0}.
\end{equation}
To simplify the notation, let $\phi=\vf''$. Our objective is then to compute, for $w,w'\in\{0,1\}$,
\[
\begin{split}
\sum_{x \in \bZ_N} \phi\left(\frac{x}{N}\right)  {\bs S}^{(q)}_{x+w,x-w'}
 =  \sum_{x \in \bZ_N}\sum_{j,j',j''\in \bZ_N} \hat\phi_j \widehat {\bs S}^{(q)}_{j',j''} \psi_j(x)\psi^*_{j'}(x+w)\psi_{j''}(x-w')
\end{split}
\]
where we have set $\phi_x=\phi\left(\frac{x}{N}\right)$ for $x\in\bZ_N$. We must then compute
\begin{equation}\label{eq:S-F}
\begin{split}
&N^{-1} \sum_{x \in \bZ_N} \phi\left(\frac{x}{N}\right)  {\bs S}^{(q)}_{x+w,x-w'}=\\
&=N^{-1}\sum_{j,j',j''\in \bZ_N} \hat\phi_j \widehat {\bs S}^{(q)}_{j',j''} \sum_{x \in \bZ_N}\psi_{j-j'+j''}(x)\psi^*_{j'}(w)\psi^*_{j''}(w')\\
&=N^{-\frac 12}\sum_{j,j'\in \bZ_N} \hat\phi_{j} \widehat {\bs S}^{(q)}_{j',j'-j} \psi^*_{j'}(w)\psi^*_{j'-j}(w').
\end{split}
\end{equation}
The key step is to relate $\widehat {\bs S}^{(q)}_{j',j''}$ to ${\bs S}^{(p)}_{y,y}$. This is achieved in Proposition \ref{prop:Sqhat} which allows us to write
\begin{equation}
  \label{eq:7}
   {\hat {{\bs S}}}^{(q)}_{j,j'} = \widehat H_{j,j'} + \widehat B_{j,j'},
 \end{equation}
 where
 \begin{equation}
  \label{eq:theta}
  \begin{split}
 &\hat H_{j,j'} =N^{-\frac 12}\sum_{y\in\bZ_N} \frac{2\Theta(v_j,v_{j'})}{v_j+v_{j'}} \psi_{j-j'}(y) {{\bs S}^{(p)}_{y,y}}\\
  & \Theta(v_j, v_{j'}) = \left(1 + \frac{(v_j - v_{j'})^2}{2\gamma^2 (v_j+v_{j'})}\right)^{-1}\\
 &\sup_{l\in\bZ_N}\sum_{j\in\bZ_N}|\widehat B_{j,j-l}| \le \Const R(N,\ve,t).
\end{split}
\end{equation}
The last of \eqref{eq:theta} implies that the contribution to \eqref{eq:S-F} of $\widehat B_{j,j'}$ is negligible, indeed we have
\begin{equation}\label{eq:Bsmall}
\begin{split}
\biggl|N^{-\frac 12}&\sum_{j,j'\in\bZ_N} \hat\phi_{j} \widehat B_{j',j'-j} \psi^*_{j'}(w)\psi^*_{j'-j}(w')\biggr|
\leq \frac{\Const}{ N^{\frac 32}}\sum_{j, j'\in\bZ_N} |\hat\phi_{j}| |\widehat B_{j',j'-j}|\\
&\hskip20pt\leq \Const N^{-\frac 32}   R(N,\ve,t)\sum_{j\in\bZ_N} |\hat\phi_{j}| \leq \Const N^{-1}   R(N,\ve,t)\|\phi\|_{\cC^2},
\end{split}
\end{equation}
where, in the last inequality, we have used  Lemma \ref{lem:bound}.\\
Using Equations \eqref{eq:7}, \eqref{eq:theta} and \eqref{eq:Bsmall} into Equation \eqref{eq:S-F} yields
\begin{equation}\label{eq:S-F1}
\begin{split}
&N^{-1} \!\!\sum_{x \in \bZ_N} \!\! \phi\left(\frac{x}{N}\right)  {\bs S}^{(q)}_{x+w,x-w'}=N^{-\frac 12}\!\!\!\!\sum_{j,j' \in \bZ_N}\!\! \hat\phi_{j'} \widehat H_{j,j-j'} \psi^*_{j}(w)\psi^*_{j-j'}(w')\\
&+\cO(N^{-1}R(N,\ve,t)\|\phi\|_{\cC^2})\\
=&\frac 1N\!\!\sum_{y\in\bZ_N}
\sum_{j'\in\bZ_N}\!\! \hat\phi_{j'} \psi_{j'}(y) \psi_{j'}(w') {{\bs S}^{(p)}_{y,y}}\!\left[\sum_{j\in\bZ_N}\!\!\frac{2\Theta(v_{j},v_{j-j'})}{v_{j}+v_{j-j'}}\psi^*_{j}(w+w')\right]\\
&+\cO(N^{-1}R(N,\ve,t)\|\phi\|_{\cC^2}).
\end{split}
\end{equation}
Note that, recalling \eqref{eq:theta} and \eqref{eq:29}, for $j\in\{-N/2,\dots,N/2\}\mod N$ and since \(\Theta(v_{j},v_{j})=1\),
\begin{equation}\label{eq:theta_app}
\begin{split}
&\left|\sum_{j\in\bZ_N}\frac{2\Theta(v_{j},v_{j-j'})}{v_{j}+v_{j-j'}}\psi^*_{j}(w+w')- \sum_{j\in\bZ_N}\frac{\psi^*_{j}(w+w')}{v_{j}}\right|\leq \\
&\leq N^{-\frac 12}\sum_{j\in\bZ_N}\left|\frac{2\Theta(v_{j},v_{j-j'})}{v_{j}+v_{j-j'}}- \frac{2\Theta(v_{j},v_{j})}{v_{j}+v_{j}}\right|\leq \Const |j'| N^{-\frac 12}.
\end{split}
\end{equation}
It is then natural to define $\bG_{\omega_0,N}(x)= N^{-\frac 12}\sum_{j}\frac{\psi_{j}(x)}{v_{j}}=N^{-\frac 12}\sum_{j}\frac{\psi^*_{j}(x)}{v_{j}}$ (see also \eqref{eq:GRN}).
Which allows us to write \eqref{eq:S-F1} as
\[
\begin{split}
&N^{-1} \sum_{x \in \bZ_N} \phi\left(\frac{x}{N}\right)  {\bs S}^{(q)}_{x+w,x-w'}=N^{-1}\sum_{y\in\bZ_N} \sum_{j'\in\bZ_N} \hat\phi_{j'} \psi_{j'}(y+w')  {{\bs S}^{(p)}_{y,y}}\bG_{\omega_0,N}(w+w') \\
&+\cO\left(N^{-\frac 52} \sum_{j'\in\bZ_N} |\hat\phi_{j'} |  |j'|\sum_{y\in\bZ_N}{{\bs S}^{(p)}_{y,y}}\right)+\cO(R(N,\ve,t)\|\phi\|_{\cC^2})\\
&=N^{-1}\sum_{y\in\bZ_N} \phi\left(\frac{y+w'}N\right)  {{\bs S}^{(p)}_{y,y}}\bG_{\omega_0,N}(w+w') +\cO(R(N,\ve,t)\|\phi\|_{\cC^3})
\end{split}
\]
where we have used Lemma \ref{lem:bound}. To conclude the proof note that 
\[
\biggl|\phi\left(\frac{y+w'}N\right)-\phi\left(\frac{y}N\right)\biggr|\leq N^{-1}\|\phi\|_{\cC^1}
\]
so that, recalling \eqref{eq:riemann-approx}, \(|w+w'| \le 2\) and that $\phi=\vf''$,
\[
\begin{split}
N^{-1} \sum_{x \in \bZ_N} \vf''\left(\frac{x}{N}\right)  {\bs S}^{(q)}_{x+w,x-w'} =& N^{-1}\sum_{y\in\bZ_N} \vf''\left(\frac{y}N\right)  {{\bs S}^{(p)}_{y,y}}\bG_{\omega_0}(w+w')\\
&+ \cO(R(N,\ve,t)\|\vf\|_{\cC^5}).
\end{split}
\]
The above, recalling Equations \eqref{eq:14}, \eqref{S1ts1}, proves the proposition.
 \end{proof}

 %%%%%%%%%%
\subsection{The dynamics of the covariance matrix}\ \\
\label{sec:dynam-covar-matr}
{ Since the final time $t$ is fixed in this section, we remove the dependence of $\bs S$ on $t$. Before starting our study of the covariance matrix, we need to introduce some notation.}\\
In analogy to \eqref{eq:unperturbed-matrices}, we define the matrix $2N\times 2N$ 
\begin{equation}
\label{A}
\boldsymbol{A}=
\left(
  \begin{array}{cc}
    \zerob&-\Id\\
    -\Delta_{\bZ_N} +\om_0^2 \Id& \ga \Id
  \end{array}
\right),
\end{equation}
Also, we define
\begin{equation}
\label{S2}
\Xi (S^{(p)}) =\left[\begin{array}{cc}
\zerob&\zerob\\
\zerob&{\rm Diag}(S^{(p)})
\end{array}\right],
\end{equation}
where $\zerob$ is $N\times N$ null matrix and for each $N\times N$ matrix $M$
 \begin{equation}
\label{D2}
{\rm Diag}(M) =
       \begin{bmatrix}
 M_{1,1}& 0 & 0 &\dots&0\\
                     0&   M_{2,2} &  0 &\dots&0\\
                     0 & 0 &  M_{3,3} &\dots&0\\
                     \vdots & \vdots & \vdots & \vdots&\vdots\\
 0& 0 & 0 & \dots &  M_{N,N}
                        \end{bmatrix}.
                      \end{equation}
The following result controls the Fourier transform of the covariance matrix.
\begin{proposition}\label{prop:Sqhat}
We have the representation\footnote{ See Equation \eqref{eq:theta} for the definition of $\Theta(v_j,v_{j'})$, \eqref{eq:29} for the definition of $v_j$, and Equations \eqref{eq:29}, \eqref{eq:fourier} and \eqref{eq:fourier-matrix} for the Fourier coefficients.}
\[
   {\hat {{\bs S}}}^{(q)}_{j,j'} = \frac{2\Theta(v_j,v_{j'})}{v_j+v_{j'}}  \hat F_{j,j'} + \widehat B_{j,j'}=:\hat H_{j,j'} + \widehat B_{j,j'}
\]
where 
\[
\hat F_{j,j'} = \sum_{y\in\bZ_N} \psi_j(y) \psi^*_{j'}(y) {{\bs S}^{(p)}_{y,y}}
\]
and the matrix $\widehat B$ satisfies
 \begin{equation}
   \label{eq:bB}
   \sup_{l\in\bZ_N}\sum_{j\in\bZ_N}|\widehat B_{j,j-l}| \le \Const R(N,\ve,t).
 \end{equation}
 \end{proposition}
 \begin{proof}
 By a straightforward computation:
\begin{equation}
  \label{eq:15}
\cG S_N(t) = -\boldsymbol{A} S_N(t) - S_N(t) \boldsymbol{A}^T + 2\gamma \Xi(S^{(p)}_N(t)).
\end{equation}
It is then convenient to define, recall also definition \eqref{eq:14},
\begin{equation}\label{eq:R_def}
\cR:=\boldsymbol{A} {\bs S} +{\bs S} \boldsymbol{A}^T -
      2 \gamma \Xi \Big({\bs S}^{(p)}\Big),
\end{equation}
which, remembering \eqref{S1ts}, can be written in components as
     \begin{equation}
       \label{eq:18}
       \begin{split}
         & \cR^{(q)} = -{\bs S}^{(p,q)} - {\bs S}^{(q,p)} \\
         &\cR^{(q,p)}={\bs S}^{(q)} (\omega_0^2 \Id - \Delta_{\bZ_N}) + \gamma {\bs S}^{(q,p)} - {\bs S}^{(p)}
         \\
         &\cR^{(p,q)}=(\omega_0^2 \Id - \Delta_{\bZ_N})  {\bs S}^{(q)} + \gamma {\bs S}^{(p,q)} - {\bs S}^{(p)}
         \\
        &  \cR^{(p)}=(\omega_0^2 \Id - \Delta_{\bZ_N})  {\bs S}^{(q,p)}+  {\bs S}^{(p,q)} (\omega_0^2 \Id - \Delta_{\bZ_N})
          \\
         &\phantom{\cR^{(p)}=}
         + 2\gamma {\bs S}^{(p)} -  2 \gamma{\rm Diag}\Big({\bs S}^{(p)}\Big).
       \end{split}
     \end{equation}
Using the notation \eqref{eq:fourier-matrix}, we can rewrite Equations \eqref{eq:18} as
\begin{equation}
  \label{eq:31}
  \begin{split}
  &\widehat \cR_{j,j'}^{(q)} = - {\hat {{\bs S}}}^{(q,p)}_{j,j'} - {\hat {{\bs S}}}^{(p,q)}_{j,j'}\\
  &\widehat \cR_{j,j'}^{(q,p)} +\widehat \cR_{j,j'}^{(p,q)} + \gamma\widehat \cR_{j,j'}^{(q)} = 
  \left(\upsilon_j +\upsilon_{j'}\right) {\hat {{\bs S}}}^{(q)}_{j,j'}
    - 2  {\hat {{\bs S}}}^{(p)}_{j,j'}\\
 &\widehat \cR_{j,j'}^{(q,p)} - \widehat \cR_{j,j'}^{(p,q)}-\gamma \widehat \cR_{j,j'}^{(q)}
 = 2 \gamma  \hat {{\bs S}}^{(q,p)}_{j,j'} -  \hat {{\bs S}}^{(q)}_{j,j'}\left(\upsilon_j  - \upsilon_{j'}\right)\\
  & \widehat \cR_{j,j'}^{(p)} +\upsilon_{j'}\widehat \cR_{j,j'}^{(q)}  =
    \left(\upsilon_j  - \upsilon_{j'}\right) 
    \hat {{\bs S}}^{(q,p)}_{j,j'} -2\gamma \hat F_{j,j'}
    + 2\gamma  {\hat {{\bs S}}}^{(p)}_{j,j'}.
  \end{split}
\end{equation}
In fact, the second line of Equation \eqref{eq:31} follows from the first three lines of \eqref{eq:18} and the fact that, for any $N\times N$ matrix $\mathcal A$,
\begin{equation}
\label{eq:diagonalization}
\begin{split}
&\left(U^{-1}(\omega_0^2\Id-\Delta_{\bZ_N})\mathcal A U\right)_{j,j'} = \upsilon_j \hat{\mathcal A}_{j,j'}\\ 
&\left(U^{-1}\mathcal A (\omega_0^2\Id-\Delta_{\bZ_N}) U\right)_{j,j'} = \upsilon_{j'} \hat{\mathcal A}_{j,j'}.
\end{split}
\end{equation}
 The third line of Equation \eqref{eq:31} is obtained similarly by the first three lines of \eqref{eq:18}, and the fourth equation follows from the fourth line of \eqref{eq:18} and the first of \eqref{eq:31}.\\
The next Lemma, whose proof is briefly postponed, shows that the matrices ${\widehat\cR}^{\alpha}$ are negligible in the limit.

\begin{lemma}\label{lem:Restimate}
For any $\alpha \in \{(q), (p), (q,p), (p,q)\}$ we have the estimates
\[
  \sup_{l\in\bZ_N}\sum_{j\in\bZ_N}|\widehat\cR^{\alpha}_{j,j-l}| \le \Const R(N,\ve,t).
\]
\end{lemma}
By simple algebraic manipulations, Equation \eqref{eq:31} implies
\[
{\hat {{\bs S}}}^{(q)}_{j,j'} = \frac{2\Theta(v_j,v_{j'})}{v_j+v_{j'}}  \hat F_{j,j'} + \widehat B_{j,j'},
\]
where the matrix $\widehat B$ is a linear combination of the $\widehat \cR^\alpha$ with uniformly bounded coefficients. Hence, the bound \eqref{eq:bB} follows then from Lemma \ref{lem:Restimate}, concluding the proof.
\end{proof}

\begin{proof}[\it \bfseries Proof of Lemma \ref{lem:Restimate}]
We need to treat appropriate linear combinations of the elements of ${\bs S}$, to this end we consider $2N\times 2N$ constant real matrices 
  \begin{equation}
\bs\beta_N=\begin{pmatrix} \beta^{(q)}_N&\beta_N^{(q,p)}\\ 
\beta_N^{(p,q)}&\beta_N^{(p)},
\end{pmatrix}\label{eq:17}
\end{equation}
where $\beta_N^{\alpha}$ are $N\times N$ matrices such that \(\sup_{N}\|\bol \beta_N\| < \infty\). By \eqref{eq:norm},
\begin{equation}\label{eq:en_control}
\trn{\operatorname{Tr}(S_N{\bs \beta_N})}\leq \Const  \|{\bs \beta_N}\|.
\end{equation}

Integrating in time and averaging with respect to the initial conditions Equation \eqref{eq:15}, using Corollary \ref{cor:accellerated_generator} with $A(t) =  \operatorname{Tr}(S_N(t){\bs \beta_N})$, and Equations \eqref{eq:en_control}, \eqref{eq:R_def}, we obtain,
     \begin{equation}
       \label{eq:23}
 \begin{split}
      & \operatorname{Tr}\left({\bs\beta}_N\cR\right)=\operatorname{Tr}\left({\bs\beta}_N\left[\boldsymbol{A} {\bs S} 
      +{\bs S}\boldsymbol{A}^T -
      2 \gamma \Xi \Big({\bs S}^{(p)}\Big)\right]\right) \\
      &= \frac{\operatorname{Tr}\left({\bs\beta}_N\bE_{\mu_N}\left(S_N(0)-S_N(N^2t)\right)\right)}{N^2}
     +\cO\left(\|{\bs \beta}_N\| R(N,\ve,t)\right).
\end{split}
     \end{equation}
{ Denote
\begin{equation}\label{eq:symmetric}
\begin{split}
    \cR^{\alpha,s}_{x,y} &=  \frac 12 \left(\cR^{\alpha}_{x,y} + \cR^{\alpha}_{-x,-y}\right)\\
    \widehat \cR^{\alpha,s}_{l,j}&=\frac 12 \left(\widehat\cR^{\alpha}_{l,j}+\widehat\cR^{\alpha}_{-l,-j}\right)=\Re \mathfrak{e}(\widehat\cR^{\alpha}_{l,j}).
\end{split}
\end{equation}
To deal with the absolute value of the Fourier coefficients, we must make the appropriate choices of ${\bs \beta}_N$. For any fixed $l\in\bZ$ and $\alpha \in \{(q), (p), (q,p), (p,q)\}$, we choose ${\bs \beta}_N={\bs \beta}^{l ,\alpha}$, with the following
block components
\[
  \beta^{l,\alpha,\alpha}_{x,y}=\sum_{j \in \bZ_N}\psi^*_{-j}(x)\psi_{-j+l}(y)\operatorname{sign}(\widehat \cR^{\alpha,s}_{l-j,-j}),
  \quad x,y \in \bZ_N,
\]
for the $\alpha$ block and $\beta^{l,\alpha,\alpha'}=0$ for \(\alpha ' \neq \alpha\).} Note that, for all $w,v\in\bR^N$,
\[
\begin{split}
  |\langle w,\beta^{l,\alpha,\alpha} v\rangle|& =
\left|\sum_{j\in\bZ_N}\sum_{x,y\in\bZ_N}\psi^*_{-j}(x)\psi_{-j+l}(y)w_xv_y
    \operatorname{sign}(\widehat {\cR}^{\alpha,s}_{l-j,-j})\right|\\
&\leq \sum_{j\in\bZ_N}|\hat w_{-j}|\,|\hat v_{j-l}|\leq \|\hat w\|\|\hat v\|=\|w\|\|v\|,
\end{split}
\]
which implies that $\|\bs\beta^{l,\alpha}\|\leq 1$. Moreover, for any \(\cR^{\alpha}\) in \eqref{eq:18},
\[
\begin{split}
 \operatorname{Tr}\left(\cR^{\alpha,s} { \beta^{l,\alpha, \alpha}}\right) &= \!\!\!
 \sum_{x,y \in \bZ_N} \sum_{j,j',j''\in \bZ_N} \!\!\!\!\!\psi^{*}_{j'}(x)  \psi_{j''}(y) \widehat \cR^{\alpha,s}_{j',j''}\psi^*_{-j}(y)  \psi_{-j+l}(x)  \operatorname{sign}(\widehat \cR^{\alpha,s}_{l-j,-j})
 \\
&= \frac 1N \sum_{j,j',j''\in \bZ_N}\sum_{x,y\in \bZ_N}  \psi^*_{j+j'-l}(x)  \psi_{j+j''}(y)  \operatorname{sign}(\widehat \cR^{\alpha,s}_{l-j,-j})
\widehat \cR^{\alpha,s}_{j',j''}\\
&= \sum_{j,j',j''\in \bZ_N} \delta_{j,l-j'}  \delta_{j,-j''}  \operatorname{sign}(\widehat \cR^{\alpha,s}_{l-j,-j})
\widehat \cR^{\alpha,s}_{j',j''}\\
&=  \sum_{j\in \bZ_N}  \operatorname{sign}(\widehat \cR^{\alpha,s}_{l-j,-j}) \widehat \cR^{\alpha,s}_{l-j,-j}
 = \sum_{j\in\bZ_N}|\widehat\cR^{\alpha,s}_{j,j-l}|, 
\end{split}
\]
where in the line next to the last we have used \eqref{eq:krone}.\\
Denote $\check \cR^{\alpha}_{x,y} =\cR^{\alpha}_{-x,-y}$, $\check \beta^{l,\alpha,\alpha}_{x,y} =\beta^{l,\alpha,\alpha}_{-x,-y}$ and $\check {\bs\beta}^{l,\alpha}$ the matrix with all zero blocks apart from $\check \beta^{l,\alpha,\alpha}$. Then, by the above equality and \eqref{eq:symmetric} we obtain
\begin{equation}
  \label{eq:19}
  \begin{split}
\sum_{j\in\bZ_N}|\widehat\cR^{\alpha,s}_{j,j-l}|  &= \operatorname{Tr}\left(\cR^{\alpha,s} { \beta^{l,\alpha, \alpha}}\right) \\
&= \frac 12 \operatorname{Tr}\left(\cR^{\alpha} { \beta^{l,\alpha,\alpha}}\right)
  + \frac 12 \operatorname{Tr}\left(\check\cR^{\alpha} { \beta^{l,\alpha,\alpha}}\right)\\
 &=  \frac 12 \operatorname{Tr}\left(\cR^{\alpha} { \beta^{l,\alpha,\alpha}}\right)
  + \frac 12 \operatorname{Tr}\left(\cR^{\alpha} { \check\beta^{l,\alpha,\alpha}}\right)\\
  &= \frac 12 \operatorname{Tr}\left(\cR { \bs\beta^{l,\alpha}}\right)
  + \frac 12 \operatorname{Tr}\left(\cR {\check{\bs\beta}}^{l,\alpha}\right).
\end{split}
\end{equation}
By \eqref{eq:23}, and \(\|\bs \beta^{l,\alpha}\| \le 1\), \(\|\check{\bs \beta}^{l,\alpha}\| \le 1\), the right-hand side is less than
\[
\begin{split}
  \sup_{l\in\bZ_N} &\frac{\operatorname{Tr}\left(\frac 12 
  ({\bs\beta}^{l,\alpha} +\check{\bs\beta}^{l,\alpha})
    \bE_{\mu_N}\left(S_N(0)-S_N(N^2t)\right)\right)}{N^2} +\cO\left(R(N,\ve,t)\right).
\end{split}
\]
Note that, by  \eqref{eq:en_control}, for any \(t \in \bR_+\), 
\[
\begin{split}
  \sup_{l\in\bZ_N}\left| \operatorname{Tr}\left({\bs\beta_N^{l,\alpha}}\bE_{\mu_N}\left(S_N(t))\right)\right)\right|&
  =\sup_{l\in\bZ_N}\left| \bE_{\mu_N}\left(\operatorname{Tr}\left({\bs\beta_N^{l,\alpha}}S_N(t))\right)\right)\right|
\leq \Const  N.
\end{split}
\]
The Lemma follows by estimating the imaginary part in a similar manner.
\end{proof}

%%%%%
\subsection{Equipartition}\ \\
\label{sec:equipartition}
In this section, we present an equipartition result that connects the expectation of $p_x^2$ to that of $\fe_x$.

\begin{proposition}
\label{cor012912-21} \relax
For any $\varphi\in \cC^3(\bT)$, $N\in\bN$, $\ve\in (0,N^{-6})$, \(t>s\in \bR^+\), and initial probability measure $\mu_N\in\cM_{\textrm{init}}$ we have, 
\[
\begin{split}
\Biggl|\frac{1}{N}\sum_{x \in \bZ_N}\varphi''\left(\frac{x}{N}\right) 
\int_{s}^{t} \bE_{\mu_N}\Big[p_x^2(N^2u)- \fe_x(N^2u)\Big] &du\Biggr|  \leq\Const \|\varphi\|_{\cC^3}
 R(N,\ve,(t-s)).
\end{split}
\]
  \end{proposition}
\begin{proof}
By Remark \ref{rem:grid} we can assume $N^2s\in\ve\bN$. Then, by Theorem \ref{thm:standard-pair-invariance}, $\phi^{N^2s}_\star\mu_N\in\cM_{\textrm{init}}$, so it is enough to prove the statement for $s=0$. A direct computation yields, setting $Q_x=q_{x-1} ( q_{x} -  q_{x-1})$,
\begin{equation}
  \label{eq:4}
  \begin{split}
   2p_x^2 - 2\fe_x =&
  -(\nabla_{\bZ_N}Q)_x + \cG \left( q_x p_x + \frac{\gamma}{2}  q_x^2\right).
\end{split}
\end{equation}
Since 
\[
\trn{\frac{1}{N}\sum_{x\in\bZ_N} \varphi''\left(\frac{x}{N}\right)\left[ q_x  p_x + \frac{\gamma}{2}  q_x^2\right]}\leq \Const\|\vf\|_{\cC^2},
\]
we can apply again Corollary \ref{cor:accellerated_generator} and write
\begin{equation}\label{010905-22}
\begin{split}
& \frac{2}{N}\sum_{x\in\bZ_N} \varphi''\left(\frac{x}{N}\right)\int_0^t \mu_N \Big[ p_x^2(N^2s) - \fe_x(N^2 s)\Big] d s\\
&=\frac{1}{N}\sum_{x\in\bZ_N} \nabla_{\bZ_N}\varphi''\left(\frac{x}{N}\right) \int_0^t \bE_{\mu_N}\left[ Q_{x}(N^2s) \right] ds\\
 & \phantom{=}
 +\frac{1}{N^3}\sum_{x\in\bZ_N} \varphi''\left(\frac{x}{N}\right) \bE_{\mu_N}  \Big[ q_x(N^2t) p_x(N^2t) + \frac{\ga}{2}  q_x^2(N^2t) \Big]\\
&  \phantom{=}  - \frac{1}{N^3}\sum_{x=0}^n \varphi''\left(\frac{x}{N}\right) \bE_{\mu_N} \Big[ q_x(0)  p_x(0) +\frac{\ga}{2}  q_x^2(0) \Big]\\
&  \phantom{=}
+\|\vf''\|_\infty\cO\left(R(N,\ve,t)\right).
\end{split}
\end{equation}
We have, for each $s\in\bR_+$,
\[
\left|\frac{1}{N^3}\sum_{x\in\bZ_N} \varphi''\left(\frac{x}{N}\right) \bE_{\mu_N}
\Big[ q_x(N^2s) p_x(N^2s) + \frac{\ga}{2}  q_x^2(N^2s) \Big]\right| \le\Const
\frac{\|\varphi\|_{\cC^2}} {N^2}.
 \]
 While for the term in the second line of \eqref{010905-22} we have
 \begin{equation}
   \label{eq:53}
   \begin{split}
  &  \left|\int_0^{t}   \frac{1}{N}\sum_{x\in\bZ_N} \nabla_{\bZ_N} \varphi''\left(\frac{x}{N}\right)
    \bE_{\mu_N}\left[ Q_{x}(N^2s))\right] ds\right|  \le \Const\frac{\|\varphi\|_{\cC^3} t}{N}.
   \end{split}
 \end{equation}
 \end{proof}

%%%%%%%%%%%
\section{Proof of Main Theorem \ref{thm:heat-equation}}\label{sec:final}
This section is devoted to the proof of Theorem \ref{thm:heat-equation}.\\
Given a sequence $(\ve_N)_{n\in\bN}$ such that $\ve_N\leq N^{-\alpha}$, $\alpha>6$,
Theorem \ref{prop:fundamental} implies that, for all $\vf\in\cC^6$,  \(t \in \bR_{+}\),
\begin{equation}
\label{eq:54}
\lim_{N\to\infty}  \left[\xi_N(t,\varphi) -  \xi_N(0,\varphi) -
\frac{D}{2\gamma} \int_0^t  \xi_N(s,\varphi'') ds \right] = 0.
\end{equation}

For each $t_*>0$, let  \(\cC^0_{*} := \cC^0([0,t_{*}], \mathcal P(\bT))\) where \(\mathcal P(\bT)\) is endowed with the metric defined in \eqref{eq:33}, 
which metrizes the weak topology. 
Note that, for all $N\in\bN$, $\xi_N\in \cC^0_{*}$. Let \((\mu_N)_{N\in\bN} \in \cM^{\star}_{\textrm{init}}\), \(\Ene_0\) be given by \eqref{eq:initial-conditions} and \(\Ene(\cdot) \in \cC^0_{*}\) be solution of the heat equation with initial condition \(\Ene_0\).\\
Since $\mathcal P(\bT)$ is compact, the relative compactness of the sequence $(\xi_N)$ follows from the
equicontinuity in time of $\xi_N$ proved in Appendix \ref{sec:compactness} and by Ascoli's theorem \cite[ pag. 223, Theorem 17]{Kelley}.
Consequently, there exist accumulation points $\xi_*(\cdot) \in \cC^0_{*}$ of the sequence $(\xi_N(\cdot))$. By \eqref{eq:54} such accumulation points satisfy the weak heat Equation \eqref{eq:PDE}. Furthermore, by \eqref{eq:initial-conditions}, \(\lim_{N \rightarrow \infty}\xi_N(0) = \Ene_0\). Therefore, since the heat equation has a unique solution,\footnote{ Just compute its Fourier coefficients $\Ene_k(t):=\int_{\bT} e^{2\pi i kx} \Ene(t, dx)$.} one has
\begin{equation}\label{eq:limit}
\lim_{N\to\infty}\xi_{N}(t)=\Ene(t)
\end{equation}
uniformly in $[0,t_*]$ in the weak topology.\\
To extend the result in time, note that 
$$
(\bar q(t), \bar p(t),\bar \theta(t)) :=
(q(t+N^2 t_*), p(t+N^2 t_*),\theta(t+N^2 t_*))
$$ 
satisfies Equations \eqref{eq:det}. Moreover, the distribution $\mu^*_N$ of 
$(\bar q(0), \bar p(0),\bar \theta(0))$ belongs to
$\cM_{\textrm{init}}$ by Theorem \ref{thm:standard-pair-invariance}, provided $N^2t_*\in\ve\bN$. By Remark \ref{rem:grid}, this can be assumed by replacing $t_*$ with $N^{-2}\ve\lceil N^2t_*\ve^{-1}\rceil$, at the price of a negligible error.
In addition, we have just seen that $\mu_N^*$ satisfies \eqref{eq:initial-conditions}, that is, its energy profile converges to $\Ene(t_*)$. It follows that we can apply our results to such a process and extend it for another time $t_*$. Iterating the argument yields
\begin{equation}\label{eq:limit-infy}
\lim_{N\to\infty}\xi_{N}(t)=\Ene(t)\qquad\forall t\in\bR_+.
\end{equation}
%%%%%%%%%%%%%%%%%%%%%%%%%%

\subsection{Error estimate}\label{sec:atlast}\ \\
To conclude the proof of Theorem \ref{thm:heat-equation}, we must compute the difference
\[
\eta_N(t,\varphi)=\xi_N(t,\varphi) - \Ene(t,\varphi)\;,\quad  t \in \bR_+.
\]
Theorem \ref{prop:fundamental} yields
\begin{equation}\label{eq:main_remind}
\begin{split}
 \eta_N(t,\varphi)  = \eta_N(s,\varphi)+\frac{D}{2\gamma} \int_s^t  \eta_N(u,\varphi'') du 
+ \cO\left(\|\vf\|_{\cC^5} R(N,\ve,t-s)
\right).
\end{split}
\end{equation}
 Evaluating \eqref{eq:main_remind} on $\bs\psi_k(u)=e^{2i\pi k u}$, $k\neq 0$, we obtain\footnote{Note that, by linearity, Equation \eqref{eq:main_remind} holds as well for complex functions.}
\begin{equation}\label{eq:error_eq}
\begin{split}
 \eta_N(t,\bs\psi_k)  = \eta_N(s,\bs\psi_k)- &2\pi^2k^2\frac{D}{\gamma} \int_s^t  \eta_N(\tau,\bs\psi_k) d\tau\\
  &+\cO\left(k^5 R(N,\ve,t-s)\right),
\end{split}
\end{equation}
while $\eta_N(t,1)=0$ for all times.\\
Since \(\ve \in (0, N^{-\alpha})\), \(\alpha > 6\), and, from Equation \eqref{eq:err}, $ \beta_{\star} = \min\left\{1,\frac{\alpha}6-1\right\}$, one has, for  \(N \ge 2\),
\begin{equation}
\label{eq:epsilontoN} 
R(N,\ve,t) \le \Const \max \left\{t,1\right\}  N^{-\beta_*}\left(\ln N\right)^2.
\end{equation} 
Iterating once \eqref{eq:error_eq} for $t-s\leq\min \left\{ 1, \frac{\gamma}{2\pi^2 D}k^{-2}\right\}$, yields
\begin{equation}\label{eq:error_eq2}
\begin{split}
\left| \eta_N(t,\bs\psi_k) \right| \leq& \left(1- 2\pi^2k^2\frac{D}{\gamma}(t-s)\right)|\eta_N(s,\bs\psi_k)| \\
&+4\pi^4k^4\frac{D^2}{\gamma^2}\int_s^t (t-\tau) |\eta_N(\tau,\bs\psi_k)| d\tau
\\
  &+\Const (k^7(t-s)+k^5) N^{-\beta_{\star}}(\ln N)^2.
\end{split}
\end{equation}
Let $\zeta(s,u, k)=\sup_{\tau\in [s,s+u]}|\eta_N(\tau,\bs\psi_k)|$, and $t_k= \min \left\{1, \frac{\gamma}{4\pi^2 D}k^{-2}\right\}$,  and $a_k=2\pi^2k^2\frac D\gamma$, so that $a_kt_k\le\frac 12$. By \eqref{eq:error_eq}, for all $\tau\in[s,s+t_k]$ we have $|\eta_N(\tau,\bs\psi_k)|\le|\eta_N(s,\bs\psi_k)|+a_kt_k\,\zeta(s,t_k,k)+\Const k^5N^{-\beta_\star}(\ln N)^2$, hence
\[
\zeta(s, t_k, k) \leq 2|\eta_N(s,\bs\psi_k)| +\Const k^5 N^{-\beta_{\star}}(\ln N)^2.
\]
Using the above in \eqref{eq:error_eq2} with $t=s+t_k$, and since $x^2\le \frac x2$ for $x=a_kt_k$, yields
\begin{equation}\label{eq:error_eq3new}
\begin{split}
|\eta_N(s+t_k,\bs\psi_k)|&\leq \left(1-a_kt_k+a_k^2t_k^2\right)|\eta_N(s,\bs\psi_k)|+\Const k^5 N^{-\beta_{\star}}(\ln N)^2\\
&\leq \left(1-\frac{a_kt_k}2\right)|\eta_N(s,\bs\psi_k)|+\Const k^5 N^{-\beta_{\star}}(\ln N)^2.
\end{split}
\end{equation}
 For each $t\in\bR_+$, we write $t=mt_k+\tau$, $m\in\bN$, $\tau\in [0,t_k)$. Note that $a_kt_k\ge\min\{2\pi^2\frac D\gamma,\frac 12\}$, so the factor $(1-\frac{a_kt_k}2)$ in \eqref{eq:error_eq3new} is bounded away from $1$ uniformly in $k$, and $1-\frac{a_kt_k}2\le e^{-\frac{a_kt_k}2}$. Iterating \eqref{eq:error_eq3new} $m$ times, then using the bound on $\zeta(mt_k,t_k,k)$ and $a_kt_k\le\frac 12$, and setting $\tilde D=\frac{\pi^2D}{\gamma}$ (so that $\frac{a_k}2=\tilde Dk^2$), we obtain
\begin{equation}\label{eq:error_expdec}
|\eta_N(t,\bs\psi_k)|\leq  2e^{\frac 14}e^{-\tilde D k^2 t}|\eta_N(0,\bs\psi_k)|+\Const k^5 N^{-\beta_{\star}}(\ln N)^2.
\end{equation}
Recall that $\left|\eta_N(0,\bs\psi_k)\right|=\varpi (N,k, \mu_N)$.\footnote{  See \eqref{eq:varpi} for the definition of \(\varpi\).}
Thus, for any $\vf\in\cC^{7}(\bT, \bR)$ and $L\in\bN$,\footnote{ Setting \(\hat \vf (k) = \int_{\bT} e^{-2\pi i k \theta}\vf(\theta)d\theta\) for the Fourier coefficients.}
\begin{equation}\label{eq:final}
\begin{split}
&\sup_{t \in \bR_+}|\xi_N(t,\vf)-\Ene(t,\vf)|\leq \sup_{t \in \bR_+} \sum_{k\in\bZ\setminus \{0\}} |\hat\vf(k)||\eta_N(t,\bs\psi_k)|\\
&\leq   \Const \|\vf\|_{\cC^{7}}\left\{\sum_{k\in\bZ\setminus\{0\}}k^{-7}e^{-\tilde D t}\varpi (N,k, \mu_N)+  N^{-\beta_{\star}}(\ln N)^2\right\}\\
&\leq \Const \|\vf\|_{\cC^{7}}\left\{e^{-\tilde D t}\Psi (N; \mu_N)+  N^{-\beta_{\star}}(\ln N)^2\right\}
\end{split}
\end{equation}
which proves the theorem. 

%%%%%%%%%%%%
\newpage

\appendix

%%%%%%%%%
\section{Decay of correlations for expanding map}\label{sec:TO}
In this section, we quickly recap the facts used in this paper pertaining to the statistical properties of smooth expanding maps. For simplicity, we confine ourselves to one dimensional maps, but the same theory holds in higher dimensions and for hyperbolic maps, see \cite{DKL21} for a quick exposition or \cite{Bal00, Bal18} for exhaustive details.\\
Expanding maps have lots of invariant probability measures, but they have the remarkable property of having only one absolutely continuous invariant probability measure (see Remark \ref{rem:uniqueabs}). Moreover, correlations with respect to such a measure decay exponentially if tested on smooth functions.
Such properties follow from Theorem \ref{thm:decay} and its consequence Lemma \ref{lem:dec_cor}.\\
Let $f\in \cC^2(\bT,\bT)$ such that $f'\geq \lambda>1$. A trivial consequence of expansivity is the following.
\begin{lemma}
\label{lem:dens-preimages}
For any \(\theta \in \bT\), \((f^{-n}(\theta) )_{n \in \bN}\) is dense.
\end{lemma}
\begin{proof}
Let \((a, b) \subset \bT\). Note that \(f^n\), \(n \in \bN\), belongs to \(\cC^2(\bT,\bT)\), and \(\inf_{[0,1)}(f^n)'  \ge \lambda^n\). Therefore, for any \(n \in \bN\), either \(f^n((a,b)) = \bT\) or
\[
1>\left|f^n\left((a,b)\right)\right| = \int_a^b (f^n)'(\theta)d\theta \geq |b-a|\lambda^n,
\]
which can happen only for finitely many \(n\). Hence, there exists $\bar n\in\bN$ such that \(f^{\bar n}((a,b))= \bT\), so \(f^{-\bar n}(\theta ) \cap (a,b) \neq \emptyset\), thus proving the Lemma.
\end{proof}
We are interested in studying the decay of correlation for observables of bounded variation. To this end, it is useful to introduce the concept of transfer operator defined as, for each $g\in\cC^0, \rho\in L^1$,\footnote{ By $g\circ f(\theta)$ we mean, as usual, the composition $g(f(\theta))$.}
\begin{equation}\label{eq:L_deriv}
\int_{\bT} \rho(\theta)  (g\circ f)(\theta)d\theta=\int_{\bT} g(\theta) \cL \rho(\theta)d\theta.
\end{equation}
Note that, by definition, $\int_{\bT} \cL \rho(\theta)d\theta=\int_{\bT}\rho(\theta)d\theta$ for each $\rho\in L^1$.
A simple change of variables shows that
\[
\cL \rho(\theta)=\sum_{\theta'\in f^{-1}(\theta)}\frac{\rho(\theta')}{f'(\theta ')}
\]
this is usually called the SRB (Sinai, Ruelle, Bowen) transfer operator.
Note also that $|\cL \rho(\theta)|\le \cL |\rho|(\theta)$, that implies $\|\cL \rho\|_{L^1} \le  \|\rho\|_{L^1}$. All the properties used in this paper are summarized as follows.  We recall that the space BV of functions from \(\bT\) to \(\bC\) of bounded variation is defined as
\[
BV = \biggl\{ g \in L^{1}(d\theta): \sup_{\{\vf \in \cC^1(\bT, \bC): \text{ }\|\vf\|_{\cC^0} \leq 1\}} \int_{\bT} \vf'(\theta) g(\theta) < + \infty \biggr\},
\]
and, for \(g \in\) BV, the BV norm is defined as
\[
\|g\|_{BV} = \|g\|_{L^1} + \sup_{\{\vf \in \cC^1(\bT, \bC): \text{ }\|\vf\|_{\cC^0} \leq 1\}} \int_{\bT} \vf'(\theta) g(\theta).
\]
\begin{theorem} \label{thm:decay}
$\cL$ is a bounded operator when acting on BV (we write \(\cL\in L(BV, BV)\)). Moreover, there exists $\nu\in (0,1)$ such that the spectrum 
\[
\sigma_{BV}(\cL)\subset \{1\}\cup\{z\in \bC\;:\; |z|< \nu\},
\]
and $1$ is a simple eigenvalue of $\cL$ with eigenvector $\rho_\star\in\cC^1$, $\rho_\star>0$.
\end{theorem}

The above theorem, together with Equation \eqref{eq:L_deriv}, implies that there exists a unique $\rho_\star\in BV$ such that $\cL \rho_\star=\rho_\star$,  \(\int \rho_{\star} = 1\), and that $\cL^n$, $n\in\bN$, has the spectral decomposition
\begin{equation}\label{eq:spectral}
\cL^n \rho=\rho_\star\int_{\bT} \rho(\theta)d\theta+Q^n\rho
\end{equation}
where $Q\in L(BV,BV)$ has spectral radius strictly smaller than $\nu$, hence 
\begin{equation}
\label{eq:spectral-gap-Q}
\|Q^n\|_{BV}\leq \Const \nu^n.
\end{equation}
The above decomposition allows us to establish the following basic result on the decay of correlations.
\begin{lemma}\label{lem:dec_cor}
For each $\rho\in BV$, $g\in L^1$ and $n\in\bN$, we have
\[
\begin{split}
\left|\int_{\bT} \rho(\theta) (g\circ f^n)(\theta)d\theta-\int_{\bT} g(\theta) \rho_\star(\theta)d\theta\int_{\bT} \rho(\theta')d\theta'\right|\leq \Const \nu^n \|g\|_{L^1} \|\rho\|_{BV}.
\end{split}
\]
Moreover, if $\rho, g\in L^2$, then
\[
\begin{split}
\lim_{n\to\infty}\left|\int_{\bT} \rho(\theta) (g\circ f^n)(\theta)d\theta-\int_{\bT} g(\theta) \rho_\star(\theta)d\theta\int_{\bT} \rho(\theta')d\theta'\right|=0.
\end{split}
\]
\end{lemma}
\begin{proof}
Let \(\rho \in BV\) and \(g \in L^1\). We just compute, recalling \eqref{eq:spectral},
\[
\begin{split}
\biggl|\int_{\bT} \rho(\theta)(g\circ f^n)&(\theta)d\theta-\int_{\bT} g(\theta) \rho_\star(\theta)d\theta\int_{\bT} \rho(\theta)d\theta\biggr|\\
&=\left|\int_{\bT} g(\theta)\cL^n \rho(\theta)d\theta-\int_{\bT} g(\theta) \rho_\star(\theta)d\theta\int_{\bT} \rho(\theta)d\theta\right|\\
&=\left|\int_{\bT} g(\theta) Q^n \rho(\theta)d\theta\right|\leq \Const \nu^n \|g\|_{L^1} \|\rho\|_{BV}.
\end{split}
\]
To prove the second statement, first note that, for any \(n \in \bN\), \(\rho \in L^2\),
\[
    \|\cL^n \rho\|^2_{L^2} = \int_{\bT}  \rho (\theta) \left(\cL^n \rho\right) \circ f^n (\theta)d\theta \le \|(\cL^n \rho) \circ f^n\|_{L^2}\|\rho\|_{L^2},
\]
and, since \eqref{eq:spectral} and \eqref{eq:spectral-gap-Q} imply
\[
\|\cL^n 1\|_{\infty} \le  \|\cL^n 1 - \rho_{\star}\|_{BV} + \|\rho_{\star}\|_{BV} \!= \! \|Q^n 1\|_{BV} + \|\rho_{\star}\|_{BV}\!\le\! \Const,
\]
we have
\[
  \|(\cL^n \rho) \circ f^n\|_{L^2}^2 = \int_{\bT}|\cL^n \rho|^2  \cL^n 1 \le \|\cL^n 1\|_{\infty} \|\cL^n \rho\|^2_{L^2}\leq \Const \|\cL^n \rho\|^2_{L^2}.
\]
These above inequalities imply that \(\cL\) is power-bounded on \(L^2\), i.e.,
\begin{equation}
\label{eq:power-bounded}
  \|\cL^n \rho\|_{L^2} \le \Const \|\rho\|_{L^2}, \quad \text{for any } n \in \bN.
\end{equation}
For each \(g,\rho \in  L^2\), \(\ve >0\), let \(\rho_{\ve} \in \!\!BV\) such that \(\|\rho - \rho_{\ve}\|_{L^2}\le  \ve\). Using \eqref{eq:power-bounded}, and the first part of the Lemma,
\[
\begin{split}
\lim_{n\rightarrow \infty}\left|\int_{\bT} g\cL^n \rho - \int_{\bT}\rho\int_{\bT} g\rho_{\star}\right| &\le
\lim_{n\rightarrow \infty}\left|\int_{\bT} g\cL^n \rho_\ve - \int_{\bT}\rho_\ve\int_{\bT} g\rho_{\star}\right| + \Const \|g\|_{L^2} \ve\\
&=\Const \|g\|_{L^2} \ve,
\end{split}
\]
which by the arbitrariness of \(\ve\) concludes the proof.
\end{proof}
Similar results, that will be used in the main text, can be obtained for multiple times correlations.

\begin{remark}\label{rem:uniqueabs} If we apply Lemma \ref{lem:dec_cor} to $\rho=1$, it implies that $\rho_\star\geq 0$, hence $\rho_\star d\theta$ is indeed a probability measure. Moreover, if $\rho\in L^1$ is the density of an invariant probability measure, then it can be approximated in $L^1$ by functions $\rho_\ve\in \cC^1(\bT,\bR)$. Then, Lemma \ref{lem:dec_cor} implies
\[
\begin{split}
\int_\bT \rho \vf=\lim_{n\to\infty}\int_\bT \rho(\vf\circ f^n)&=
\lim_{n\to\infty}\int_\bT \rho_\ve(\vf\circ f^n)+\cO(\ve\|\vf\|_\infty)
\\
&=\int_\bT \rho_\star\vf+\cO(\ve\|\vf\|_\infty).
\end{split}
\]
That is, by the arbitrariness of $\ve$, $\rho=\rho_\star$. In other words, $f$ has a unique absolutely continuous invariant probability measure.
\end{remark}

Theorem \ref{thm:decay} is well-known (e.g. see \cite[Chapter 1]{DKL21}). However, we provide a brief proof for the reader's convenience.

\begin{proof}[\bf \em Proof of Theorem \ref{thm:decay}]For all $\rho\in \cC^1$ we have, by a direct computation,
\begin{equation}\label{eq:pre_LY}
(\cL \rho)'=\cL\left( \frac{\rho'}{f'}\right)-\cL\left(\frac{f''}{(f')^2}\rho\right).
\end{equation}
This implies
  \[
    \|(\cL \rho)'\|_{L^1}\leq\left \|\cL\left( \frac{\rho'}{f'}\right)\right\|_{L^1}
    + \left\|\cL\left(\frac{f''}{(f')^2}\rho\right)\right\|_{L^1} \le \lambda^{-1}\|\rho'\|_{L^1}
    + B\|\rho\|_{L^1}
    \]
where 
\begin{equation}\label{eq:distorsion}
B=\left\| \frac{f''}{(f')^2}\right\|_\infty. 
\end{equation}
More generally, Equation \eqref{eq:pre_LY} implies the following basic property of $\cL$: for each $\rho\in\cC^1$, we have\footnote{ Recall that $\|\rho\|_{W^{1,1}}=\|\rho\|_{L^1}+\|\rho'\|_{L^1}$.}
\begin{equation}\label{eq:LY}
\begin{split}
&\|\cL \rho\|_{L^1}\leq \|\rho\|_{L^1}\\
&\|\cL \rho\|_{W^{1,1}} = \|(\cL \rho)'\|_{L^1} + \|(\cL \rho)\|_{L^1}
\leq \lambda^{-1}\|\rho\|_{W^{1,1}} + \left(B +1\right)\|\rho\|_{L^1},
\end{split}
\end{equation}
these are the so-called Lasota-Yorke, or Doeblin-Fortet, inequalities.
Hence, $\cL\in L(L^1,L^1)\cap L(W^{1,1},W^{1,1})$ and, iterating, Equations \eqref{eq:LY} imply that 
\begin{equation}\label{eq:power_bounded}
  \begin{split}
\|\cL^n \rho\|_{W^{1,1}}&\leq \lambda^{-n} \|\rho\|_{W^{1,1}}
  + (B+1) (1-\lambda^{-1})^{-1}\|\rho\|_{L^1}\\
&\leq \left(\lambda^{-n}+ (B+1) (1-\lambda^{-1})^{-1} \right)\|\rho\|_{W^{1,1}}.
\end{split}
\end{equation}
 Thus $\sigma_{W^{1,1}}(\cL)\subset \{z\in\bC\;:\; |z|\leq 1\}$, and by \eqref{eq:power_bounded} the sequence
\begin{equation}
\label{eq:krylov-sequence}
\frac{1}n\sum_{k=0}^{n-1}\cL^k 1
\end{equation}
is compact in $L^1$ and hence has accumulation points,
which belong to $W^{1,1}$.
Since $\cL^k 1\geq 0$,
if $\rho_\star$ is one such accumulation point, then $\rho_\star\geq 0$
and $\cL \rho_\star=\rho_\star$, that is
$1\in \sigma_{W^{1,1}}(\cL)$.  To prove that \(\rho_{\star}\) is strictly positive, note that if there exists \(\overline \theta\) such that  \(\rho_{\star}(\overline \theta) = 0\), then, for any \(n \in \bN\),
\begin{equation}
\label{eq:preimages-trick}
     0 = \rho_{\star}(\overline \theta) = \cL^{n}\rho_{\star}(\overline \theta) = \sum_{y \in f^{-n}(\overline \theta)} \frac{\rho_{\star}(y)}{(f^n)'(y)}.
\end{equation}
The right-hand side is a sum of non-negative terms, hence \(\rho_{\star}(y) = 0\) for any \(y \in f^{-n}(\overline \theta)\), \(n \in \bN\), which, by Lemma \ref{lem:dens-preimages}, forms a dense set. Since \(\rho_{\star} \in W^{1,1}\subset \cC^0\), this means that \(\rho_{\star} \equiv 0\), which is a contradiction. 

The above also implies that $1$ must be a simple eigenvalue. Indeed, suppose that $\cL\rho=\rho$, $\rho\in W^{1,1}\setminus\{\alpha \rho_\star\;:\;\alpha\in\bR\}$.\footnote{ Note that if $\cL(\alpha+i\beta)=\alpha+i\beta$, then it must be $\cL\alpha=\alpha$ and $\cL\beta=\beta$, so we can restrict to real functions.} Then there exists $b>0$ such that $\rho_{\star}- b\rho\geq 0$ but there exists $\overline \theta\in\bT$ such that $\rho_{\star}(\overline \theta)- b\rho(\overline \theta)=0$. Then the same argument as before yields the contradiction $\rho_{\star}= b\rho$.

To study the rest of the spectrum, note that if $\bV=\left\{\rho\in W^{1,1}\;:\; \int_{\bT}\rho =0\right\}$, then, by \eqref{eq:L_deriv}, $\cL \bV\subset \bV$. Moreover, for all $\rho\in\bV$, $\|\rho\|_{L^1}\leq \|\rho'\|_{L^1}$. This implies that, for each $\rho\in\bV$,
\begin{equation}\label{eq:gap_easy}
\|(\cL \rho)'\|_{L^1}\leq \left(\lambda^{-1}+B\right)\|\rho'\|_{L^1}.
\end{equation}
The above implies that the spectrum $\sigma_{\bV}(\cL)\subset\{z\in\bC\;:\; |z|\leq \lambda^{-1}+B\}$. If $\lambda^{-1}+B=:\nu<1$, which is the case if the map is close to linear, then Equation \eqref{eq:gap_easy} implies
\begin{equation}
\label{eq:gap-W}
\sigma_{W^{1,1}}(\cL)\subset \{1\}\cup \{z\in\bC\;:\; |z|\leq \nu\}.
\end{equation}

This proves the theorem apart from the fact that $\cL$ is acting on $W^{1,1}$ rather than $BV$, an issue that we will discuss shortly.

It is necessary to make a more sophisticated argument to prove the theorem in the general case, in which $B$ may be large.

From \eqref{eq:LY} and Hennion's theorem (\cite{He93}, but see \cite[Appendix B]{DKL21} for a self-contained exposition) it follows that $\cL$, when acting on $W^{1,1}$, has essential spectral radius bounded by $\lambda^{-1}$.\\
If $\cL\rho=e^{i\theta_\star}\rho$, $\int_{\bT}|\rho|=1$, $\rho\neq\rho_\star$.
Then,
\[
|\rho|\leq \cL|\rho|
\]
implying $ |\rho|= \cL|\rho|$,  since $\int_{\bT}\left(\cL|\rho| -|\rho|\right)=0$ and the integrand is non-negative. 
Hence, $\rho=\rho_\star e^{i\psi}$ for some function $\psi\in\cC^0([0,1],\bR)$, $\psi(0)-\psi(1)=2k\pi$, and $k\in\bN$. Moreover, \(\cL(\rho_{\star} e^{i \psi}) = e^{i\theta_{\star}}e^{i\psi}\rho_{\star}\) i.e., \(\cL\left(e^{i\psi-i\psi\circ f - i \theta_{\star}}\rho_\star\right)-\rho_\star=0\). Integrating and taking the real part one has
\[
\int_{\bT}\rho_\star\left(1-\cos(\psi-\psi\circ f - \theta_{\star})\right)=0,
\]
which is possible only if $\psi-\psi\circ f - \theta_{\star}=2k\pi$. Since $f$ has at least one fixed point, we have \(\theta_{\star} = 2k\pi\) and \(1\) is the only eigenvalue on the unit circle. 

Thus, there exists $\nu\in (\lambda^{-1},1)$ such that $\sigma_{W^{1,1}}\cL\setminus\{1\}\subset \{z\in\bC\;:\; |z|\leq \nu\}$ (see \cite[Proposition 1.3]{DKL21} for a more detailed discussion). Accordingly, we have the spectral decomposition \eqref{eq:spectral} (for more details, see, e.g., \cite[Section 1.2]{DKL21}). The fact that $\rho_\star\in \cC^1$ can be found in \cite[Problem 1.7]{DKL21}. 
 
To extend $\cL$ to $BV$ note that for each $\rho\in BV$ there exists $\{\rho_n\}\subset W^{1,1}$ such that $\| \rho_n\|_{W^{1,1}}\leq \|\rho\|_{BV}$, $\rho_n$ converges to $\rho$ in $L^1$, and $\liminf \|\rho_n\|_{W^{1,1}} \ge  \|\rho\|_{BV}$ (e.g. see \cite[Section 5.2.2 Theorem 2]{EG}). This and \eqref{eq:LY} immediately imply that, for each $\rho\in BV$,
\begin{equation}\label{eq:LYBV}
\begin{split}
&\|\cL \rho\|_{L^1}\leq \|\rho\|_{L^1}\\
&\|\cL \rho\|_{BV}\leq \lambda^{-1}\|\rho\|_{BV}+( B +1)\|\rho\|_{L^1},
\end{split}
\end{equation}
and by Hennion's theorem the spectral decomposition \eqref{eq:spectral} holds also for $\cL\in L(BV,BV)$. 
\end{proof}
Note that, iterating \eqref{eq:LYBV}, we have, for all \(n \in \bN\),
\begin{equation}\label{eq:LYBV2}
\begin{split}
&\|\cL^n \rho\|_{L^1}\leq \|\rho\|_{L^1}\\
&\|\cL^n \rho\|_{BV}\leq \lambda^{-n}\|\rho\|_{BV}+\frac{B+1}{1-\lambda^{-1}}\|\rho\|_{L^1}.
\end{split}
\end{equation}
Another very useful fact is that densities with bounded logarithmic derivatives constitute an invariant set: let 
\begin{equation}
\label{eq:cones}
C_a=\left\{\rho\in\cC^1\;:\; \rho(\theta)>0, \quad \frac{|\rho'(\theta)|}{\rho(\theta)}\leq a,\; \forall \theta\in\bT\right\}.
\end{equation}
Then, since $|\rho'| \le a  \rho$ if $\rho\in C_a$, by \eqref{eq:pre_LY}, we have
\begin{equation}\label{eq:cone}
\left|  \frac{d}{d\theta}\cL \rho(\theta)\right|\leq \lambda^{-1}\cL|\rho'|(\theta)+B\cL \rho(\theta)\leq
  (\lambda^{-1}a+B)\cL \rho(\theta)    
\end{equation}
so $\cL C_a\subset C_a$, provided $a>\frac B{1-\lambda^{-1}}$.\\

The next result is used multiple times in the text.
\begin{lemma}
\label{lem:decay-of-correlations}
For all $K\in\bN$, \(\{q_{j}\}_{j=0}^{K-1} \in BV\), $\max_j\|q_j\|_{BV}\geq 2$,
\[
\begin{split}
   &  \sum_{j=0}^{K-1}\left\|\cL^j\left(q_j - \rho_{\star}\int_{\bT}q_j(\theta)d\theta\right)\right\|_{L^1}   \\
&\hspace{1.3cm}\le \Const\left(\max_j\|q_j\|_{L^1}\right) \ln \left(\max_j\|q_j\|_{BV}\right)+\Const.
\end{split}
\]
\end{lemma}
\begin{proof}
Denote \(\widetilde{q}_j = q_j - \rho_{\star}\int_{\bT}q_{j}(\theta)d\theta\).
Let $M= \lf{C_1 \ln\left(\max_j\|q_j\|_{BV}\right)} \ge \lf{C_1 \ln 2}$, where $C_1\in\bR_+$ will be chosen later. Then,
\[
\begin{split}
 &\sum_{j=0}^{K-1}\left\|\cL^j \widetilde{q}_j\right\|_{L^1}  = \sum_{j=0}^{\min\{M,K-1\}}\left\|\cL^j\widetilde{q}_j\right\|_{L^1} +\!\!\!\!\!\! \sum_{j=\min\{M+1,K\}}^{K-1}\left\|\cL^j\widetilde{q}_j\right\|_{BV}.
\end{split}
\]
Using the first of \eqref{eq:LY} we have
\[
  \sum_{j=0}^{\min\{M,K-1\}}\left\|\cL^j\widetilde{q}_j\right\|_{L^1} \le 2M\max_j\|q_j\|_{L^1}.
\]
On the other hand, by Equation \eqref{eq:spectral-gap-Q},
\[
\begin{split}
 \sum_{j=\min\{M+1,K\}}^{K-1}&\left\|\cL^j \widetilde{q}_j\right\|_{BV} = \sum_{j=\min\{M+1,K\}}^{K-1}\left\|Q^j \widetilde{q}_j\right\|_{BV} \leq \Const \sum_{j=M+1}^{K-1}\nu^j \|q_j\|_{BV}\\
&\le \Const \nu^{M+1}\max_j\|q_j\|_{BV} \le \Const \left(\max_j\|q_j\|_{BV}\right)^{1-C_1\ln\nu^{-1}},
\end{split}
\]
and choosing $C_1 \ge \left(\ln \nu^{-1}\right)^{-1}$  yields the result.
\end{proof}
Note that, for \(g \in L^{\infty}\) and \(q_j\) as in Lemma \ref{lem:decay-of-correlations}, we have, for all \(N\),
\begin{equation}
\label{eq:corollary-lemma.decay}
\begin{split}
\sum_{j=0}^{N-1}\biggl| \int_{\bT}&g(f^j \theta)q_j(\theta)d\theta - \int_{\bT}g(\theta)\rho_{\star}(\theta)d\theta\sum_{j=0}^{N-1}\int_{\bT}q_j(\theta')d\theta'\biggr|\\
&\le  \Const\|g\|_\infty\left(\max_j\|q_j\|_{L^1}\right) \ln \left(\max_j\|q_j\|_{BV}\right)+\Const\|g\|_{\infty}.
\end{split}
\end{equation}

%%%%%%%%%%%%%%%%%%%%%%%%%%%%%%
\section{Coboundaries}
\label{subsec:coboundaries}
In this appendix, we provide a criterion to check that $\gamma\neq 0$. See Equation \eqref{eq:green-Kubo} for the definition of $\gamma$.\\
We say that a function \(\mathfrak c \in \cC^0(\bT, \bR)\) is a continuous coboundary if there exists \(\zeta \in \cC^0(\bT, \bR)\) such that \(\mathfrak c =  \zeta-\zeta \circ f \). We call \(L^2(\rho_{\star})\) the \(L^2\) functions on \(\bT\) with respect to the measure \(\rho_{\star}(\theta)d\theta\). 
\begin{lemma}
\label{lem:L2menogamma}
For all \(n \in \bN\),
\[
\left|\; \left\|\frac 1{\sqrt{n}}\sum_{k=0}^{n-1}b\circ f^k\right\|^2_{L^2(\rho_{\star})} -2\gamma\right| \le \frac{\Const}n.
\]
\end{lemma}
\begin{proof}
By a direct computation, and recalling \eqref{eq:green-Kubo},
\[
\begin{split}
& \left\|\sum_{k=0}^{n-1}b\circ f^k\right\|^2_{L^2(\rho_{\star})}= \sum_{k,j=0}^{n-1}\int_{\bT} (b\circ f^k)(\theta)(b \circ f^j)(\theta) \rho_{\star}(\theta)d\theta\\
&=\sum_{k=0}^{n-1}\int_{\bT}b^2(\theta) \rho_{\star}(\theta)d\theta + 2\sum_{j=0}^{n-2}\sum_{k>j}^{n-1}\int_{\bT}(b\circ f^{k-j})(\theta) b(\theta) \rho_{\star}(\theta)d\theta\\
&= n \int_{\bT}b^2(\theta) \rho_{\star}(\theta)d\theta  + 2\sum_{s = 1}^{n-1}(n-s)\int_{\bT} (b \circ f^{s})(\theta) b(\theta) \rho_{\star}(\theta)d\theta\\
&= 2n\gamma -2n\sum_{s = n}^{\infty}\int_{\bT} (b \circ f^{s})(\theta) b(\theta) \rho_{\star}(\theta)d\theta-2 \sum_{s=1}^{n-1}s \int_{\bT} (b \circ f^{s})(\theta) b(\theta) \rho_{\star}(\theta)d\theta.
\end{split}
\]
Finally, by Lemma \ref{lem:dec_cor}, \(\left| \int_{\bT} (b \circ f^{s}) b \rho_{\star}\right| \le \Const \nu^{s}\), for some \(\nu < 1\).
\end{proof}
\begin{lemma}
\label{lem:coboundaries}
\(\gamma = 0\) if and only if \(b\) is a continuous coboundary.
\end{lemma}
\begin{proof}
If \(b\) is a continuous coboundary, then Lemma \ref{lem:L2menogamma} implies
\[
 2\gamma=\lim_{n\to\infty} \left\|\frac{1}{\sqrt n}\sum_{k=0}^{n-1}b\circ f^k\right\|^2_{L^2(\rho_{\star})}
=\lim_{n\to\infty} \left\|\frac{1}{\sqrt n}\left(\zeta\circ f^n-\zeta\right)\right\|^2_{L^2(\rho_{\star})}.
\]
Conversely, if \(\gamma = 0\), then Lemma \ref{lem:L2menogamma} implies \(\left\|\sum_{k=0}^{n-1}b\circ f^k\right\|_{L^2(\rho_{\star})} \le \Const\),  for all \(n \in \mathbb N\). Thus, the sequence \(\left\{\sum_{k=0}^{n-1}b\circ f^k\right\}_n\) admits weakly convergent subsequences in \(L^2(\rho_{\star})\). Let \(\mathfrak \zeta \in L^2\) be an accumulation point. Then, recalling Lemma \ref{lem:dec_cor} and that $b$ is zero average, for each \(\psi \in W^{1,1} \subset \cC^0\), and some \((n_j)_j\), we have
\[
\begin{split}
&\int_{\bT} \zeta \circ f (\theta) \psi(\theta)\rho_{\star}(\theta)d\theta=\int_{\bT} \zeta (\theta) \cL (\psi\rho_{\star})(\theta)d\theta\\
&=\lim_{j \rightarrow \infty}\sum_{k=0}^{n_j - 1}\int_{\bT} b \circ f^{k+1} (\theta)\psi(\theta)\rho_{\star}(\theta)d\theta\\
&= - \int_{\bT}b(\theta)\psi(\theta)\rho_{\star}(\theta)d\theta + \lim_{j \rightarrow \infty}\sum_{k=0}^{n_j }\int_{\bT}b \circ f^k(\theta)\psi(\theta)\rho_{\star}(\theta)d\theta \\
& = - \int_{\bT}b(\theta)\psi(\theta)\rho_{\star}(\theta)d\theta + \int_{\bT} \zeta (\theta) \psi(\theta)\rho_{\star}(\theta)d\theta.
\end{split}
\]
Since \(W^{1,1}\) is dense in \(L^2\) and $\rho_\star>0$ (see Theorem \ref{thm:decay}) it follows that
\begin{equation}\label{eq:cob-L2}
        b =  \zeta- \zeta \circ f ,
\end{equation}
as elements of \(L^2\), and we assume without loss of generality that \(\int \zeta \rho_{\star} = 0\). It remains to show that \(\zeta \in \cC^0\). Multiplying by $\rho_\star$ and applying \(\cL\) to \eqref{eq:cob-L2} yields
\[
 \cL (  b \rho_{\star} ) = - (\mathbbm{1}-\cL)(  \zeta \rho_{\star} ).
\]
Hence, for all $\psi\in L^2$, using the decay of correlation in Lemma \ref{lem:dec_cor}, it follows
\[
\sum_{k=0}^{\infty}\int \psi\cL^{k+1} (  b \rho_{\star} ) =\lim_{n\to\infty}\int \psi\left[\cL^n(  \zeta \rho_{\star} )-(  \zeta \rho_{\star} )\right]=-\int \psi\zeta \rho_{\star}.
\]
Since \( \int_{\bT}b(\theta)\rho_{\star}(\theta) d\theta = 0\),  \(\sum_{k=0}^{\infty}\cL^{k+1} (  b \rho_{\star} )=(\Id - \cL)^{-1}\cL(  b \rho_{\star} )\). Hence,
\[
(\Id - \cL)^{-1}\cL(  b \rho_{\star} )= - \zeta \rho_{\star}
\]
as elements of $L^2$. Finally, since \(\rho_{\star}>0\), we have
\[
\zeta  = -\frac{1}{\rho_{\star}}(\mathbbm{1} - \cL)^{-1}\cL(b\rho_{\star})\in W^{1,1} \subset \cC^0.
\]
\end{proof}

The importance of Lemma \ref{lem:coboundaries} is that being a continuous coboundary is relatively easy to check. Indeed, if \(\mathfrak c\) is a continuous coboundary, for any \(f\)-invariant measure \(\nu\), we have \(\int_{\bT} \mathfrak c d\nu = 0\). In particular, if there exists just one periodic orbit \(\{p,f(p),...,f^{M}(p) = p\}\), \(M \in \bN\), such that \(\sum_{k=0}^{M-1}b(f^k (p)) \neq 0\), then \(b\) cannot be a coboundary and \(\gamma>0\). 

%%%%%%%%%%%%%%%
\section{Standard pairs invariance}
\label{sec:Standard-Pairs}

This section is mainly devoted to the proof of Lemma \ref{lem:density} and Theorem \ref{thm:standard-pair-invariance} and some other technical lemmata that are needed in the main text.
\begin{proof}[\bf{\em Proof of Lemma \ref{lem:density}}]
For each $N\in\bN$, we can partition $\bT^N$ in domains $\bI(a_{j}, b_{j})$ where $|b_{j,x}-a_{j,x}|\in [\delta_\ve,2\delta_\ve]$.\footnote{ See just after Equation \eqref{eq:delta_choice} for the definition of $\bI$.}
Accordingly, given any configuration $(\bar{\bol q}, \bar{\bol p})$ and density $\rho_N$ as in Definition \ref{def:initial0}, we can define the probability measure
\[
\begin{split}
\mu_N (A) =& \int_{\bT^N} A(\bar{\bol q}, \bar{\bol p},\theta) \rho_N(\theta) d\theta\\
=&\sum_j\left[\int_{\bI(a_j,b_j)}\rho_N(\theta)d\theta\right]\int_{\bI(a_j,b_j)} A(\bar{\bol q}, \bar{\bol p},\theta) \frac{\rho_N(\theta)}{\int_{\bI(a_j,b_j)}\rho_N(\theta)d\theta} d\theta.
\end{split}
\]
That is, it can be written as a standard family made of standard pairs with a constant 
$G(\theta) = (\bar{\bol q}, \bar{\bol p})\in \mathfrak S^*(N, \ve, C_0)$ for any $\ve>0$. 
Since any measure $\bar\mu_N$ on $\Sigma_N$ is a weak limit of convex combinations of delta measures, the statement follows. 
\end{proof}

\begin{proof}[{\bf \emph{Proof of Theorem \ref{thm:standard-pair-invariance}}}]
Since $\phi_t$ is a continuous map, it suffices to prove that $\phi^t_\star(\mathfrak F_{N,\ve, C_0})\subset \mathfrak F_{N,\ve, C_0}$.
We first prove the statement for \(t = \ve\).
Without loss of generality, we can assume that \(\mu\in \mathfrak F_{N,\ve, C_0}\) consists of a single standard pair \(\ell \in \mathfrak S^*(N,\ve, C_0)\). \\
Recall that
\[
\begin{split}
z(\ve, G_\ell(\theta), \theta) &= e^{(-\bA + \ve^{-\frac 12} \bB_{\ve}(0,\theta))\ve}G_\ell(\theta),\\
\theta_x(\ve, \theta) &= f(\theta_x)
\end{split}
\]
are the solutions of \eqref{eq:det-matrix} and \eqref{eq:thetadyn} at time \(\ve\) starting at \((G_{\ell}(\theta),\theta)\) at time zero. Thus, for any \(g \in \cC^0( \Sigma_N\times\bT^N, \bR)\):
    \[
        \begin{split}
        &\phi^\ve_\star\mu_{\ell}(g) = \mu_{\ell}(g \circ \phi^\ve)\\
        &= \int_{\bI_\ell} g(z(\ve, G_{\ell}(\theta), \theta), f(\theta))\rho_{\ell}(\theta)d\theta\\
        &=\int_{\bI_\ell}g\biggl( e^{\bigl(-\bA + \ve^{-\frac 12}\bB_{\ve}(0,\theta)\bigr)\ve}G_{\ell}(\theta),f(\theta)\biggr)\prod_{x \in \bZ_N}\rho_{\ell,x}(\theta_x)d\theta_{x}.
        \end{split}
    \]
We call \(\varphi_x\) the inverse of \( f|_{[a_x,b_x]}\), where $\bI_\ell=\bI(a,b)$. Setting \(\widetilde \theta_x= f(\theta_x)\),
    \begin{equation}
\label{eq:standard-pairs}
        \int_{\bT^N} g\left( e^{-\bA\ve + \ve^{\frac 12}\bB_{\ve}(0,\varphi(\widetilde \theta))}G_{\ell}(\varphi(\widetilde{\theta})),\widetilde{\theta} \right)\prod_{x \in \bZ_N}(\cL(\rho_{\ell,x} \Id_{[a_x,b_x]}))(\widetilde \theta_x)d\widetilde{\theta}_{x},
    \end{equation}
 where  \(\cL\) is the SRB transfer operator associated to the map \(f\): for $\rho\in BV$,
 \[
\cL \rho(\theta)=\sum_{\theta'\in f^{-1}(\theta)}\frac{\rho(\theta')}{f'(\theta ')},
\]
see Appendix \ref{sec:TO} for a detailed explanation. \\
Note that $\cL(\rho_{\ell,x} \Id_{[a_x,b_x]})=\cL(\rho_{\ell,x}) \Id_{f([a_x,b_x])}$.
 By Equation \eqref{eq:cone}, choosing\footnote{Recall that $B=\left\| \frac{f''}{(f')^2}\right\|_\infty$, see definition \eqref{eq:distorsion}.} \(C_0>\frac{B}{1-\lambda^{-1}}\), $\cL\rho_{\ell,x}$ belongs to $\cD_{C_0}(\bar f(a),\bar f(b))$, $\bar f(a)_x=f(a_x)$ and $\bar f(b)_x=f(b_x)$. Since \(f\) is expanding,  \(\left|f\left([a_x,b_x]\right)\right| > \delta_\ve/2\). However, it can happen that \(\left|f\left([a_x,b_x]\right)\right| > \delta_\ve\), so the domain becomes too large to be a standard pair. In this case we can break the image in finitely many pieces indexed by \(q \in \fQ\), such that \( |b^{(q)}_x-a^{(q)}_x| \in [\frac{\delta}{2},\delta]\) and define, for all $\tilde \theta\in \bI(a^{(q)}, b^{(q)})$,
\[
\begin{split}
&Z_{q} = \int_{\bI(a^{(q)},b^{(q)})}\prod_{x \in \bZ_N} \cL\rho_{\ell,x}(\theta_x)d\theta_x \\
&\rho^{(q)}(\widetilde\theta) = Z_{q}^{-1}\prod_{x \in \bZ_N}\left(\cL\rho_{\ell,x}(\widetilde \theta_x)\right)\Id_{\left[a^{(q)}_x,b^{(q)}_x\right]}(\widetilde\theta_x)\\
&\widetilde{G}^{(q)}(\widetilde{\theta}) =  e^{\bigl(-\bA + \ve^{-\frac 12}\bB_{\ve}(0,\varphi (\widetilde \theta))\bigr)\ve}G_{\ell}(\varphi(\widetilde{\theta})).
\end{split}
\]
This allows to rewrite \eqref{eq:standard-pairs} as
\[
\begin{split}
 \sum_{q \in \fQ} Z_q \int_{\bI(a^{(q)},b^{(q)})}& \rho^{(q)}(\widetilde \theta)
 g( \widetilde{G}^{(q)}(\widetilde{\theta}),\widetilde{\theta} )d\widetilde{\theta},
\end{split}
\]
where \(\rho^{(q)}_x \in \cD_{C_0}(a^{(q)},b^{(q)})\), and \(\sum_{q \in \fQ}Z_q = 1\). It thus remains to check that $\widetilde{G}^{(q)}\in\fS(N,\ve,C_0)$ for each $q\in\fQ$.\\
By the conservation of energy, \(\widetilde{G}^{(q)} \in \cC^2(\bI(a^{(q)},b^{(q)}), \Sigma_N)\). 
For all \(\vartheta\in\bT^N\),
\begin{equation}\label{eq:derivative-matrix-standardpairs-inv}
    \begin{split}
\left\|\frac{\partial e^{\bigl(-\bA + \ve^{-\frac 12}\bB_{\ve}(0,\vartheta)\bigr)\ve}}{\partial \vartheta_x}\right\| &\le  \sum_{k=1}^{\infty} \left\|\frac{\partial (-\bA+ \ve^{-\frac 12}\bB_{\ve}(0, \vartheta))^k}{\partial \vartheta_x}\right\|\frac{\ve^k}{k!}\\
&\leq  \sum_{k=1}^{\infty} k \Const^k\frac{\ve^{\frac k2}}{k!} = \Const\sqrt{\ve}e^{\Const\sqrt{\ve}}.
    \end{split}
\end{equation}
Using the above inequality, since \(\left\|\frac{\partial \varphi_x^{(q)}}{\partial \widetilde \theta_x}\right\|_\infty \le \lambda^{-1}\), the chain rule yields
\[
    \begin{split}
    &\biggl\|\frac{\partial \widetilde{G}^{(q)}}{\partial \widetilde{\theta}_{x}}\biggr\| \le  \biggl\|\frac{\partial e^{\bigl(-\bA + \ve^{-\frac 12}\bB_{\ve}(0,\vf(\widetilde \theta))\bigr)\ve}}{\partial \widetilde \theta_x}G_{\ell}^{(q)}(\vf(\widetilde \theta)) \\
&\hspace{5cm}+ e^{\bigl(-\bA + \ve^{-\frac 12}\bB_{\ve}(0,\vf(\widetilde\theta))\bigr)\ve}\frac{\partial G_{\ell}^{(q)}(\vf(\widetilde \theta))}{\partial \widetilde \theta_x} \biggr\|\\
    &\le   \lambda^{-1}\Const \sqrt{\ve}e^{\Const \sqrt{\ve}}\biggl\|G_{\ell}(\vf(\widetilde \theta))\biggr\|+\lambda^{-1}e^{\Const \sqrt{\ve}}  \biggl\|\frac{\partial G_{\ell}(\theta)}{\partial \theta_x}\biggr\|\\
    &\leq\lambda^{-1}\Const e^{\Const \sqrt{ \ve}} \sqrt{N\ve}+ \lambda^{-1}e^{\Const \sqrt{\ve}} C_0\sqrt{N\ve}.
    \end{split}
\]
If \(\ve\) is small enough, we obtain the claim by choosing $C_0$ sufficiently large.\\
The statement for $t=K\ve$, $K\in\bN$, follows by applying $K$ times the above result, since $\phi^{K\ve}=(\phi^\ve)^K$, see Remark \ref{rem:grid}.
\end{proof}
We end this section with the following useful lemmata.
  \begin{lemma}\label{lemmBV-new}
    Let $\rho$ a probability density on $\bT$,   $\operatorname{supp}(\rho)= I$, where $I=[a,b]$, $|b-a| = \delta<1/2$
    such that  $\rho\in \mathcal C^1(\mathring I,\bR)\cap \cC^0(I,\bR)$, and
    $\int _I |\rho'(u)| du \le C$, for some $C>0$. Then
    \begin{equation}
      \label{eq:9}
      \| \rho \|_{BV} \le 2\delta^{-1} + 3C  + 1.
    \end{equation}
  \end{lemma}
  \begin{proof}
 Since $\rho$ can only have jumps at the boundaries of $[a,b]$,  its total variation is bounded by $C$
 plus $\rho(a) + \rho(b)$. On the other hand, by the mean value theorem, there exists $x_*\in I$ such that $\rho(x_*)=\delta^{-1}$, hence for $x\in \{a,b\}$
    \[
    |\rho(x)|=\delta^{-1}+\left|\int_{x_*}^x\rho'(y)dy\right|\leq \delta^{-1}+C.
    \]
and
\begin{equation}\label{eq:12}
\| \rho\|_{BV} \le \int \rho(x)dx+\rho(a)+\rho(b)+\int|\rho'(x)|dx\leq  1+  2\delta^{-1}  + 3C. 
\end{equation}
\end{proof}

\begin{lemma}
\label{lem:BV-standard-pair}
Let \(\ell \in \fS^*(N, \ve)\) and \(\rho_{\ell, x} \in \cD_{C_0}\)
be the associated density. Then,
\[
\|\rho_{\ell, x} \Id_{I^x_{\ell}}\|_{BV} \le  (6C_0 + 3)\delta_\ve^{-1}.
\]
\end{lemma}
\begin{proof}
  We have
  $|\rho'_{\ell, x}| \le 2 C_0 \rho_{\ell, x}$, so that $\int _I |\rho'_{\ell, x}(\theta)| d\theta \le 2C_0$, so the statement follows directly from Lemma \ref{lemmBV-new} remembering \eqref{eq:delta_choice}.
\end{proof}

\begin{lemma}
\label{lem:BV-estimate}
Let \(\ell \in \fS^*(N,\ve)\) and let \(\mathfrak h \in \mathfrak E_N\). Then,
\[
\sup_{x\in \bZ_N}\sup_{\theta \in \bI_{\ell}}\left|\frac{\partial}{\partial \theta_x}\mathfrak h(G(\theta))\right| \le C_0 \ve^{\frac 12}N\trn{\mathfrak h}.
\]
\end{lemma}
\begin{proof}
Using the Cauchy-Schwarz inequality and the definition of standard pairs,
\[
\begin{split}
\left|\frac{\partial}{\partial \theta_x}\mathfrak h(G_{\ell}(\theta))\right| &= \left|\sum_{k=1}^{4N} \frac{\partial}{\partial z_k}\mathfrak h(G_{\ell}(\theta))\frac{\partial G_{\ell}(\theta)_k}{\partial \theta_x}\right|\\
&\le \left(\sum_{k=1}^{4N} \left(\frac{\partial}{\partial z_k}\mathfrak h (G_{\ell}(\theta))\right)^2\right)^{\frac 12}\left(\sum_{k=1}^{4N} \left(\frac{\partial G_{\ell}(\theta)_k}{\partial \theta_x} \right)^2\right)^{\frac 12}\\
&\le \|\partial_\theta G_{\ell}\|_{\infty}\left(\trn{\mathfrak h}\sqrt{N}\right) \le C_0\sqrt{\ve} \trn{\mathfrak h}N.
\end{split}
\]
\end{proof}

%%%%%%%%%%%
\section{Correlation estimates}\label{sec:techncial}\ \\
To avoid disrupting the main text, we prove here two Lemmata used in Section~\ref{sec:averages}.

\begin{lemma}
\label{lem:covariance}
 Let $I$ be an interval of size $\delta$,  $\rho\in\cC^1(I,\bR)$ be a probability density such that $\|\frac{\rho'}{\rho}\|_\infty\leq 2C_0$. Then for \(h \in (\ve,1)\), and  \(\ve\) small enough,
 \[
\begin{split}
&\left|\ve^{-1}\int_0^h \int_0^s\int_{I}b(f^{\lf{u\ve^{-1}}}\theta)b(f^{\lf{s\ve^{-1}}}\theta)\rho(\theta)d\theta dsdu -  h \gamma\right| \le \Const\ve (\ln\delta^{-1})^2.
\end{split}
\]
\end{lemma}
\begin{proof}
We first divide the integral into two terms,
\begin{equation}
\label{eq:first}
\int_0^h \int_0^{\ve \lf{\ve^{-1}s}}\int_{I}b(f^{\lf{u\ve^{-1}}}\theta)b(f^{\lf{s\ve^{-1}}}\theta)\rho(\theta) d\theta dsdu,
\end{equation}
and
\begin{equation}
\label{eq:second}
\int_0^h \int_{\ve \lf{\ve^{-1}s}}^s\int_{I}b(f^{\lf{u\ve^{-1}}}\theta)b(f^{\lf{s\ve^{-1}}}\theta)\rho(\theta) d\theta dsdu.
\end{equation}
Recalling that $\rho_\star$ is the invariant measure of $f$, we can write \eqref{eq:second} as,
\begin{equation}
\label{eq:zero}
\begin{split}
&\int_0^{\ve\lceil h\ve^{-1}\rceil} (s-\ve \lf{\ve^{-1}s})\int_{\bT}b(f^{\lf{s\ve^{-1}}}\theta)^2\rho(\theta) \Id_{I}(\theta) d\theta ds+\cO(\ve^2)\\
&=\frac{\ve^2}{2}\sum_{j=0}^{\lf{h\ve^{-1}}}\int_{\bT}b(f^{j}\theta)^2\left(\rho(\theta) \Id_{I}(\theta)-\rho_\star(\theta)\right)d\theta ds + \cO(\ve^2)\\
&\phantom{=}
+\frac{\ve^2}{2}\sum_{j=0}^{\lf{h\ve^{-1}}}\int_{\bT}b(f^{j}\theta)^2 \rho_\star(\theta) d\theta ds.
\end{split}
\end{equation}
Note that \(\|\rho \Id_{I}\|_{L^1} =1\) and,  by Lemma \ref{lem:BV-standard-pair}, $\|\rho \Id_{I}\|_{BV} \leq \Const\delta^{-1}$. Hence, by \eqref{eq:corollary-lemma.decay} with $q_j=\rho \Id_{I}-\rho_\star$, for all \(j\) (hence $\int_{\bT}q_j=0$), and $g=b^2$, we have that the first term of \eqref{eq:zero} is \(\cO(\ve^2\ln\delta^{-1})\). Hence, by the invariance of $\rho_\star$, the quantity \eqref{eq:zero} is equal to
\begin{equation}
\label{eq:zero1}
    \frac{\ve h}2\int_{\bT}b(\theta)^2 \rho_\star (\theta)d\theta+\cO(\ve^2\ln\delta^{-1}).
\end{equation}
The term \eqref{eq:first}  is equal to
\begin{equation}
\label{eq:third}
\begin{split}
&\sum_{j=1}^{\lf{h\ve^{-1}}-1}\int_{j\ve}^{(j+1)\ve}\hskip-12pt ds\sum_{k=0}^{j-1}\int_{k\ve}^{(k+1)\ve}\hskip-12ptdu\int_{I}b(f^{\lf{u\ve^{-1}}}\theta)b(f^{\lf{s\ve^{-1}}}\theta)\rho(\theta) d\theta \\
&+ \int_{\ve\lf{h\ve^{-1}}}^{h}\hskip-15pt ds \sum_{k=0}^{\lf{h\ve^{-1}}-1}\int_{k\ve}^{(k+1)\ve}du\int_{I}b(f^{\lf{u\ve^{-1}}}\theta)b(f^{\lf{s\ve^{-1}}}\theta)\rho(\theta) d\theta\\
&=\ve^2\sum_{j=1}^{\lf{h\ve^{-1}}-1}\sum_{k=0}^{j-1}\int_{I}b(f^k\theta)b(f^j \theta)\rho(\theta) d\theta \\
&+ (h-\ve\lf{h\ve^{-1}})\ve \sum_{k=0}^{\lf{h\ve^{-1}}-1}\int_{I}b(f^{\lf{h\ve^{-1}}}\theta)b(f^k \theta)\rho(\theta) d\theta.
\end{split}
\end{equation}
The last line of \eqref{eq:third} can be written as
\begin{equation}\label{eq:car1}
\begin{split}
       &(h-\ve\lf{h\ve^{-1}})\ve \sum_{k=0}^{\lf{h\ve^{-1}}-1}\int_{\bT}b(f^{\lf{h\ve^{-1}}-k}\theta)b(\theta)\left[\cL^k\left(\rho \Id_{I}\right)\right]\!\!(\theta)\, d\theta\\
&=(h-\ve\lf{h\ve^{-1}})\ve \sum_{j=1}^{\lf{h\ve^{-1}}}\int_{\bT}b(f^{j}\theta)b(\theta)\left[\cL^{\lf{h\ve^{-1}}-j}\left(\rho \Id_{I}\right)\right]\!\!(\theta)\, d\theta.
\end{split}
\end{equation}
Recalling \eqref{eq:LYBV2}, and Lemma \ref{lem:BV-standard-pair},
\[
\begin{split}
&\|b\cL^{\lf{h\ve^{-1}}-j}\left(\rho \Id_{I}\right)\|_{L^1} \leq\|b\|_\infty ,\\
&\|b\cL^{\lf{h\ve^{-1}}-j}\left(\rho \Id_{I}\right)\|_{BV} \leq \|b\|_{BV}\|\cL^{\lf{h\ve^{-1}}-j}\left(\rho \Id_{I}\right)\|_{BV}\leq \Const \delta^{-1}_{\ve}.
\end{split}
\]
Using  Equation \eqref{eq:corollary-lemma.decay}, with $q_j=b\cL^{\lf{h\ve^{-1}}-j}\left(\rho \Id_{I}\right)$ and $g=b$, and since \((h-\ve\lf{h\ve^{-1}}) \le \ve\), allows writing the expression \eqref{eq:car1} as
\begin{equation}\label{eq:D7}
\begin{split}
&(h-\ve\lf{h\ve^{-1}})\ve \!\!\!\sum_{j=1}^{\lf{h\ve^{-1}}}\!\!\!\int_{\bT}b(f^{j}\theta)\rho_\star (\theta)d\theta\int_{\bT}b(\theta')\left[\cL^{\lf{h\ve^{-1}}-j}\left(\rho \Id_{I}\right)\right]\!\!(\theta'_x)\, d\theta'\\
&+\cO(\ve^2 \ln \delta^{-1})=\cO(\ve^2 \ln \delta^{-1})
\end{split}
\end{equation}
where we have used the invariance of $\rho_\star$ and the fact that $b$ is zero average.

The line next to the last in \eqref{eq:third} can be written as
\begin{equation}
\label{eq:four}
\begin{split}
 &\ve^2\sum_{j=1}^{\lf{h\ve^{-1}}-1}\sum_{k=0}^{j-1}\int_{\bT}b(f^k\theta)b(f^j \theta)\rho_{\star}(\theta) d\theta \\
&-\ve^2\sum_{j=1}^{\lf{h\ve^{-1}}-1}\sum_{k=0}^{j-1}\int_{\bT}b(f^{j-k}\theta)b(\theta)\left[\cL^k\left(\rho_{\star} - \rho_{x}\Id_{I}\right)\right]\!\!(\theta)\, d\theta=\\
&=\ve^2\sum_{j=1}^{\lf{h\ve^{-1}}-1}\left(\lf{h\ve^{-1}}-1-j\right)\int_{\bT}b(f^{j}\theta)b(\theta)\rho_{\star}(\theta) d\theta \\
 &   -\ve^2\sum_{j=1}^{\lf{h\ve^{-1}}-1}\int_{\bT}b(f^j\theta)b(\theta)\hskip-.2cm\sum_{k=0}^{\lf{h\ve^{-1}}-1 -j}\hskip-.2cm\left[\cL^{k}\left(\rho_{\star}- \rho_{x}\Id_{I}\right)\right]\!\!(\theta) \,d\theta.
\end{split}
\end{equation}
For any \(p \in \bN\), using Lemma \ref{lem:decay-of-correlations} with $q_j=\rho_{\star} - \rho_{x}\Id_{I}$,
\[
\biggl\|b\sum_{k=0}^{\lf{h\ve^{-1}}-1 -j}\cL^{k}\left(\rho_{\star} - \rho_{x}\Id_{I}\right)\biggr\|_{L^1} \le \Const \ln \delta^{-1},
\]
where we have used Lemma \ref{lem:BV-standard-pair} and \(\|\rho_{\star} - \rho_{x}\Id_{I}\|_{L^1} \le 2\).
Furthermore, using the spectral decomposition \eqref{eq:spectral},
\[
\begin{split}
   \biggl\|b\sum_{k=0}^{\lf{h\ve^{-1}}-1 -j}\cL^{k}\left(\rho_{\star} - \rho_{x}\Id_{I}\right)\biggr\|_{BV} &\le \|b\|_{BV} \Const \delta^{-1} \sum_{k=0}^{\lf{h\ve^{-1}}-1 -j}\nu^k\\ &\le \Const \delta^{-1}. 
\end{split}
\]
Therefore, using again \eqref{eq:corollary-lemma.decay} with $q_j=b\sum_{k=0}^{\lf{h\ve^{-1}}-1 -j}\cL^{k}\left(\rho_{\star} - \rho_{x}\Id_{I}\right)$ and $g=b$, we can rewrite \eqref{eq:four} as
\begin{equation}
\label{eq:omega}
\begin{split}
& h \ve \sum_{j=1}^{\infty}\int_{\bT}b(f^j\theta)b(\theta)\rho_\star(\theta)d\theta + \cO(\ve^2(\ln \delta^{-1})^2).
\end{split}
\end{equation}
The Lemma follows by summing \eqref{eq:zero1}, \eqref{eq:D7} and \eqref{eq:omega} and recalling the definition of $\gamma$ in \eqref{eq:green-Kubo}. Indeed, the sum equals 
\[
\ve h\left(\frac 12\int_\bT b^2\rho_\star+\sum_{j=1}^\infty\int_\bT (b\circ f^j)\, b\,\rho_\star\right)+\cO(\ve^2(\ln\delta^{-1})^2)=\ve h\gamma+\cO(\ve^2(\ln\delta^{-1})^2).
\]
\end{proof}

\begin{lemma}
\label{lem:covariance2}
 Let $I$ be an interval of size $\delta$,  $\rho\in\cC^1(I,\bR)$ be a probability density such that $\|\frac{\rho'}{\rho}\|_\infty\leq 2C_0$.Then, for all \(\ve\) small enough and \(h \in (\ve,1)\),
\[
\begin{split}
&\left|\ve^{-1}\int_{0}^{h} \int_{0}^h  \int_{I}b(f^{\lf{(h + s)\ve^{-1}}}\theta)b(f^{\lf{u\ve^{-1}}}\theta)\rho(\theta)d\theta dsdu\right|  \le \Const\ve (\ln\delta^{-1})^2.
\end{split}
\]
\end{lemma}
\begin{proof}
We assume that $h\ve^{-1} =K\in\bN$ as the general case can be treated similarly to Lemma \ref{lem:covariance}. This is always the case in our applications, since $h\in\ve\bN$ in Theorem \ref{thm:basic_fact} and Lemma \ref{lem:annoying-term}. By adding and subtracting \(\rho_{\star}\), the integral in the statement can be written as
\begin{equation}
\label{eq:correlation-at-distance}
\begin{split}
&\ve \sum_{k=0}^{K-1}\sum_{j=0}^{K-1}\int_{\bT}b(f^{\lf{h\ve^{-1}}+k}\theta)b(f^j \theta)(\rho(\theta)\Id_{I}(\theta)- \rho_{\star}(\theta))d\theta+\\
&+\ve \sum_{k=0}^{K-1}\sum_{j=0}^{K-1}\int_{\bT}b(f^{\lf{h\ve^{-1}}+k}\theta)b(f^j \theta) \rho_{\star}(\theta)d\theta.
\end{split}
\end{equation}
Let us consider the first term of \eqref{eq:correlation-at-distance}. Setting \(p=k-j+K\),
\[
\begin{split}
&\ve \sum_{p=1}^{2K-1}\int_{\bT}b(f^p\theta)b(\theta)\hskip-12 pt\sum_{j=\max\{K-p,0\}}^{\min\{K-1, 2K-1 -p\}}\hskip-12 pt\cL^j\left(\rho(\theta)\Id_{I}(\theta)-\rho_{\star}(\theta)\right)d\theta.
\end{split}
\]
By Lemma \ref{lem:decay-of-correlations}, we have,  for any \(p \in \{0,1,..., 2 K -1\}\),
\[
\left\|b\sum_{j=\max\{K-p,0\}}^{\min\{K-1, 2K-1 -p\}}\cL^j\left(\rho\Id_{I}-\rho_{\star}\right)\right\|_{L^1} \le \Const  \ln \delta^{-1},
\]
and, using the spectral decomposition \eqref{eq:spectral},
\[
\begin{split}
\Biggl\|b\sum_{j=\max\{K-p,0\}}^{\min\{K-1, 2K-1 -p\}}&\cL^j\left(\rho \Id_{I}-\rho_{\star}\right) \Biggr\|_{BV} \le \Const \|b\|_{BV}\delta^{-1}\sum_{j=0}^{2K-1}\nu^{j} \le \Const\delta^{-1}.
\end{split}
\]
Therefore, applying \eqref{eq:corollary-lemma.decay}, the first term of \eqref{eq:correlation-at-distance} is of order \(\cO(\ve(\ln \delta^{-1})^2)\). Setting \(p=k-j+K\) and noting that 
\[
\min\{K-1, 2K-1 -p\}- \max\{K-p,0\}\leq p-1,
\]
the second term of \eqref{eq:correlation-at-distance} can be bounded by
\begin{equation}
\label{eq:bb}
\begin{split}
&\ve\sum_{p=1}^{2K-1}\biggl(p-1\biggr)\left|\int_{\bT}b(f^p\theta)b(\theta)\rho_{\star}(\theta)d\theta \right|.
\end{split}
\end{equation}
Note that by Lemma \ref{lem:dec_cor}, \(\left|\int_{\bT}b(f^p\theta)b(\theta)\rho_{\star}(\theta)d\theta\right| \le \Const \nu^{p}\), from which the Lemma follows.
\end{proof}

%%%%%%
\section{Discrete gradient, Laplacian and Fourier}
\label{sec2.7}\ \\
For any $g:\bZ_N\to\R$, $N\in\bN\cup\{\infty\}$, $x\in\bZ_N$, the lattice gradient, and the Laplacian are  defined  as
\begin{equation}
  \label{012901-23}
(\nabla_{\bZ_N}g)_x=f_{x+1}- f_{x}\;,\quad (\Delta_{\bZ_N}g)_x=g_{x+1}+g_{x-1}-2 g_{x}
\end{equation}
Recall the orthonormal base of eigenvectors $\{\psi_j\}$ and eigenvalues $\upsilon_j$ 
of $- \Delta_{\Z_N}$ defined in \eqref{eq:29}.
\begin{lemma}
  \label{lem:bound}
  Let $\varphi\in \mathcal C^k(\bT)$  and \(\hat \vf_j\) be defined by \eqref{eq:fourier} with \(\vf_x = \vf\left(\frac{x}{N}\right)\). Then, for each $N\in\bN$ and $|j|\leq N/2$,
  \begin{equation}
    \label{eq:37}
   | \hat\varphi_j| \le \frac{\Const  \|\vf\|_{\cC^k} N^{\frac 12}}{|j|^k}.
  \end{equation}
\end{lemma}
\begin{proof}
Note that,
  \begin{equation*}
    \begin{split}
      |N^{-\frac12}\hat\varphi_j| &=\left| \frac{1}{N(e^{-i2\pi j/N}-1)}
       \sum_{x\in\bZ_N} e^{-i2\pi j x/N} \nabla_{\bZ_N}\varphi\left(\frac {x}{N}\right) \right|\\
       &\le \frac{1}{N|e^{-i2\pi j/N}-1|}
       \sum_{x\in\bZ_N} \frac{\|\varphi'\|_\infty}{N}.
    \end{split}
  \end{equation*}
Since for $\theta\in[-\pi/2,\pi/2]$ we have $|1-e^{i\theta}|^2=2|1-\cos\theta|\geq\frac{\theta^2}4$, we have
    \begin{equation*}
   N^{-\frac 12}  |\hat\varphi_j|\leq \frac 1{\pi j}\|\varphi'\|_\infty,
      \end{equation*}
that proves the statement for $k=1$. Then iterate for general $k$.
\end{proof}

%%%%%%%%%%%%
\section{On the diffusion coefficient $D$}
\label{sec:the-diff-expr}
The Green's function of $\omega^2_0-\Delta_{\Z_N}$
is given by
\begin{equation}
    \label{eq:GRN}
 \begin{split}   \left(\omega^2_0-\Delta_{\Z_N}\right)^{-1}(x,y) &= 
    \sum_j v_j^{-1} \psi_j(x)\psi_j^*(y) \\
    &= \frac 1N \sum_j v_j^{-1} e^{i 2\pi (x-y)j/N} := \bG_{\omega_0,N}(x-y).
    \end{split}
\end{equation}
Since $v_j = v_{-j}$, this can be written as
\begin{equation}
    \label{eq:GRN1}
    \bG_{\omega_0,N}(x) = \frac 1N \sum_j v_j^{-1} \cos( 2\pi x j/N).
\end{equation}
On the other hand, the Green function of $-\Delta_{\Z} + \omega^2_0$,
where $\Delta_{\Z}$ is the Laplacian on the integer lattice $\Z$,
that is given by, for $x\in\bZ$,
\begin{equation}
\label{GR}
\begin{split}
&\bG_{\om_0}(x) := \left(-\Delta_{\Z} + \omega^2_0 \right)^{-1}(y, y+x) 
=\int_0^1\frac{\cos(2\pi ux)}{4\sin^2(\pi u)+\om_0^2}du.
\end{split}
\end{equation}
As $\bG_{\omega_0,N}(x)$ is a simple Riemann sum approximation of $\bG_{\omega_0}(x)$, 
\begin{equation}\label{eq:riemann-approx}
\left|\bG_{\om_0,N}(x)-\bG_{\om_0}(x)\right|\leq \Const N^{-1}.
\end{equation}
Using \eqref{GR} the diffusion coefficient appearing in \eqref{eq:D1} reads:
\begin{equation}
  \label{eq:28}
  \begin{split}
    D &= 1 + 2 \bG_{\om_0}(1) - (1+\omega_0^2) \bG_{\om_0}(0) - \bG_{\om_0}(2) \\
   & = 1 + \left(\omega_0^2 - \Delta\right) \bG_{\om_0}(1) - \omega_0^2 \left(\bG_{\om_0}(0) + \bG_{\om_0}(1)\right) \\
     & =  1 - \omega_0^2\left(\bG_{\om_0}(0) + \bG_{\om_0}(1)\right)\\
     &=\int_0^1  \frac{ 1  -  \cos(4\pi u)}{4\sin^2(\pi u)+\om_0^2} du.
  \end{split}
\end{equation}
From the last line of  \eqref{eq:28} an explicit computation yields
\begin{equation}
  \label{eq:32}
  D = \frac{2}{2+\omega_0^2 + \omega_0\sqrt{4+\omega_0^2}}.
\end{equation}
The thermal diffusivity $D$ can also be written differently.
There is a kinetic expression that relates it to the dispersion relation of the harmonic lattice 
$\omega(k) = \sqrt{\omega_0^2 + 4\sin^2(\pi k)}$. In the kinetic limit, the velocity of a \emph{phonon}
of wave number $k$ is $\frac{\omega'(k)}{2\pi}$, and the diffusion coefficient is related to the kinetic energy of the phonons by (see Section 6 in \cite{Basile16} for more details)
\begin{equation}
  \label{eq:Dkin}
  D = \frac {1}{2\pi^2} \int_0^1 \left[\omega'(k)\right]^2 dk,
\end{equation}
that is valid for a generic harmonic lattice. For nearest-neighbor interaction
\begin{equation}
  \label{eq:26}
  \begin{split}
  D&= \frac {1}{2\pi^2} \int_0^1 \left[\omega'(k)\right]^2 dk
  =  \int_0^1 \frac{8\sin^2(\pi k)\cos^2(\pi k)}{\omega(k)^2} dk\\
&  = \int_0^1 \frac{1-\cos(4\pi k)}{\omega(k)^2} dk.
\end{split}
\end{equation}
%%%%%%%%%%%%
\section{Compactness}
\label{sec:compactness}
Recall that \(\mathcal P(\bT)\) is the space of probability measures on \(\bT\) endowed with the weak topology. Under this topology $\mathcal P(\bT)$ is complete, separable, metrizable, and compact.
We choose the metric
\begin{equation}
  \label{eq:33}
  d(\mu,\mu') = \sum_{k\in\bZ} \frac 1{2^{|k|}}
  \min\left\{1, \left|\int \bs\psi_k d\mu - \int \bs\psi_k d\mu'\right|\right\},
\end{equation}
where $\bs\psi_k(\theta)=e^{2\pi i\theta k}$, \(k \in \bZ\), {is a countable set whose span is dense} in $\mathcal C^0(\mathbb T)$ 
by the Stone-Weierstrass theorem.

We have that $\xi_N\in \mathcal C([0,t_*],  \mathcal P(\bT))$, that we endow with the uniform topology
induced by the metric \eqref{eq:33}. 
Since  \(\mathcal P(\bT)\) is compact, relative compactness of the sequence $(\xi_N)$ follows from the following equicontinuity in time of $\xi_N$, by Ascoli's theorem \cite[ pg. 223, Theorem 17]{Kelley}.

\begin{proposition}
Provided $\ve\leq N^{-\alpha}$, $\alpha>6$, we have, for all $t_*>0$,
  \begin{equation}
    \lim_{\delta\to 0} \limsup_{N\to\infty} \sup_{0\le s,t\le t_*, |t-s|\le \delta} d(\xi_N(t),\xi_N(s)) = 0
    \label{eq:56}
\end{equation}
\end{proposition}

\begin{proof}
  Let $\varphi \in \mathcal C^6(\mathbb T)$.
 Recalling the definition of $\xi_N$ in \eqref{eq:6}, by Theorem \eqref{prop:fundamental}, we have, for all $0\le s < t\leq t_*$, 
\[
\begin{split}
&\| \xi_N(t,\varphi) \|\leq  \|\vf\|_\infty,\\
&\xi_N(t,\varphi) -  \xi_N(s,\varphi) =\frac{D}{2\gamma} \int_s^t  \xi_N(u,\varphi'') du + \|\vf\|_{\cC^6}\cO\left((t-s)  +  N^{-1}\right).
  \end{split}
\]
It follows
  \begin{equation}
    \label{eq:58}
    |\xi_N(t,\varphi) - \xi_N(s,\varphi) | \le \Const (|t-s| +N^{-1})\|\varphi\|_{C^6}.
  \end{equation}
Then using the metric \eqref{eq:33} we have
  \begin{equation}
    \label{eq:34}
    \begin{split}
      d(\xi_N(t),\xi_N(s)) \le \sum_{k \in \bZ} \frac 1{2^{|k|}} \min\left\{ 1, \Const  k ^6 (|t-s|+N^{-1})\right\},
    \end{split}
  \end{equation}
  and, for any \(\delta \le 1\),
  \begin{equation}
    \label{eq:35}
 \begin{split}
& \sup_{|t-s|\le \delta} d(\xi_N(t),\xi_N(s)) \le
 \Const  \left(  \sum_{k=0}^{\infty} \frac 1{2^k} \min\left\{ 1, \delta k^6\right\}+ N^{-1} \right)\\
 &\leq   \Const \left( \sum_{k=0}^{\delta^{-\frac {1}{12}}} \frac 1{2^k} \sqrt\delta+ \sum_{k=\delta^{-\frac {1}{12}}}^\infty \frac 1{2^k} +N^{-1} \right)
 \leq \Const  (\sqrt\delta+N^{-1}),
\end{split}
  \end{equation}
 from which the proposition follows.
\end{proof}
\newpage
\section*{List of Notations and Symbols}
\label{sec:notation}
{\footnotesize
\begin{description}[leftmargin=1.8cm]
\item[$\langle v, w \rangle$] Euclidean scalar product in $\mathbb{R}^n$. Page \pageref{page1}

\item[$\|v\|$] $\sqrt{\langle v, v \rangle}$. Euclidean norm. Page \pageref{eq:enenorm}

  \item[$\|M\|$] $\sup_{\|v\|=1}\|Mv\|$. Operator norm of the matrix $M$. Page \pageref{eq:functions-controlled-by-energy-norms}

    \item[$(\nabla_{\mathbb{Z}_N} g)_x$] $g_{x+1} - g_x$. Discrete gradient on $\mathbb{Z}_N$. Pages \pageref{eq:intermediate-heat-step}, \pageref{012901-23}

    \item[$(\Delta_{\mathbb{Z}_N} g)_x$] $g_{x+1} + g_{x-1} - 2g_x$. Discrete Laplacian on $\mathbb{Z}_N$. Pages \pageref{eq:det}, \pageref{012901-23}

    \item[$A_N(t)$] $A_N(z(t,\bar z,\bar\theta))$. Deterministic evolution of $A \in \mathfrak{E}_N$. Page \pageref{eq:45}

 \item[$\lvvvert A_N\rvvvert$] Energy-controlled seminorm for $A_N$. Page \pageref{eq:norm-a}
    \item[$\mathbb{A}$] Generator of the dynamics on $\mathbb{Z}_N$, harmonic force. Page \pageref{eq:unperturbed-matrices}

    \item[$b(\theta)$] External field; Pages \pageref{eq:det}, \pageref{eq:zero-average}

    \item[$\mathbb{B}_\varepsilon(t,\theta)\varepsilon^{-\frac{1}{2}}$] Generator of the time-varying external field. Page \pageref{eq:unperturbed-matrices}

    \item[$C_\#$] Generic computable constant, depending only on the model parameters $(\omega_0,b,f)$ in \eqref{eq:det} and on $C_0$ in \eqref{eq:Czero}. Page \pageref{eq:enenorm}

    \item[$D$] $\displaystyle \frac{2}{2+\omega_0^2+\omega_0\sqrt{\omega_0^2+4}}$. Effective diffusion coefficient. Pages \pageref{eq:PDE}, \pageref{sec:the-diff-expr}

    \item[$DA$, $D^2A$] Gradient and Hessian of a $\cC^2$ function $A$. Page \pageref{eq:functions-controlled-by-energy-norms}

    \item[$\cD_{C_0}(a,b)$] Admissible densities of standard pairs. Page \pageref{eq:density_def}

    \item[$\delta_0$, $\delta_{\ve}$] Size of the support of standard pairs. Page \pageref{eq:delta_choice}

    \item[$\fe_x(z)$] $\displaystyle \frac{p_x^2}{2} + \frac{(q_x - q_{x-1})^2}{2} + \frac{\omega_0^2 q_x^2}{2}$. Local energy of particle $x$. Page \pageref{eq:energy}

    \item[$\mathcal{E}_0$] Initial macroscopic energy profile (limit of microscopic profiles). Page \pageref{eq:initial-conditions}

    \item[$\mathcal{E}(t)$] Macroscopic energy profile at time $t$. Page \pageref{eq:wheat}

    \item[$\mathcal{E}_N(z)$] $\sum_{x\in\mathbb{Z}_N} e_x(z)$ total energy (conserved and equal to $N$). Page \pageref{eq:energy}

    \item[$\mathfrak{E}_N$] Set of degree-two polynomials $\Sigma_N \to \mathbb{R}$. Page \pageref{def:contro-energy}

    \item[$f:\mathbb{T}\to\mathbb{T}$] Fast dynamics, uniformly expanding map. Page \pageref{eq:thetadyn}

    \item[$\ff$, $\mu_{\ff}$] A standard family and its associated measure. Page \pageref{defin:standardf}

    \item[$\bF_{N,\ve, C_0}$] Set of standard families. Page \pageref{eq:sf_rep}

    \item[$\fF_{N,\ve, C_0}$] Set of probability measures determined by a standard family. Page \pageref{eq:sf_rep}

    \item[$\mathcal{G}$] $\mathcal{A} + \gamma \mathcal{S}$. Generator of the limiting stochastic dynamics. Page \pageref{eq:5}

\item[$\phi^t$] Flow generated by the microscopic dynamics. Page \pageref{not:2}

\item[$\phi_{\star}^t$] Push forward of measures by the flow $\phi^t$. Page \pageref{not:2}

\item[$G(\theta)$] Graph of standard pairs. Page \pageref{eq:delta_choice}

    \item[$\bG_{\omega_0}$] Green function of the operator $\left(\omega_0^2 - \Delta\right)^{-1}$. Page \pageref{eq:gaussian-equilibrium-chain}

    \item[$\gamma$] $\displaystyle \frac 12\int_{\mathbb{T}} b(\theta)^2\rho_\star(\theta)\,d\theta
    +  \sum_{k=1}^\infty \int_{\mathbb{T}} b(\theta)b(f^k\theta)\rho_\star(\theta)\,d\theta$.
    Self-correlations. Page \pageref{eq:green-Kubo}

    \item[$\bI(a,b)$] Supports for the graphs of standard pairs. Page \pageref{eq:delta_choice}

    \item[$j_{x,x+1}$] $p_x(q_{x}-q_{x+1})$. Microscopic energy current. Page \pageref{eq:current}

    \item[$\ell$, $\mu_{\ell}$]  A standard pair and its associated measure. Page \pageref{defin:standardp}

    \item[$\scriptstyle\widetilde\cM_{\mathrm{init}}(N,C_0)$] Simple example of initial measures on $\Sigma_N \times \mathbb{T}^N$. Sometimes abbreviated  $\widetilde\cM_{\mathrm{init}}$. Page \pageref{def:initial0}
    \item[$\scriptstyle\cM_{\mathrm{init}}(N, \ve,C_0)$] General class of initial measures on $\Sigma_N \times \mathbb{T}^N$ for which the results in paper hold. Sometimes abbreviated as $\cM_{\mathrm{init}}(N)$  or $\cM_{\mathrm{init}}$, see Remark \ref{rem:C0_choice}. Pages \pageref{nota:2}, \pageref{def:initial}
    \item[$\cM^\star_{\mathrm{init}}$] Initial data sequences associated with a macroscopic profile $\mathcal{E}_0$. Page \pageref{def:initialstar}
    \item[$\mu_N$] Initial probability measure belonging to $\mathcal{M}_{\mathrm{init}}(N,C_0)$. Page \pageref{def:initial}

    \item[$\mathbb{N}$] $\{1,2,\dots\}$. Page \pageref{sec:determinsitic}

    \item[$\mathcal{O}(a)$] Quantity bounded in absolute value by $C_\# |a|$ (uniformly in relevant parameters). Page \pageref{lem:Taylor-0}

    \item[$\mathcal{P}(\mathbb{T})$] Space of probability measures on $\mathbb{T}$ endowed with the weak topology. Page \pageref{def:initialstar}

    \item[$p_x$, $q_x$] $p_x = (p_{x,1}, p_{x,2})$, \; $q_x = (q_{x,1}, q_{x,2})$. Momentum and position of particle at site $x \in \mathbb{Z}_N$. Page \pageref{eq:det}

    \item[$\boldsymbol{p}$, $\boldsymbol{q}$, $z$] $\boldsymbol{p} = (p_x)_{x\in\mathbb{Z}_N}$, \;
    $\boldsymbol{q} = (q_x)_{x\in\mathbb{Z}_N}$, \;
    $z = (\boldsymbol{q}, \boldsymbol{p})$. Vectors of momenta and positions. Page \pageref{notation:z}

    \item[$\Psi(N;\mu_N)$] Quantitative convergence rate of the initial profile. Page \pageref{eq:varpi}
    \item[$\mathbb{R}_+$]  $\{a \in \mathbb{R} : a \ge 0\}$. Page \pageref{cor:long_time}
    \item[$\rho_{x,N}$, $\rho_N$] Marginal of $\mu_N \in \mathcal{M}_{\mathrm{init}}$ on $\mathbb{T}$ associated to particle $x$ and product of the marginals. Page \pageref{eq:random}

    \item[$\rho_\star$] Unique absolutely continuous invariant density for $(f,\mathbb{T})$. Page \pageref{eq:inve_den}

    \item[$S_N$] Covariance matrix with $[S_N^{(q)}]_{x,y} = \langle q_x, q_y\rangle$, etc. Page \pageref{S1ts}

    \item[$\Sigma_N$] $\{z \in \mathbb{R}^{4N} : \mathcal{E}_N(z) = N\}$. Energy shell of total energy $N$. Page \pageref{eq:shell}

    \item[$\mathfrak S(N, \ve, C_0)$] Set of admissible graphs of standard pairs. Page \pageref{eq:Graphdef}

    \item[$\mathfrak S^*(N, \ve, C_0)$] Set of standard pairs. Page \pageref{defin:standardp}

    \item[$\mathbb{T}^n$] $\mathbb{R}^n / \mathbb{Z}^n$ (and $\mathbb{T} = \mathbb{T}^1$). $n$-dimensional torus. Page \pageref{sec:determinsitic}

    \item[$\omega_0 > 0$] Pinning parameter. Page \pageref{eq:det}

    \item[$z(t,\bar z,\bar\theta)$] $(\boldsymbol{q}(t,\bar z,\bar\theta), \boldsymbol{p}(t,\bar z,\bar\theta))$. Solution of the deterministic ODE. Page \pageref{lem:energy}

    \item[$\mathbb{Z}_N$] $\mathbb{Z}/N\mathbb{Z}$. One-dimensional discrete periodic lattice with $N$ sites. Page \pageref{eq:det}
\end{description}
}


\newpage
\begin{thebibliography}{100}
\bibitem{Bal00} V. Baladi, 
\emph{Positive Transfer Operators and Decay of Correlations}, 
Advanced Series in Nonlinear Dynamics, Vol.~16, World Scientific, River Edge, NJ, 2000.

\bibitem{Bal18} V. Baladi, 
\emph{Dynamical Zeta Functions and Dynamical Determinants for Hyperbolic Maps: A Functional Approach}, 
Ergebnisse der Mathematik und ihrer Grenzgebiete, Vol.~68, Springer, Cham, 2018.

\bibitem{Bardos93} C. Bardos, F. Golse, C.D. Levermore, 
\emph{Fluid Dynamic Limits of Kinetic Equations II: Convergence Proofs for the Boltzmann Equation}, 
Comm. Pure Appl. Math. \textbf{46} (1993), 667–753.

\bibitem{Bardos97} C. Bardos, F. Golse, J.F. Colonna, 
\emph{Diffusion Approximation and Hyperbolic Automorphism of the Torus}, 
Physica D \textbf{104} (1997), 32–70.

\bibitem{Basile16} G. Basile, C. Bernardin, M. Jara, T. Komorowski, S. Olla, 
\emph{Thermal Conductivity in Harmonic Lattices with Random Collisions}, 
in \emph{Thermal Transport in Low Dimensions: From Statistical Physics to Nanoscale Heat Transfer}, 
S. Lepri (ed.), Lecture Notes in Physics, Vol.~921, Springer, 2016. 
%\href{https://doi.org/10.1007/978-3-319-29261-8_5}{doi:10.1007/978-3-319-29261-8\_5}.

\bibitem{bernardin07} C. Bernardin, 
\emph{Hydrodynamics for a System of Harmonic Oscillators Perturbed by a Conservative Noise}, 
Stochastic Process. Appl. \textbf{117} (2007), 487–513. 
%\href{https://doi.org/10.1016/j.spa.2006.08.006}{doi:10.1016/j.spa.2006.08.006}.

\bibitem{BHLLO} C. Bernardin, F. Huveneers, J.L. Lebowitz, C. Liverani, S. Olla, 
\emph{Green–Kubo Formula for Weakly Coupled Systems with Noise}, 
Comm. Math. Phys. \textbf{334} (2015), 1377–1412.

\bibitem{Bernardin05} C. Bernardin, S. Olla, 
\emph{Fourier Law and Fluctuations for a Microscopic Model of Heat Conduction}, 
J. Stat. Phys. \textbf{118} (2005), 271–289. 
%\href{https://doi.org/10.1007/s10955-005-7578-9}{doi:10.1007/s10955-005-7578-9}.

\bibitem{Bernardin11} C. Bernardin, S. Olla, 
\emph{Transport Properties of a Chain of Anharmonic Oscillators with Random Flip of Velocities}, 
J. Stat. Phys. \textbf{145} (2011), 1224–1255.

\bibitem{bhat22} J.M. Bhat, G. Cane, C. Bernardin, A. Dhar, 
\emph{Heat Transport in an Ordered Harmonic Chain in Presence of a Uniform Magnetic Field}, 
J. Stat. Phys. \textbf{186} (2022), 2. 
%\href{https://doi.org/10.1007/s10955-021-02848-5}{doi:10.1007/s10955-021-02848-5}.

\bibitem{BGS-R} T. Bodineau, I. Gallagher, L. Saint-Raymond, 
\emph{The Brownian Motion as the Limit of a Deterministic System of Hard Spheres}, 
Invent. Math. \textbf{203} (2016), 493–553.

\bibitem{Bodineau20} T. Bodineau, I. Gallagher, L. Saint-Raymond, S. Simonella, 
\emph{Fluctuation Theory in the Boltzmann–Grad Limit}, 
J. Stat. Phys. \textbf{180} (2020), 873–895.

\bibitem{Bodineau19} T. Bodineau, I. Gallagher, L. Saint-Raymond, 
\emph{A Microscopic View of the Fourier Law}, 
C. R. Physique \textbf{20} (2019), 402–418. 
%\href{https://doi.org/10.1016/j.crhy.2019.08.002}{doi:10.1016/j.crhy.2019.08.002}.

\bibitem{TGSS} Thierry Bodineau, Isabelle Gallagher, Laure Saint-Raymond, Sergio Simonella,
\emph{Derivation of the Boltzmann equation from hard-sphere dynamics (after Y. Deng, Z. Hani, and X. Ma)}. Preprint arXiv:2602.04407.

\bibitem{cane21} G. Cane, J.M. Bhat, A. Dhar, C. Bernardin, 
\emph{Localization Effects Due to a Random Magnetic Field on Heat Transport in a Harmonic Chain}, 
J. Stat. Mech. (2021), 113204. 
%\href{https://doi.org/10.1088/1742-5468/ac32b8}{doi:10.1088/1742-5468/ac32b8}.

\bibitem{CL22} R. Castorrini, C. Liverani, 
\emph{Quantitative Statistical Properties of Two-Dimensional Partially Hyperbolic Systems}, 
Adv. Math. \textbf{409} (2022), Paper No.~108625.

\bibitem{CD09} N. Chernov, D. Dolgopyat, 
\emph{Brownian Brownian Motion. I}, 
Mem. Amer. Math. Soc. \textbf{198} (2009), no.~927.

\bibitem{CFKMZ} I. Chevyrev, P.K. Friz, A. Korepanov, I. Melbourne, H. Zhang, 
\emph{Multiscale Systems, Homogenization, and Rough Paths}, 
in \emph{Probability and Analysis in Interacting Physical Systems}, 
Springer Proc. Math. Stat. \textbf{283}, Springer, Cham, 2019, 17–48.

\bibitem{Demasi88} A. De Masi, N. Ianiro, A. Pellegrinotti, E. Presutti, 
\emph{A Survey of the Hydrodynamical Behavior of Many-Particle Systems}, 
NASA STI/Recon Technical Report A \textbf{85} (1984), 123–294.

\bibitem{Demasi89} A. De Masi, R. Esposito, J.L. Lebowitz, 
\emph{Incompressible Navier–Stokes and Euler Limits of the Boltzmann Equation}, 
Comm. Pure Appl. Math. \textbf{42} (1989), 1189–1214.

\bibitem{DKL21} M.F. Demers, N. Kiamari, C. Liverani, 
\emph{Transfer Operators in Hyperbolic Dynamics—An Introduction}, 
33º Colóq. Bras. Mat., IMPA, Rio de Janeiro, 2021.

\bibitem{DHX} Y. Deng, Z. Hani, X. Ma, 
\emph{Long Time Derivation of the Boltzmann Equation from Hard Sphere Dynamics}, to appear in Annals of Mathematics. Preprint
arXiv:2408.07818

\bibitem{DHX2} Y. Deng, Z. Hani, X. Ma, 
\emph{Hilbert's sixth problem: derivation of fluid equations via Boltzmann's kinetic theory}. Preprint 	arXiv:2503.01800.

\bibitem{DSL18} J. De Simoi, C. Liverani, 
\emph{Limit Theorems for Fast–Slow Partially Hyperbolic Systems}, 
Invent. Math. \textbf{213} (2018), 811–1016.

\bibitem{DSL15} J. De Simoi, C. Liverani, 
\emph{The Martingale Approach after Varadhan and Dolgopyat}, 
in \emph{Hyperbolic Dynamics, Fluctuations and Large Deviations}, 
Proc. Sympos. Pure Math. \textbf{89}, Amer. Math. Soc., Providence, RI, 2015, 311–339.

\bibitem{DSL16} J. De Simoi, C. Liverani, 
\emph{Statistical Properties of Mostly Contracting Fast–Slow Partially Hyperbolic Systems}, 
Invent. Math. \textbf{206} (2016), 147–227.

\bibitem{DLPV} J. De Simoi, C. Liverani, C. Poquet, D. Volk, 
\emph{Fast–Slow Partially Hyperbolic Systems versus Freidlin–Wentzell Random Systems}, 
J. Stat. Phys. \textbf{166} (2017), 650–679.

\bibitem{Dobru83} C. Boldrighini, R.L. Dobrushin, Y.M. Suhov, 
\emph{One-Dimensional Hard Rod Caricature of Hydrodynamics}, 
J. Stat. Phys. \textbf{31} (1983), 577–616.

\bibitem{Dobru90}
  C. Boldrighini, R.L. Dobrushin, Y. M. Suhov.
  \emph{One-dimensional hard-rod caricature of hydrodynamics: Navier-Stokes correction},
  Technical report, Dublin Institute for Advances Studies, Preprint (1990).

\bibitem{Bodri97} C. Boldrighini, Y.M. Suhov, 
\emph{One-Dimensional Hard-Rod Caricature of Hydrodynamics: “Navier–Stokes Correction” for Local Equilibrium Initial States}, 
Comm. Math. Phys. \textbf{189} (1997), 577–590.

\bibitem{BS80} L.A. Bunimovich, Ya.G. Sinaĭ, 
\emph{Statistical Properties of Lorentz Gas with Periodic Configuration of Scatterers}, 
Comm. Math. Phys. \textbf{78} (1980/81), 479–497.

\bibitem{Do05} D. Dolgopyat, 
\emph{Averaging and Invariant Measures}, 
Mosc. Math. J. \textbf{5} (2005), 537–576.

\bibitem{DoLi11} D. Dolgopyat, C. Liverani, 
\emph{Energy Transfer in a Fast–Slow Hamiltonian System}, 
Comm. Math. Phys. \textbf{308} (2011), 201–225.

\bibitem{Ei05} A. Einstein, 
\emph{\"Uber die von der molekularkinetischen Theorie der W\"arme geforderte Bewegung von in ruhenden Fl\"ussigkeiten suspendierten Teilchen}, 
Ann. Phys. \textbf{17} (1905), 549–560.

\bibitem{EG} L.C. Evans, R.F. Gariepy, 
\emph{Measure Theory and Fine Properties of Functions}, 
CRC Press, Boca Raton, 1992.

\bibitem{Pablo23} P. Ferrari, S. Olla, 
\emph{Macroscopic Diffusive Fluctuations for Generalized Hard Rod Dynamics}, 
Ann. Appl. Probab. \textbf{35} (2025), 1125–1142. 
%\href{https://doi.org/10.1214/24-AAP2137}{doi:10.1214/24-AAP2137}.

\bibitem{Fritz88} J. Fritz, 
\emph{On the Hydrodynamic Limit of a One-Dimensional Ginzburg–Landau Lattice Model: The a Priori Bounds}, 
J. Stat. Phys. \textbf{47} (1987), 551–572.

\bibitem{FFL94}
  J. Fritz, T. Funaki, and J. L. Lebowitz. \emph{Stationary states of random hamiltonian systems.}
  Probability Theory and Related Fields, \textbf{99}(2) (1994), 211–236.

\bibitem{FLO97} J. Fritz, C. Liverani, S. Olla, 
\emph{Reversibility in Infinite Hamiltonian Systems with Conservative Noise}, 
Comm. Math. Phys. \textbf{189} (1997), 481–496.

\bibitem{giardina05} C. Giardinà, J. Kurchan, 
\emph{The Fourier Law in a Momentum-Conserving Chain}, 
J. Stat. Mech. (2005), P05009. 
%\href{https://doi.org/10.1088/1742-5468/2005/05/P05009}{doi:10.1088/1742-5468/2005/05/P05009}.

\bibitem{GPV} M.Z. Guo, G.C. Papanicolaou, S.R.S. Varadhan, 
\emph{Nonlinear Diffusion Limit for a System with Nearest Neighbor Interactions}, 
Comm. Math. Phys. \textbf{118} (1988), 31–59.

\bibitem{He93} H. Hennion, 
\emph{Sur un Théorème Spectral et son Application aux Noyaux Lipschitziens}, 
Proc. Amer. Math. Soc. \textbf{118} (1993), 627–634.

\bibitem{Hi902} D. Hilbert, 
\emph{Mathematical Problems}, 
Bull. Amer. Math. Soc. \textbf{8} (1902), 437–479.

\bibitem{KM17} D. Kelly, I. Melbourne, 
\emph{Deterministic Homogenization for Fast–Slow Systems with Chaotic Noise}, 
J. Funct. Anal. \textbf{272} (2017), 4063–4102.

\bibitem{Kelley} J.L. Kelley, 
\emph{General Topology}, 
Graduate Texts in Mathematics, Vol.~27, Springer-Verlag, New York, 1975.

\bibitem{KL99} C. Kipnis, C. Landim, 
\emph{Scaling Limits of Interacting Particle Systems}, 
Grundlehren Math. Wiss., Vol.~320, Springer, Berlin, 1999.

\bibitem{klo22} T. Komorowski, J.L. Lebowitz, S. Olla, 
\emph{Heat Flow in a Periodically Forced, Thermostatted Chain II}, 
J. Stat. Phys. \textbf{190} (2023), 87. 
%\href{https://doi.org/10.1007/s10955-023-03103-9}{doi:10.1007/s10955-023-03103-9}.

\bibitem{kos1} T. Komorowski, S. Olla, M. Simon,  \emph{Heat flow in a periodically forced, unpinned thermostatted chain},
    {Electronic Journal of Probability},  \textbf{30} (2025)
    no. 66, 1-48.
%\href{https://doi.org/10.1214/25-EJP1326}{doi:10.1214/25-EJP1326}

\bibitem{La75} O.E. Lanford III, 
\emph{Time Evolution of Large Classical Systems}, 
in \emph{Dynamical Systems, Theory and Applications}, Lecture Notes in Physics, Vol.~38, Springer, Berlin, 1975, 1–111.

\bibitem{Lepri03} S. Lepri, R. Livi, A. Politi, 
\emph{Thermal Conduction in Classical Low-Dimensional Lattices}, 
Phys. Rep. \textbf{377} (2003), 1–80.

\bibitem{Li18} C. Liverani, 
\emph{Transport in Partially Hyperbolic Fast–Slow Systems}, 
Proc. ICM 2018, Vol.~III, World Scientific (2018) 2643–2667.

\bibitem{LO96} C. Liverani, S. Olla, 
\emph{Ergodicity in Infinite Hamiltonian Systems with Conservative Noise}, 
Probab. Theory Relat. Fields \textbf{106} (1996), 401–445.

\bibitem{LO12} C. Liverani, S. Olla, 
\emph{Toward the Fourier Law for a Weakly Interacting Anharmonic Crystal}, 
J. Amer. Math. Soc. \textbf{25} (2012), 555–583.

\bibitem{Morrey55} C.B. Morrey, 
\emph{On the Derivation of the Equations of Hydrodynamics from Statistical Mechanics}, 
Comm. Pure Appl. Math. \textbf{8} (1955), 279–326.

\bibitem{OVY} S. Olla, S.R.S. Varadhan, H.T. Yau, 
\emph{Hydrodynamical Limit for a Hamiltonian System with Weak Noise}, 
Comm. Math. Phys. \textbf{155} (1993), 523–560.

\bibitem{sasada13} S. Olla, M. Sasada, 
\emph{Energy Diffusion in an Anharmonic Chain with Conservative Noise}, 
Probab. Theory Relat. Fields \textbf{157} (2013), 721–775.

\bibitem{PS83} Ya.B. Pesin, Ya.G. Sinaĭ, 
\emph{Gibbs Measures for Partially Hyperbolic Attractors}, 
Ergodic Theory Dynam. \textbf{2} (1983), 417–438.

\bibitem{Laure09} L. Saint-Raymond, 
\emph{Hydrodynamic Limits of the Boltzmann Equation}, 
Lecture Notes in Mathematics, Vol.~1971, Springer, 2009.

\bibitem{saito18} K. Saito, M. Sasada, 
\emph{Thermal Conductivity for Coupled Charged Harmonic Oscillators with Noise in a Magnetic Field}, 
Comm. Math. Phys. \textbf{361} (2018), 951–995.

\bibitem{Sp} H. Spohn, 
\emph{Large Scale Dynamics of Interacting Particles}, 
Texts and Monographs in Physics, Springer, 1991.

\bibitem{Spohn23} H. Spohn, 
\emph{Hydrodynamic Scales of Integrable Many-Particle Systems},
World Scientific, 2024.
%\href{https://doi.org/10.1142/13600}{doi:10.1142/13600} 
%arXiv:2301.08504 (2023).

\bibitem{Va} S.R.S. Varadhan, 
\emph{Nonlinear Diffusion Limit for a System with Nearest Neighbor Interactions II}, 
in \emph{Asymptotic Problems in Probability Theory: Stochastic Models and Diffusions on Fractals}, 
Pitman Res. Notes Math. Ser., Vol.~283, Longman Sci, 1993.
\end{thebibliography}
\end{document}